\documentclass[10pt,reqno]{amsart}
\usepackage{xcolor}
\usepackage[
    colorlinks=true,       
    citecolor=magenta,        
    linkcolor=blue,         
    urlcolor=teal          
]{hyperref}

\usepackage[utf8]{inputenc}
\usepackage{setspace}
\usepackage[cmtip,arrow]{xy}
\usepackage{lscape}
\usepackage{latexsym}
\usepackage[activeacute,english]{babel}
\usepackage{graphicx}
\usepackage{amsmath}
\usepackage{adjustbox} 
\usepackage{amsfonts}
\usepackage{pst-node}
\usepackage{amsthm}
\usepackage{comment}
\usepackage{amsmath}
\usepackage{enumitem}
\usepackage{extpfeil}
\usepackage{extarrows} 
\usepackage{amssymb}
\usepackage{amscd}
\usepackage{mathrsfs}
\usepackage{tikz}
\usepackage{pb-diagram,pb-xy}
\xyoption{all}
\usepackage{graphics}
\usepackage{stmaryrd}
\usepackage{yfonts,bbm}
\usepackage{tikz-cd}
\usetikzlibrary{decorations.pathmorphing}
\allowdisplaybreaks 

\newtheorem{theorem}{Theorem}[section]
\newtheorem{lemma}[theorem]{Lemma}
\newtheorem{proposition}[theorem]{Proposition}
\newtheorem{corollary}[theorem]{Corollary}
\newtheorem{definition}[theorem]{Definition}
\newtheorem{example}[theorem]{Example}

\newtheorem{Remark}[theorem]{Remark}

\usepackage[square,sort,comma,numbers]{natbib}
\usepackage{url}
\usepackage{relsize}

\newcommand{\A}{\mathcal{A}}
\newcommand{\B}{\mathcal{B}}
\newcommand{\C}{\mathcal{C}}
\newcommand{\D}{\mathcal{D}}
\newcommand{\E}{\mathcal{E}}
\newcommand{\F}{\mathcal{F}}
\newcommand{\Ll}{\mathcal{L}}
\newcommand{\N}{\mathcal{N}}
\newcommand{\R}{\mathcal{R}}
\newcommand{\T}{\mathcal{T}}

\newcommand{\X}{\mathcal{X}}
\newcommand{\Y}{\mathcal{Y}}

\newcommand{\tf}{\mathfrak{t}}
\newcommand{\lf}{\mathfrak{l}}

\newcommand{\Ee}{\mathbb{E}}
\newcommand{\Ff}{\mathbb{F}}
\newcommand{\Gg}{\mathbb{G}}

\newcommand{\s}{\mathfrak{s}}
\newcommand{\rr}{\mathfrak{r}}
\newcommand{\uu}{\mathfrak{u}}

\newcommand{\xdashrightarrow}[1]{\overset{#1}{\dashrightarrow}} 
\newcommand{\di}{\backslash}
\newcommand{\Ho}{\operatorname{Hom}}
\newcommand{\Hom}{\operatorname{Hom}}
\newcommand{\Mor}{\operatorname{Mor}}
\newcommand{\Ker}{\operatorname{Ker}}

\newcommand{\im}{\operatorname{Im}}

\newcommand{\Ebb}{\mathbb{E}}        
\newcommand{\Fbb}{\mathbb{F}}        

\newcommand{\wC}{\wt{\C}}             
\newcommand{\wE}{\wt{\Ebb}}           
\newcommand{\ws}{\wt{\sfr}}           

\newcommand{\op}{\mathrm{op}}       

\newcommand{\id}{\mathrm{id}}         

\newcommand{\co}{\colon}              

\newcommand{\lla}{\longleftarrow}     
\newcommand{\lra}{\longrightarrow}    
\newcommand{\dra}{\dashrightarrow}  

\newcommand{\wt}{\widetilde}          
\newcommand{\FR}[3]{{#1}\backslash {#2}/{#3}}

\newcommand{\sfr}{\mathfrak{s}}  
\newcommand{\tfr}{\mathfrak{t}}  
\newcommand{\Ssc}{\mathscr{S}}   

\newcommand{\Acal}{\mathcal{A}} 
\newcommand{\Bcal}{\mathcal{B}} 
\newcommand{\Ccal}{\mathcal{C}} 
\newcommand{\Dcal}{\mathcal{D}} 
\newcommand{\Mcal}{\mathcal{M}} 
\newcommand{\Xcal}{\mathcal{X}} 
\newcommand{\Ycal}{\mathcal{Y}} 
\newcommand{\ovl}{\overline}    
\newcommand{\ov}{\overset}      

\newcommand{\bsm}{\begin{smallmatrix}}
	\newcommand{\esm}{\end{smallmatrix}}

\newcommand{\lam}{\lambda}       
\newcommand{\del}{\delta}        

\begin{document}
	\title{$s$-recollements  in extriangulated categories}
	
	\author[Mendoza-Santiago-Velasco]{  Octavio Mendoza \hspace{1cm} Valente Santiago \hspace{1cm} Edgar Omar Velasco  }

	\thanks{The authors would like to thank the SECIHTI Project CBF-2023-2024-2630 and the Research Support Office of the Metropolitan Autonomous University (DAI-UAM) for the support granted. The third author also acknowledges the financial support provided by the Secretaría de Ciencia, Humanidades, Tecnología e Innovación (SECIHTI) through the 2025 Initial Postdoctoral Fellowship Program.}
	\subjclass{2000]{Primary 18A25, 18E05, 18E10; Secondary 16D90,16G10}}
	
	\keywords{Functor categories, Recollements, Extriangulated  categories}
	\dedicatory{}
	\maketitle
	
	\begin{abstract}
	In this paper, we investigate recollements in the setting of extriangulated categories which simultaneously generalize exact and triangulated categories. Building on the notion of recollement given in the extriangulated setting by \cite{AMP25}, we introduce and develop the concept of s-recollement and establish a completion theorem for fully faithful functors between extriangulated categories. More precisely, we prove that a fully faithful functor admitting both a left and a right adjoint, with all three functors involved being extriangulated, can be completed to an s-recollement. This result extends the classical theorem of Cline-Parshall-Scott from triangulated categories to the broader framework of extriangulated categories. Moreover, our detailed proof highlights the role of the extriangulatedness of the adjoint functors and explains why the choice of the underlying notion of s-recollement is essential for this extension. Our approach also draws on \(s\)-torsion pairs and localization theory in extriangulated categories, providing the categorical tools needed to establish the main result.
	\end{abstract}
	
	\section{Introduction}
	The notion of \emph{recollement} was introduced for triangulated categories by Beilinson-Bernstein-Deligne in \cite{Beilinson-Bernstein-Deligne}. This was done as a part of their axiomatization of the six Grothendieck functors for
derived categories of sheaves, within the framework of the study of
stratifications of spaces. On the other hand, one of the first fundamental examples of \emph{recollement} in the context of abelian categories appeared in the construction of perverse sheaves, developed by MacPherson and Vilon in \cite{MV}.
	
Towards the end of the 1980's, in the context of representation theory, Cline-Parshall-Scott highlighted the close relationship between recollements of derived categories, highest weight categories and  quasi-hereditary algebras \cite{Cline-Parshall-Scott, Parshall-Scott}. Subsequently, Kuhn employed recollements of certain abelian categories in the study of the generic representation theory of the general linear group  \cite{Kuhn}. We also point out that several authors, among them Happel, Beligiannis and Reiten, as well as Angeleri Hügel, Koenig, and Liu, have explored in depth the connections between recollements of triangulated categories and various aspects of representation theory and homological algebra. In particular, their work addresses the relationship that these structures have with tilting theory, homological conjectures and the stratification of derived categories of rings (see \cite{AKL1}, \cite{AKL2}, \cite{BR} and \cite{Happel}).
	
More recently, with the development of extriangulated categories \cite{NP19}, which serve as a common generalization of both exact and triangulated categories, there has been growing interest in extending the notion of recollement to this new setting. The first definition of recollement for extriangulated categories was introduced by 
Wang-Wei-Zhang in \cite{Wang-Wei-Zhang}. Among other results, their work recovers a theorem by Chen concerning recollements of triangulated categories and highlights new phenomena that emerge when considering abelian categories.	
	
It is worth noting that in \cite[Def. 3.1]{Wang-Wei-Zhang}, the existence of certain conflations is required, as stipulated in conditions (R4) and (R5). This requirement is natural when considering the definition of recollement in triangulated categories \cite[Def. 1.4.3]{Beilinson-Bernstein-Deligne}. However, in the triangulated setting, it is possible to formulate an equivalent definition without explicitly requiring the existence of such conflations, as explained in \cite[p.~96]{Psaroudakis}. Furthermore, in the context of abelian categories there is a bijection  between recollements of module categories and certain torsion pairs \cite[Cor. 6.6]{Psaroudakis} called ttf-triples. Inspired by this correspondence, a second notion of recollement for extriangulated categories was introduced in \cite[Def. 5.3]{AMP25}. This is precisely the definition we adopt in the present work to formulate the notion of s-recollement as can be seen in Definition \ref{def: recollement alejandro 1}.
	
The main objective of this work is to provide a detailed proof that a fully faithful functor between extriangulated categories, which admits both left and right adjoints and for which both of these adjoints are themselves extriangulated functors, can be extended to a s-recollement. This result constitutes a generalization, to the extriangulated context, of the result obtained by Cline-Parshall-Scott in \cite{Cline-Parshall-Scott} (see Theorem 2.1) for triangulated categories.
	
Now, we say how this paper is organized. In Section 2, we introduce preliminaries on extriangulated categories, as well as the basic concepts of extriangulated functors, extriangulated natural transformations and extriangulated adjoint pairs. We also recall extriangulated categories with negative first extension group. This is needed to define   $s$-torsion pairs and the induced $s$-torsion and $s$-torsion-free functors. The main result of this section is Proposition \ref{pro:torsion_funtor_extriangulado}.
	
Section 3 addresses the localization theory in extriangulated categories developed in \cite{NakaokaOgawaSakai2022}. The primary purpose of this section is to establish the notation and preliminaries which are necessary for the development of the subsequent sections of the paper. It is worth noting that special emphasis is placed on describing the elements of the abelian group $\widetilde{\Ee}(C,A)$ in terms of the so-called \emph{coends}, an aspect not explicitly developed in \cite{NakaokaOgawaSakai2022}. This description is very important for understanding the origin of the slight abuse of notation employed when representing elements $[\FR{\overline{t}}{\overline{\delta}}{\overline{s}}]\in \widetilde{\Ee}(C,A)$. Finally, the construction of the category $(\widetilde{\C},\widetilde{\Ee},\widetilde{\s})$ is recalled, along with its universal property.	
	
In Section 4, we present a brief review of the notion of recollement in additive categories, as well as some of its basic properties. Furthermore, we introduce the definition of $s$-recollement (see Definition \ref{def: recollement alejandro 1}) and establish several preliminary results that serve as fundamental components for the proof of our main result.	
	
Finally, in Section 5, we present our main result Theorem \ref{teo: completar recoll}, along with a fully detailed proof. This proof sheds light on why the definition of  recollement which is introduced in \cite{Wang-Wei-Zhang} is inadequate for extending the result of Cline-Parshall–Scott to the context of extriangulated categories. In particular, it highlights the necessity of requiring the functors of the recollement $(i^*, i_*, i^!)$ to be an extriangulated adjoint triple.

\section{Preliminaries}

	Exact categories and triangulated categories are valuable tools in contemporary mathematics. Their usefulness lies in the fact that they provide suitable settings for homological algebra. More recently, H. Nakaoka and Y. Palu introduced in \cite{NP19} the notion of an extriangulated category, a concept that encompasses both triangulated and exact categories.
	
	The precise definition of an extriangulated category can be found in \cite[Def.~2.12]{NP19}. For an intuitive perspective, extriangulated categories can be understood as subcategories of triangulated categories that are closed under extensions.
	
	An extriangulated category consists of a triple $(\mathcal{C}, \mathbb{E}, \mathfrak{s})$, where $\mathcal{C}$ is an additive category, $\mathbb{E} : \mathcal{C}^{op} \times \mathcal{C} \to \mathbf{Ab}$ is an additive bifunctor, and $\mathfrak{s}$ is a function that associates to each element $\delta \in \mathbb{E}(C,A)$ an equivalence class $[A \xrightarrow{x} B \xrightarrow{y} C]$ satisfying a collection of axioms  \cite[Def.~2.12]{NP19}.
For example, let $\mathcal{C}$ be an abelian category and consider $\mathbb{E} := \operatorname{Ext}^1_{\mathcal{C}}(-,-)$, or let $\mathcal{C}$ be a triangulated category with $\mathbb{E} = \operatorname{Hom}_{\mathcal{C}}(-,-[1])$. Moreover, whenever we refer to an abelian or triangulated category as an extriangulated category, it will always be understood via an $\mathbb{E}$-bifunctor of one of these forms.
	
	Let $(\mathcal{C}, \mathbb{E}, \mathfrak{s})$ be an extriangulated category, as in \cite{NP19}, we recall the following notions. For $A,B \in \mathcal{C},$ an element $\delta \in \mathbb{E}(C,A)$ is called \textbf{$\mathbb{E}$-extension}, and the zero element $0 \in \mathbb{E}(B,A)$ is called the \textbf{split $\mathbb{E}$-extension}.
	If $\mathfrak{s}(\delta) = [A \xrightarrow{x} B \xrightarrow{y} C]$, we say that $A \xrightarrow{x} B \xrightarrow{y} C$ \textbf{realizes} $\delta$, or that it is an \textbf{$\mathfrak{s}$-conflation} (or simply a \textbf{conflation}), also known as an \textbf{$\mathbb{E}$-triangle}; and in this case, we use the notation $\delta : A \xrightarrow{x} B \xrightarrow{y} C$ or $A \xrightarrow{x} B \xrightarrow{y} C \overset{\delta}{\dashrightarrow}$.
	
	For each $\delta \in \mathbb{E}(C,A)$, $a \in \operatorname{Hom}_{\mathcal{C}}(A,A')$, and $c \in \operatorname{Hom}_{\mathcal{C}}(C',C)$, we can define
	$a \cdot \delta := \mathbb{E}(C,a)(\delta)\;$ and $\;\delta \cdot c := \mathbb{E}(c,A)(\delta)$. Note that each $\delta \in \mathbb{E}(C,A)$ induces morphisms
	\[
	\delta \cdot - : \operatorname{Hom}_{\mathcal{C}}(W,C) \to \mathbb{E}(W,A),
	\qquad
	- \cdot \delta : \operatorname{Hom}_{\mathcal{C}}(A,W) \to \mathbb{E}(C,W),
	\]
	for every $W \in \mathcal{C}$.
	
	Let $\delta \in \mathbb{E}(C,A)$ and $\delta' \in \mathbb{E}(C',A')$. A \textbf{morphism of extensions} $\delta \to \delta'$ is a pair
	$(a,c) \in \operatorname{Hom}_{\mathcal{C}}(A,A')
	\times \operatorname{Hom}_{\mathcal{C}}(C,C')$ such that
	$a \cdot \delta = \delta' \cdot c$.\\
	By the definition of an extriangulated category, for any morphism of $\mathbb{E}$-extensions $(a,c) : \delta \to \delta'$,
	with $\mathfrak{s}(\delta) = [A \xrightarrow{x} B \xrightarrow{y} C]$  and $\mathfrak{s}(\delta') = [A' \xrightarrow{x'} B' \xrightarrow{y'} C']$, there exists a morphism $b : B \to B'$ such that $b x = x' a$ and
	$y' b = c y$. In this case, it is said that $(a,b,c)$ \textbf{realizes} the morphism $(a,c) : \delta \to \delta'.$
	
	We begin by recalling some notions concerning closure properties of
	subcategories of an extriangulated category that will be used throughout
	the paper. In particular, we will consider the closure under extensions,
	cones and cocones, as well as subcategories that are compatible with the
	extriangulated structure. We recall that any class $\Xcal$ in a category $\Ccal,$ which is denoted by $\Xcal\subseteq\Ccal,$ will be considered as a full subcategory of $\Ccal.$

	Let $\mathcal{C}$ be an extriangulated category and $\mathcal{X},$ $\mathcal{Y}\subseteq\mathcal{C}.$ By following \cite{Tsukamoto}, we denote by $\mathcal{X} * \mathcal{Y}$ the
			class of the objects $M\in\mathcal{C}$  for which there exists a conflation $X \to M \to Y $
			with $X \in \mathcal{X}$ and $Y \in \mathcal{Y}.$ The  class $\mathcal{X}$ is 
			\textbf{closed under extensions} if $\mathcal{X} * \mathcal{X}\subseteq \mathcal{X}.$ It is said that $\Xcal$ is \textbf{closed under cones} if for any conflation
			$A \to B \to  C $ with $A,B \in \mathcal{X}$, we have that $C \in \mathcal{X}$.
			Dually, $\mathcal{X}$ is \textbf{closed under cocones} if 
			for any conflation $A \to B \to C $ with $B,C \in \mathcal{X}$, we have that $A \in \mathcal{X}.$ We say that $\Xcal$ is \textbf{left thick} if it is closed under extensions, cocones and direct summands in $\Ccal.$ Dually, we say that $\Xcal$ is \textbf{right thick} if it is closed under extensions, cones and direct summands in $\Ccal.$ Finally, we say that $\Xcal$ is \textbf{thick} if it is left and right thick.

	We next recall the notion of extriangulated subcategory. Let $(\mathcal C,\mathbb E,\mathfrak s)$ be an extriangulated category. A class 
		$\mathcal D \subseteq \mathcal C$ is called an \textbf{extriangulated subcategory} of $\C$ if it is closed under extensions and $0\in\Dcal.$ Notice that 
		in this case $\Dcal$ is closed under isomorphisms and finite coproducts in $\Ccal.$ Moreover,  the extriangulated structure of $\mathcal{C}$ restricts
	naturally to $\Dcal$ and hence $(\Dcal,\mathbb E|_{\Dcal},\mathfrak s|_{\Dcal})$ becomes an extriangulated category.
	
	The following long exact sequences (see \cite[Cor.~3.12]{NP19}),  which are obtained from a conflation, will be very useful throughout the paper.
	
	\begin{proposition}\label{PsucexExt}
		Let $\delta : A \xrightarrow{u} B \xrightarrow{v} C$ be a conflation in an extriangulated category $\mathcal{C}$ and let $W \in \mathcal{C}$. Then the following two sequences are exact:
		\[
		\resizebox{1.0\textwidth}{!}{$
			(W,A) \xrightarrow{(W,u)}
			(W,B) \xrightarrow{(W,v)}
			(W,C) \xrightarrow{\delta \cdot -}
			\mathbb{E}(W,A)
			\xrightarrow{\mathbb{E}(W,u)}
			\mathbb{E}(W,B) \xrightarrow{\mathbb{E}(W,v)}
			\mathbb{E}(W,C)
			$}
		\]
		and
		\[
		\resizebox{1.0\textwidth}{!}{$
			(C,W) \xrightarrow{(v,W)}
			(B,W) \xrightarrow{(u,W)}
			(A,W) \xrightarrow{- \cdot \delta}
			\mathbb{E}(C,W)
			\xrightarrow{\mathbb{E}(v,W)}
			\mathbb{E}(B,W) \xrightarrow{\mathbb{E}(u,W)}
			\mathbb{E}(A,W).
			$}
		\]
		Here we set $(W,-):=\mathcal{C}(W,-)$ and $(-,W):=\mathcal{C}(-,W)$.
	\end{proposition}
	
	\subsection*{Extriangulated functors} Following \cite[Def.~2.32]{BTS} and 
	\cite[Def.~2.11]{NakaokaOgawaSakai2022}, we recall that
	an \textbf{extriangulated functor} $(F,\beta):
	(\mathcal{A},\mathbb{E},\mathfrak{s})
	\to
	(\mathcal{B},\mathbb{G},\mathfrak{t})$  between extriangulated categories consists of an additive functor
	$F:\mathcal{A}\to\mathcal{B}$ and a natural transformation $\beta:\mathbb{E}(-,?)
	\to
	\mathbb{G}(F(-),F(?))$ such that, for every $\mathbb{E}$-triangle
	$\delta:\;X\xrightarrow{f}W\xrightarrow{g}Y$ in $\mathcal{A}$ we have that $\beta_{Y,X}(\delta):\;F(X)
	\xrightarrow{F(f)}
	F(W)
	\xrightarrow{F(g)}
	F(Y)
	$ is a $\mathbb{G}$-triangle in $\B.$ In such a case, for simplicity, we just say that the functor $F:\A\to \B$ is extriangulated.

	The composition of extriangulated functors $(F,\beta):
	(\mathcal{A},\mathbb{E},\mathfrak{s})
	\to
	(\mathcal{B},\mathbb{G},\mathfrak{t})$ and $(H,\gamma):
	(\mathcal{B},\mathbb{G},\mathfrak{t})
	\to
	(\mathcal{C},\mathbb{T},\mathfrak{u})$ is  the extriangulated functor $(H\circ F,\alpha:=\gamma\circ \beta)$, where
	$\alpha_{Y,X}
	:=
	\gamma_{F(Y),F(X)}
	\circ
	\beta_{Y,X}$
	for all $X,Y\in\mathcal{A},$ see
	\cite[Def.~3.18(ii)]{BTS}.
	
	Let $(F,\beta),(G,\gamma):
	(\mathcal{A},\mathbb{E},\mathfrak{s})
	\to
	(\mathcal{B},\mathbb{T},\mathfrak{t})$ be extriangulated functors. An \textbf{extriangulated natural transformation}
	$\eta:(F,\beta)
	\to
	(G,\gamma)$ is a natural transformation
	$\eta:F\to G$
	of additive functors such that, for all $X,Y\in\mathcal{A}$, the following diagram commutes
	\[
	\xymatrix{
		\mathbb{E}(Y,X)
		\ar[rr]^{\beta_{Y,X}}
		\ar[d]_{\gamma_{Y,X}}
		& &
		\mathbb{T}(F(Y),F(X))
		\ar[d]^{\mathbb{T}(FY,\eta_X)}
		\\
		\mathbb{T}(G(Y),G(X))
		\ar[rr]_{\mathbb{T}(\eta_Y,GX)}
		& & 
		\mathbb{T}(F(Y),G(X)).
	}
	\]
	For further details,  see
	\cite[Def.~2.11]{NakaokaOgawaSakai2022}
	and 
	\cite[Def.~4.1]{BTS}.
	
	\begin{definition}\label{extadjpD} An \textbf{extriangulated adjoint pair} is a pair $(F,G)$ of adjoint funtors $(F:\Xcal\to \Ycal,G:\Ycal \to \Xcal)$ such that: (1) $F$ and $G$ are extriangulated funtors and (2) the unit $\eta:1_\Xcal\to G\circ F$ and the co-unit $\varepsilon:F\circ G\to 1_\Ycal$ are extriangulated natural transformations.
	\end{definition}

\begin{example} \label{Eextf=tf} Any triangulated functor between triangulated categories can be seen as an extriangulated functor. Conversely,  any extriangulated functor between triangulated categories can be seen as a triangulated functor, see \cite[Ex. 3.4] {BGLS}. Moreover, in the case of an adjoint pair 
$(F:\Xcal\to \Ycal,G:\Ycal \to \Xcal)$ where $F$ and $G$ are triangulated functors, requiring that the unit and counit are extriangulated is superfluous, see 
\cite[page 375] {BGLS}.
	\end{example}	
	
\begin{definition}\label{extadjtD} An \textbf{extriangulated adjoint triple} $(F,G,H)$ is  given by the following diagram 
\[\xymatrix{
			\mathcal{X} \ar[r]^{G} &
			\mathcal{Y} 
			\ar@/^1pc/[l]^{H}
			\ar@/_1pc/[l]_{F} &			
		}
		\]
where $(F,G)$ and $(G,H)$ are extriangulated adjoint pairs. 
	\end{definition}

	\begin{Remark}\label{Rthick-sub}
		The following properties are straightforward consequences  of the above definitions.
		
		\begin{itemize}
			\item[(1)] Every left (right) thick subcategory  of an extriangulated category   is an
			extriangulated subcategory.
			
			\item[(2)] Let $(F,\varphi):
			(\mathcal{C},\mathbb{E},\mathfrak{s})
			\to
			(\mathcal{C}',\mathbb{G},\mathfrak{t})$  be an extriangulated functor between extriangulated categories.
			Then $\operatorname{Ker}(F):=\{X\in\mathcal{C}\mid F(X)=0\}$ is a thick subcategory of $\mathcal{C}.$ In particular $\Dcal:=\Ker(F)$ is an extriangulated subcategory of 
			$ \mathcal{C}.$ Moreover, we have that $(i_\D,\varphi|_{\Dcal\times\Dcal}):
			(\Dcal,\mathbb{E}|_{\Dcal},\mathfrak{s}|_{\Dcal})
			\to (\mathcal{C},\mathbb{E},\mathfrak{s})$
			is  an extriangulated functor, where $i_\D:\D\to \C$ is the inclusion functor.
		\item[(3)] 	Let $(\C,\mathbb{E},\mathfrak{s})$ be an extriangulated category and  $\X\subseteq\C$ be an extriangulated subcategory. Then $(i_\X,1):(\X,\mathbb{E}|_\X,\mathfrak{s}|_\X)\to (\C,\mathbb{E},\mathfrak{s})$ is an extriangulated functor, where $i_\X:\X\to \C$ is the inclusion funtor and $1_{A,B}:\mathbb{E}|_\X(A,B)\to \mathbb{E}(A,B)$ is the identity morphism.
		\end{itemize}
	\end{Remark}
	
	\subsection*{Extriangulated categories with negative first extension and
		s-torsion pairs}
	
	It is well known that the exact sequences in
	Proposition~\ref{PsucexExt} can be extended to the
	left in a variety of extriangulated categories. For instance, in the case
	of abelian categories, the resulting sequences extend to the left by zero,
	whereas in the case of triangulated categories, they extend further by
	terms involving the bifunctor $\mathcal{C}(-,\Sigma^{-1}?)$ (see Example \ref{triang-ext-1Ex}). Extriangulated
	categories for which such an extension exists are said to admit a
	\textbf{negative first extension}.
	
	We recall this notion below, see \cite[Def.~2.3]{Tsukamoto}. For further details on the existence of
	negative extensions in extriangulated categories, we refer the reader to
	\cite{Gorsky}.
	
	\begin{definition}\label{Dext-1}
		Let $\mathcal{C} = (\mathcal{C}, \mathbb{E}, \mathfrak{s})$ be an
		extriangulated category. We say that $\mathcal{C}$ admits a
		\textbf{negative first extension} if the following two conditions hold true.
		\begin{itemize}
		\item[(a)] There exists an additive bifunctor 
		$\mathbb{E}^{-1} :
		\mathcal{C}^{\mathrm{op}} \times \mathcal{C}
		\to
		\mathbf{Ab}.$
		\item[(b)] For each $\delta \in \mathbb{E}(C,A)$ there are natural transformations 
	\begin{center}
	$\delta^{-1}_{\sharp} :
		\mathbb{E}^{-1}(-,C)
		\to
		\mathcal{C}(-,A)\;$ and $\;\delta^{-1}_{\flat} :
		\mathbb{E}^{-1}(A,-)
		\to
		\mathcal{C}(C,-)$
\end{center}		
		 such that, for any conflation $\delta :
		A \xrightarrow{u} B \xrightarrow{v} C$
		and any $W \in \mathcal{C}$, the following sequences are exact
		\[
		\scalebox{0.95}{$
			\mathbb{E}^{-1}(W,A)
			\xrightarrow{\mathbb{E}^{-1}(W,u)}
			\mathbb{E}^{-1}(W,B)
			\xrightarrow{\mathbb{E}^{-1}(W,v)}
			\mathbb{E}^{-1}(W,C)
			\xrightarrow{\delta^{-1}_{\sharp}}
			\mathcal{C}(W,A)
			\xrightarrow{\mathcal{C}(W,u)}
			\mathcal{C}(W,B),
			$}
		\]
		
		\[
		\scalebox{0.95}{$
			\mathbb{E}^{-1}(C,W)
			\xrightarrow{\mathbb{E}^{-1}(v,W)}
			\mathbb{E}^{-1}(B,W)
			\xrightarrow{\mathbb{E}^{-1}(u,W)}
			\mathbb{E}^{-1}(A,W)
			\xrightarrow{\delta^{-1}_{\flat}}
			\mathcal{C}(C,W)
			\xrightarrow{\mathcal{C}(v,W)}
			\mathcal{C}(B,W).
			$}
		\]
	 \end{itemize}	
		In case (a) and (b) hold true, we say that $\mathcal{C}
		= (\mathcal{C},\mathbb{E},\mathfrak{s},\mathbb{E}^{-1})$ is an \textbf{extriangulated category with negative first extension}.
	\end{definition}

	\begin{example}\label{triang-ext-1Ex} Any triangulated category $\mathcal{D}$ with shift funtor $\Sigma:\Dcal\to \Dcal$ becomes an extriangulated category with negative first extension. Indeed, the necessary data to do so can be given as follows.\begin{enumerate}
\item We set  $\Ee(-,?):=\mathcal{D}(C,\Sigma ?)$ and $\Ee^{-1}(-,?):=\mathcal{C}(-,\Sigma^{-1}?).$			 
\item Let $\delta\in\Ee(C,A)$ and $A\xrightarrow{f}B\xrightarrow{g}C\xrightarrow{\delta}\Sigma A$ be a distinguished triangle. We set
				$\s(\delta):=[A\xrightarrow{f}B\xrightarrow{g}C]$ and also define 
				the natural transformations  
\begin{center}
$\delta^{-1}\colon\Ee^{-1}(-,C)\to \Ee(-,A)\;\;$ and
				$\;\;
				\delta^{\ast-1}\colon\Ee^{-1}(A,-)\to \Ee(C,-)
				$
\end{center}
			as follows. For each $W\in\mathcal{D}$, we set
				\[
				(\delta^{-1})_W\colon
				\Ee^{-1}(W,C)=\mathcal{D}(W,\Sigma^{-1}C)
				\xrightarrow{\mathcal{D}(W,\Sigma^{-1}\delta)}
				\mathcal{D}(W,A)=\Ee(W,A),
				\]
				\[
				(\delta^{\ast-1})_W\colon
				\Ee^{-1}(A,W)=\mathcal{D}(A,\Sigma^{-1}W)
				\xrightarrow{\ \sim\ }
				\mathcal{D}(\Sigma A,W)
				\xrightarrow{\mathcal{D}(\delta,W)}
				\mathcal{D}(C,W)=\Ee(C,W).
				\]
			\end{enumerate}
	\end{example}
	
	Extriangulated categories also allow for analogues of torsion pairs in
	abelian categories. 
	
	\begin{definition}\cite[Def.~3.1]{HeHuZhou}\label{def: torsion}
		Let $\mathcal{C}$ be an extriangulated category. A
		\textbf{torsion pair} in $\mathcal{C}$ is a pair
		$\mathbf{x}=(\mathcal{T},\mathcal{F})$ of classes in
		$\mathcal{C}$ which are closed under isomorphisms in $\mathcal{C}$ and such that
		$\mathcal{C}(\mathcal{T},\mathcal{F})=0$ and $\mathcal{C}=\mathcal{T}*\mathcal{F}.$
	\end{definition}
	
	Notice that the above definition is just the literal transcription of the axioms of a torsion pair in abelian categories. However, in order to get nice properties the above definition needs to be strengthened as was done in \cite[Def.~3.1]{Tsukamoto}.
	
	\begin{definition}\label{defi: s torsion}
		Let $\mathcal{C}$ be an extriangulated category with negative first
		extension. An \textbf{s-torsion pair} in $\mathcal{C}$ is
		a torsion pair
		$\mathbf{x}=(\mathcal{T},\mathcal{F})$ such that $\mathbb{E}^{-1}(\mathcal{T},\mathcal{F})=0.$
	In this case, $\mathcal{T}$ is called 
	\textbf{s-torsion class}, while $\mathcal{F}$ is called 
	\textbf{s-torsion-free class}. 
\end{definition}

We point out that the term ``s'' in the s-torsion pair stands for ``shift-closed'' and it comes from the description given in \cite[Lem. 3.3]{Tsukamoto} of the s-torsion pairs in triangulated categories.
	
	In order to state one of the very nice properties that s-torsion pairs enjoy, we need to introduce the orthogonal complements of a given class. Let $\Ccal$ be an additive category and $\Xcal\subseteq\C.$ The \textbf{left orthogonal complement} ${}^\perp\X$ of $\X$ is given by ${}^\perp\X:=\{C\in\C\;:\; \Hom_\C(C,-)|_\X=0\}.$ Dually, we have the \textbf{right orthogonal complement} $\X^\perp$ of $\X.$

	\begin{proposition} \cite[Prop.~3.2]{Tsukamoto} \label{P-ortTF}
		Let $(\mathcal{T},\mathcal{F})$ be an s-torsion pair in an extriangulated category
		$\mathcal{C}$ with negative first extension. Then $\mathcal{T}^{\perp}=\mathcal{F}$ and  ${}^{\perp}\mathcal{F}=\mathcal{T}.$ In particular, we have that 
		$\mathcal{T}$ is right thick and $\mathcal{F}$ is left thick.
	\end{proposition}

	\begin{Remark}\label{Rk-inc-ext-tp}
		Let $(\mathcal{T},\mathcal{F})$ be an s-torsion pair in the extriangulated category $(\mathcal{C},\mathbb{E},\mathfrak{s})$ with negative first extension. 
\begin{enumerate}
\item By Proposition~\ref{P-ortTF} and Remark \ref{Rthick-sub}, we have that $\mathcal{T}$ and $\mathcal{F}$ are extriangulated subcategories of $\C.$ Moreover the functors 
\begin{center}
$(i_\T,1):(\T, \mathbb{E}|_\T,\mathfrak{s}|_\T)\to (\mathcal{C},\mathbb{E},\mathfrak{s})$ and $(i_\F,1):(\F, \mathbb{E}|_\F,\mathfrak{s}|_\F)\to (\mathcal{C},\mathbb{E},\mathfrak{s})$
\end{center}
are extriangulated.

\item For each $C\in\mathcal{C},$ there exists (see \cite[Prop.~3.7]{Tsukamoto}) a canonical conflation $\delta_C:\;T_C\xrightarrow{\iota_C} C \xrightarrow{\pi_C} F^C$ with $T_C\in\mathcal{T}$ and $F^C\in\mathcal{F}$, which is unique up to isomorphism of conflations.
Moreover, for each $h:X\to Y$ in $\Ccal$ there exist unique $T_h:T_X\to T_Y$ and $F^h:F^X\to F^Y$ such the following diagram commutes
\[\xymatrix{
T_X \ar[d]_{T_h} \ar[r]^{\iota_X} & X \ar[d]^h \ar[r]^{\pi_X} & F^X \ar[d]^{F^h}\\
T_Y \ar[r]^{\iota_Y} & Y  \ar[r]^{\pi_Y} & F^Y.
}\]
Thus, the assignments $C\mapsto T_C$ and $C\mapsto F^C$ induce additive functors: 
the torsion functor $\mathfrak{t}:\mathcal{C}\to \mathcal{T}$ and the torsion-free functor
$\mathfrak{l}:\mathcal{C}\to \mathcal{F}$, which are, respectively, right and left adjoints to the inclusion functors $i_{\T}:\mathcal{T}\to \mathcal{C}$ and $i_{\F}: \mathcal{F}\to \mathcal{C}.$
\end{enumerate}		
\end{Remark}

Let $(\mathcal{T},\mathcal{F})$ be an s-torsion pair in the extriangulated category $(\mathcal{C},\mathbb{E},\mathfrak{s})$ with negative first extension. For Remark \ref{Rk-inc-ext-tp} (1), we know that the inclusion functor $i_\T:\T\to \C$ induces the 
extriangulated functor $(i_\T,1):(\T, \mathbb{E}|_\T,\mathfrak{s}|_\T)\to (\mathcal{C},\mathbb{E},\mathfrak{s}).$ However, in the case of the torsion functor $\mathfrak{t}:\mathcal{C}\to\mathcal{T},$ this does not automatically imply that this functor is extriangulated. Indeed, for this to hold, one would need a natural transformation 
$\varphi_{C,A}:\mathbb{E}(C,A)\to\mathbb{E}|_\T(\mathfrak{t}(C),\mathfrak{t}(A))$, natural in $A,C\in\mathcal{C}$ and compatible with the corresponding extriangulated realizations.

	The following proposition shows that the torsion functor $\mathfrak{t}:\mathcal{C}\to\mathcal{T}$ is extriangulated if and only if the torsion-free functor $\mathfrak{l}:\mathcal{C}\to\mathcal{F}$ is extriangulated. For this, we recall that an extriangulated category satisfies condicion (WIC), see \cite[Condition 5.8]{NP19}, if for any pair of composable morphisms $gf$ we have: $gf$ is an inflation (deflation) implies that $f$ (resp. $g$) is an inflation (deflation). The main result that we use in the proof of Proposition \ref{pro:torsion_funtor_extriangulado} is \cite[Lem. 5.9]{NP19} where condition (WIC) is required. We point out that, in the case of a triangulated category this condition is not necessary (see the proof of \cite[Lem. 5.9]{NP19}).
	
	\begin{proposition}\label{pro:torsion_funtor_extriangulado}
		Let $(\mathcal{D},\mathbb{E},\mathfrak{s})$ be an extriangulated category with negative first extension and satisfying condition $(WIC)$, and let $(\mathcal{T},\mathcal{F})$ be an s-torsion pair in $\mathcal{D}$. Then the following conditions are equivalent.
		\begin{enumerate}[topsep=1pt]
			\item[(a)] The torsion functor $\mathfrak{t}:\mathcal{D}\to\mathcal{T}$ is extriangulated.
			\item[(b)] The torsion-free functor $\mathfrak{l}:\mathcal{D}\to\mathcal{F}$ is extriangulated.
		\end{enumerate}
Moreover, if $\D$ is triangulated the condition $(WIC)$ can be deleted.		
	\end{proposition}
	
	\begin{proof}
	We will prove in detail that $(b)$ implies $(a)$, while only sketching the converse implication, since the two arguments are largely analogous. Nevertheless, we will retain certain details of the proof that will be essential in subsequent sections.
		
 $(b) \Rightarrow (a)$ Let $(\lf,\varphi):\mathcal{D}\to\mathcal{F}$ be extriangulated and $A\xrightarrow{x}B\xrightarrow{y}C\overset{\delta}{\dashrightarrow}$ be a conflation in $\mathcal{D}$. Thus we obtain the following conflation in $\mathcal{F}$
\begin{center}
  $\lf(A)\xrightarrow{\lf(x)}\lf(B)\xrightarrow{\lf(y)}
		\lf(C)\overset{\delta_1}{\dashrightarrow}$ where $\delta_1:=\varphi_{C,A}(\delta).$
\end{center}		
		
		On the other hand, for each  $X\in\mathcal{D}$, we will consider the canonical conflation $\delta_X:\;\tf(X)\xrightarrow{\iota_X}X\xrightarrow{\pi_X}\lf(X)$, associated with the $s$-torsion pair $(\mathcal{T},\mathcal{F}).$  Thus, we have the following commutative diagram
		
		\begin{equation}\label{dia: completar con tor y free}
			\begin{tikzcd}[ampersand replacement=\&]
				{\tf(A)} \& {\tf(B)} \& {\tf(C)} \& {} \\
				A \& B \& C \& {} \\
				{\lf(A)} \& {\lf(B)} \& {\lf(C)} \& {} \\
				{} \& {} \& {}
				\arrow["{\tf(x)}", from=1-1, to=1-2]
				\arrow["{\iota_A}"', from=1-1, to=2-1]
				\arrow["{\tf(y)}", from=1-2, to=1-3]
				\arrow["{\iota_B}"', from=1-2, to=2-2]
				\arrow["{\iota_C}"', from=1-3, to=2-3]
				\arrow["x", from=2-1, to=2-2]
				\arrow["{\pi_A}"', from=2-1, to=3-1]
				\arrow["y", from=2-2, to=2-3]
				\arrow["{\pi_B}"', from=2-2, to=3-2]
				\arrow["\delta"', dotted, from=2-3, to=2-4]
				\arrow["{\pi_C}"', from=2-3, to=3-3]
				\arrow["{\lf(x)}"', from=3-1, to=3-2]
				\arrow["{\delta_A}", dotted, from=3-1, to=4-1]
				\arrow["{\lf(y)}"', from=3-2, to=3-3]
				\arrow["{\delta_B}", dotted, from=3-2, to=4-2]
				\arrow["{\varphi_{C,A}(\delta)}"', dotted, from=3-3, to=3-4]
				\arrow["{\delta_C}", dotted, from=3-3, to=4-3]
			\end{tikzcd}
		\end{equation}
		
		\textit{Claim:} The sequence $\tf(A)\xrightarrow{\tf(x)}\tf(B)\xrightarrow{\tf(y)}\tf(C)$ is a conflation.
		
		Indeed, consider the following diagram of conflations
		\[
		\begin{tikzcd}[ampersand replacement=\&]
			{\tf(A)} \& {\tf(B)} \&\& {} \\
			A \& B \& C \& {} \\
			{\lf(A)} \& {\lf(B)} \& {\lf(C)} \& {} \\
			{} \& {} \& {}
			\arrow["{\iota_A}"', from=1-1, to=2-1]
			\arrow["{\iota_B}"', from=1-2, to=2-2]
			\arrow["x", from=2-1, to=2-2]
			\arrow["{\pi_A}"', from=2-1, to=3-1]
			\arrow["y", from=2-2, to=2-3]
			\arrow["{\pi_B}"', from=2-2, to=3-2]
			\arrow["\delta"', dotted, from=2-3, to=2-4]
			\arrow["{\lf(x)}"', from=3-1, to=3-2]
			\arrow["{\delta_A}", dotted, from=3-1, to=4-1]
			\arrow["{\lf(y)}"', from=3-2, to=3-3]
			\arrow["{\delta_B}", dotted, from=3-2, to=4-2]
			\arrow["{\varphi_{C,A}(\delta)}"', dotted, from=3-3, to=3-4]
		\end{tikzcd}
		\]
		
		Applying \cite[Lem. 5.9]{NP19} in the diagram above, there are $W\in\mathcal{D}$ and conflations $W\xrightarrow{p}C\xrightarrow{q}\lf(C)\overset{\mu_W}{\dashrightarrow}$ and $\tf(A)\xrightarrow{m}\tf(B)\xrightarrow{n_{y}}W
		\overset{\mu'_W}{\dashrightarrow}$
		such that the following diagram commutes
		\begin{equation}\label{dia: W completa diagrama1}
			\begin{tikzcd}[ampersand replacement=\&]
				{\tf(A)} \& {\tf(B)} \& W \& {} \\
				A \& B \& C \& {} \\
				{\lf(A)} \& {\lf(B)} \& {\lf(C)} \& {} \\
				{} \& {} \& {}
				\arrow["m", from=1-1, to=1-2]
				\arrow["{\iota_A}"', from=1-1, to=2-1]
				\arrow["{n}", from=1-2, to=1-3]
				\arrow["{\iota_B}"', from=1-2, to=2-2]
				\arrow["{\mu'_W}", dotted, from=1-3, to=1-4]
				\arrow["p", from=1-3, to=2-3]
				\arrow["x", from=2-1, to=2-2]
				\arrow["{\pi_A}"', from=2-1, to=3-1]
				\arrow["y", from=2-2, to=2-3]
				\arrow["{\pi_B}"', from=2-2, to=3-2]
				\arrow["\delta"', dotted, from=2-3, to=2-4]
				\arrow["q", from=2-3, to=3-3]
				\arrow["{\lf(x)}"', from=3-1, to=3-2]
				\arrow["{\delta_A}", dotted, from=3-1, to=4-1]
				\arrow["{\lf(y)}"', from=3-2, to=3-3]
				\arrow["{\delta_B}", dotted, from=3-2, to=4-2]
				\arrow["{\varphi_{C,A}(\delta)}"', dotted, from=3-3, to=3-4]
				\arrow["{\mu_W}", dotted, from=3-3, to=4-3]
			\end{tikzcd}
		\end{equation}
		
		and moreover
		  \begin{equation}\label{eq: digrama W }
		 i_{A} \cdot \mu'_{W}= \delta \cdot p,  \quad  \pi_A \cdot \delta = \varphi_{C,A}(\delta) \cdot q, \quad
		  m \cdot \delta_A = \delta_B \cdot \lf(x) \quad n \cdot \delta_B= \mu_{W} \cdot \lf(y).
		  \end{equation}
		
		Notice that  $\iota_B\,\tf(x)=x\,\iota_A$ and $\iota_B\, m=x\,\iota_A,$ and then  from Remark \ref{Rk-inc-ext-tp} (2) we get that $m=\tf(x).$ Moreover,  $W\in\mathcal{T}$ since $\T$ is closed under cones (see Proposition \ref{P-ortTF}). In particular, by Remark \ref{Rk-inc-ext-tp} (2), 
		for the morphism $1_{C}:C \to C$ there exist unique morphism
		$a':W\to \tf(C)$ such that  the first square the following diagram commute
		\[
		\begin{tikzcd}[ampersand replacement=\&]\label{dia:  a y b isos }
			W \& C \& {\lf(C)} \& {} \\
			{\tf(C)} \& C \& {\lf(C)} \& {}
			\arrow["p", from=1-1, to=1-2]
			\arrow["a'"', dashed , from=1-1, to=2-1]
			\arrow["q", from=1-2, to=1-3]
			\arrow[equals, from=1-2, to=2-2]
			\arrow["{\mu_W}", dotted, from=1-3, to=1-4]
			\arrow["b'", dashed ,  from=1-3, to=2-3]
			\arrow["{\iota_C}"', from=2-1, to=2-2]
			\arrow["{\pi_C}"', from=2-2, to=2-3]
			\arrow["{\delta_C}"', dotted, from=2-3, to=2-4]
		\end{tikzcd}
		\]
		 After, for \cite[(ET3)]{NP19} there exist   $b':\lf(C)\to \lf(C)$ such that 
		 $b'q=\pi_C$ and
		 \begin{equation}\label{eq: de los isos a y b}
		 	a' \cdot \mu_{W}=\delta_{C} \cdot b'
		 \end{equation}
		Similarly  there exist unique morphisms $a$ and $b$ such that the following diagram commutes
		\begin{equation}\label{diag: a y a' isomorfismos}
			\begin{tikzcd}[ampersand replacement=\&]
				W \& C \& {\lf(C)} \& {} \\
				{\tf(C)} \& C \& {\lf(C)} \& {} \\
				W \& C \& {\lf(C)} \& {}
				\arrow["p", from=1-1, to=1-2]
				\arrow["{a'}"', from=1-1, to=2-1]
				\arrow["q", from=1-2, to=1-3]
				\arrow[equals, from=1-2, to=2-2]
				\arrow["{\mu_W}", dotted, from=1-3, to=1-4]
				\arrow["{b'}", from=1-3, to=2-3]
				\arrow["{\iota_C}"', from=2-1, to=2-2]
				\arrow["{a}"', from=2-1, to=3-1]
				\arrow["{\pi_C}"', from=2-2, to=2-3]
				\arrow[equals, from=2-2, to=3-2]
				\arrow["{\delta_C}"', dotted, from=2-3, to=2-4]
				\arrow["{b}", from=2-3, to=3-3]
				\arrow["p"', from=3-1, to=3-2]
				\arrow["q"', from=3-2, to=3-3]
				\arrow["{\mu_W}"', dotted, from=3-3, to=3-4]
			\end{tikzcd}
		\end{equation}

		We assert that $a\,a'=1_W$  and $b\,b'=1_{\lf(C)}.$ Indeed, by Remark \ref{Rk-inc-ext-tp} (2) we know that 
		$1_{\lf(C)}:\lf(C)\to \lf(C)$
		and $1_W:W\to W$ are the unique morphisms satisfying $ 1_C\, p=p\, 1_W$ and  $1_{\lf(C)}\, q=q \,1_C$.
		However, from diagram \eqref{diag: a y a' isomorfismos}, we have
		$p\, a\, a' =p\,  1_W$ and $b\, b'\,  q=q\, 1_{\lf(C)}$ and thus the assertion follows. Therefore, $a$ and $b$ are
		isomorphisms. In particular, combining \eqref{eq: de los isos a y b}
		with diagram \eqref{diag: a y a' isomorfismos}, we obtain
		\begin{equation}\label{eq: a y b isos 2}
			i_{C}= p\, a  \qquad   b\, \pi_{C}=q	 \qquad \mu_{W} \cdot b = a \cdot \delta_{C}.
		\end{equation}

		Consider the conflation $\mu_W$ and the isomorphisms $a^{-1}$ and $b.$ Then \cite[Prop. 3.7]{NP19} yields the following conflation 
		$$\tf(C) \xrightarrow{pa} C \xrightarrow{b^{-1}q} \lf(C) \overset{a^{-1}\cdot\mu_W\cdot b}{\dashrightarrow}.$$ 
		In particular, we obtain the following diagram of conflations
		
			\[\begin{tikzcd}[ampersand replacement=\&]
			{\tf(C)} \&\& {C} \& {\lf(C)} \& {} \\
			{\tf(C)} \&\& {C} \& \lf(C) \& {}
			\arrow["{p\,a}", from=1-1, to=1-3]
			\arrow[equals, from=1-1, to=2-1]
			\arrow["{b^{-1}q}", from=1-3, to=1-4]
			\arrow[equals, from=1-3, to=2-3]
			\arrow[equals, from=1-4, to=2-4]
			\arrow["{i_{C}}"', from=2-1, to=2-3]
			\arrow["{\pi_C}"', from=2-3, to=2-4]
			\arrow["{\delta_C}"', dotted, from=2-4, to=2-5]
			\arrow["{ a^{-1} \cdot \mu_{W} \cdot b }", dotted, from=1-4, to=1-5]
		\end{tikzcd}
		\]
	Now consider the conflation $\mu'_{W}$ in \eqref{dia: W completa diagrama1} and the isomorphism $a: \tf(C) \to W.$ By \cite[Prop. 3.7]{NP19}, we have the following coflation
			\begin{equation}\label{diag: conflacion mu W con a}
		\tf(A)\xrightarrow{\tf(x)}\tf(B)\xrightarrow{a^{-1}n}
		\tf(C)\overset{\mu'_{W} \cdot a}{\dashrightarrow}.
		\end{equation}
By Remark \ref{Rk-inc-ext-tp} (2), the morphism $\tf(y)$ is uniquely determined by the equality $i_C\,\tf(y)=y\,\iota_B.$  On the other hand, from \eqref{eq: a y b isos 2} and \eqref{dia: W completa diagrama1} we get $i_C(a^{-1}n)=yi_B$. Hence, by the uniqueness of $\tf(y)$, we obtain $\tf(y)=a^{-1}n$.

Applying \cite[Prop. 3.7]{NP19} to the conflation $\varphi_{C,A}(\delta)$ and using the isomorphism $b^{-1}:\lf(C)\to\lf(C)$, we obtain the conflation

	\[\begin{tikzcd}[ampersand replacement=\&]
	{\lf(A)} \&\& {\lf(B)} \& {\lf(C)} \& {} 
	\arrow["{\lf(x)}", from=1-1, to=1-3]
	\arrow["{b\, \lf(y)}", from=1-3, to=1-4]
	\arrow["{ \hspace{5mm} \varphi_{C,A}(\delta)\cdot b^{-1} }", dotted, from=1-4, to=1-5]
\end{tikzcd}
\]
On the other hand, from \eqref{dia: completar con tor y free}, \eqref{eq: a y b isos 2} and \eqref{dia: W completa diagrama1} and  we get
\[ (b\,\lf(y))\pi_B=b\,\pi_C\,y=q\,y=\lf(y)\,\pi_B=\pi_C\,y. \] 
Hence, by  uniqueness of $\lf(y)$, we obtain $b\,\lf(y)=\lf(y)$. In particular, we have the following diagram of conflations:
\begin{equation}\label{diag: barphi delta b es  varphi delta}
	\begin{tikzcd}[ampersand replacement=\&]
		{\lf(A)} \&\& {\lf(B)} \& {\lf(C)} \& {} \\
		{\lf(A)} \&\& {\lf(B)} \& \lf(C) \& {}
		\arrow["{\lf(x)}", from=1-1, to=1-3]
		\arrow[equals, from=1-1, to=2-1]
		\arrow["{b\, \lf(y)}", from=1-3, to=1-4]
		\arrow[equals, from=1-4, to=2-4]
		\arrow[equals, from=1-3, to=2-3]
		\arrow["{\lf(x)}"', from=2-1, to=2-3]
		\arrow["{\lf(y)}"', from=2-3, to=2-4]
		\arrow["{\varphi_{C,A}(\delta) }"', dotted, from=2-4, to=2-5]
		\arrow["{ \hspace{5mm} \varphi_{C,A}(\delta)\cdot b^{-1} }", dotted, from=1-4, to=1-5]
	\end{tikzcd}
\end{equation}
Thus $\varphi_{C,A}(\delta)=\varphi_{C,A}(\delta)\cdot b^{-1}.$ Moreover, we also have the following commutative diagram
  	\begin{equation}\label{diag: torsión canónica  1}
  	\begin{tikzcd}[ampersand replacement=\&]
  		{\tf(A)} \& {\tf(B)} \& \tf(C) \& {} \\
  		A \& B \& C \& {} \\
  		{\lf(A)} \& {\lf(B)} \& {\lf(C)} \& {} \\
  		{} \& {} \& {}
  		\arrow["\tf(x)", from=1-1, to=1-2]
  		\arrow["{\iota_A}"', from=1-1, to=2-1]
  		\arrow["{\tf(y)}", from=1-2, to=1-3]
  		\arrow["{\iota_B}"', from=1-2, to=2-2]
  		\arrow["{\mu'_W \cdot a}", dotted, from=1-3, to=1-4]
  		\arrow["i_{C}", from=1-3, to=2-3]
  		\arrow["x", from=2-1, to=2-2]
  		\arrow["{\pi_A}"', from=2-1, to=3-1]
  		\arrow["y", from=2-2, to=2-3]
  		\arrow["{\pi_B}"', from=2-2, to=3-2]
  		\arrow["\delta"', dotted, from=2-3, to=2-4]
  		\arrow["\pi_{C}", from=2-3, to=3-3]
  		\arrow["{\lf(x)}"', from=3-1, to=3-2]
  		\arrow["{\delta_A}", dotted, from=3-1, to=4-1]
  		\arrow["{\lf(y)}"', from=3-2, to=3-3]
  		\arrow["{\delta_B}", dotted, from=3-2, to=4-2]
  		\arrow["{\varphi_{C,A}(\delta)}"', dotted, from=3-3, to=3-4]
  		\arrow["{ \delta_C= a^{-1} \cdot \mu_W \cdot b. }", dotted, from=3-3, to=4-3]
  	\end{tikzcd}
  \end{equation}
  We  assert that the vertical and horizontal morphisms in the diagram above define morphisms of the given conflations. Indeed
  
 $$\begin{aligned}
 	i_{A} \cdot (\mu_{W}' \cdot a) & = (i_{A} \cdot \mu_{W}') \cdot a & \\
 	&  = (\delta \cdot p ) \cdot a & ( \text{see } \eqref{eq: digrama W })\\
 	& = \delta \cdot (p\,a) &  ( \text{see }  \eqref{eq: a y b isos 2})\\
 	&= \delta \cdot i_{C} & 
 \end{aligned}$$
 
 $$\begin{aligned}
 	\varphi_{C,A}(\delta) \cdot \pi_{C} &=\varphi_{C,A}(\delta) \cdot (b^{-1}  q) &  ( \text{see} \ \eqref{eq: a y b isos 2}  )\\
 	& =(\varphi_{C,A}(\delta) \cdot b^{-1})  \cdot q &  ( \text{see} \ \eqref{diag: barphi delta b es  varphi delta}  )\\
 	& = \varphi_{C,A}(\delta) \cdot q & (\text{see } \eqref{eq: digrama W }) \\
 	& = \pi_A \cdot \delta. & 
 \end{aligned}$$

 $$\begin{aligned}
 	\tf(y) \cdot \delta_{B} & = (a^{-1}n) \cdot \delta_{B} & (\text{ see } \ref{diag: conflacion mu W con a} ) \\
 	& = a^{-1}\cdot ( n \cdot \delta_{B}) & \\
 	& = a^{-1} \cdot ( \mu_{W} \cdot \lf(y)) & (\text{ see } \ref{dia: W completa diagrama1}) \\
 	& = (a^{-1} \cdot  \mu_{W} ) \cdot \lf(y) &  \\
 	& = (\delta_{C} \cdot b^{-1}) \cdot \lf(y)  & (\text{ see } \ref{eq: de los isos a y b} ) \\
 	& = \delta_{C} \cdot ( b^{-1}  \lf(y)) & \\
 	& = \delta_{C} \cdot \lf(y) & (\text{ see } \ref{diag: barphi delta b es  varphi delta}).
 \end{aligned}$$

  Finally, the equality $\tf(x)\cdot\delta_A=\delta_B\cdot\lf(x)$ follows from \eqref{eq: digrama W }   since $m=\tf(x).$

Hence, by using the diagramm \eqref{diag: torsión canónica  1} we have a map 
$$\psi_{C,A}:\mathbb{E}(C,A)\to
		\mathbb{E}(\tf(C),\tf(A)),\;\delta\mapsto \mu'_W\cdot a.$$
We assert that $\psi$ is a natural transformation. 
		We prove naturality only in the first variable since the the proof of the naturality in second variable follows similarly.

		Let $\alpha:A\to A'$ in $\D.$ Consider the following diagram
		\begin{equation}\label{dia: cuadrado psi natural 1}
			\begin{tikzcd}[ampersand replacement=\&]
				{\mathbb{E}(C,A)} \&\& {\mathbb{E}(\tf(C),\tf(A))} \\
				{\mathbb{E}(C,A')} \&\& {\mathbb{E}(\tf(C),\tf(A')).}
				\arrow["{\psi_{C,A}}", from=1-1, to=1-3]
				\arrow["{\mathbb{E}(C,\alpha)}"', from=1-1, to=2-1]
				\arrow["{\mathbb{E}(\tf(C),\tf(\alpha))}", from=1-3, to=2-3]
				\arrow["{\psi_{C,A'}}"', from=2-1, to=2-3]
			\end{tikzcd}
		\end{equation}
		
		Let $\delta\in\mathbb{E}(C,A)$. Then
		$\alpha\cdot\delta\in\mathbb{E}(C,A')$ is realized by the following diagram
		\begin{equation}\label{digra: de alfa y delta 1}
			\begin{tikzcd}[ampersand replacement=\&]
				A \& B \& C \& {} \\
				{A'} \& {B'} \& C \& {}
				\arrow["x", from=1-1, to=1-2]
				\arrow["\alpha"', from=1-1, to=2-1]
				\arrow["y", from=1-2, to=1-3]
				\arrow["\beta"', from=1-2, to=2-2]
				\arrow["\delta", dotted, from=1-3, to=1-4]
				\arrow[equals, from=1-3, to=2-3]
				\arrow["{x'}"', from=2-1, to=2-2]
				\arrow["{y'}"', from=2-2, to=2-3]
				\arrow["{\alpha\cdot\delta}"', dotted, from=2-3, to=2-4]
			\end{tikzcd}
		\end{equation}
		
		Applying the construction of $\psi$ to $\delta$ and
		$\alpha\cdot\delta,$ via its corresponding diagram of type \eqref{diag: torsión canónica  1}, we get the following diagram  
		
		\begin{equation}\label{dia: nartularidad psi 2}
	\begin{tikzcd}[ampersand replacement=\&]
		{\tf(A)} \&\& {\tf(B)} \&\& {\tf(C)} \& {} \& \\
		\& {\tf(A')} \&\& {\tf(B')} \&\& {\tf(C)} \& {} \\
		A \&\& B \&\& C \& {} \\
		\& {A'} \&\& {B'} \&\& C \& {} \\
		{\lf(A)} \&\& {\lf(B)} \&\& {\lf(C)} \& {} \\
		{} \& {\lf(A')} \& {} \& {\lf(B')} \& {} \& {\lf(C)} \& {} \\
		\& {} \&\& {} \&\& {}
		\arrow["{\tf(x)}", from=1-1, to=1-3]
		\arrow["{\tf(\alpha)}"{pos=0.8}, from=1-1, to=2-2]
		\arrow["{\iota_A}"'{pos=0.2}, from=1-1, to=3-1]
		\arrow["{\tf(y)}", from=1-3, to=1-5]
		\arrow["{\tf(\beta)}"{pos=0.7}, from=1-3, to=2-4]
		\arrow["{\iota_B}"'{pos=0.2}, from=1-3, to=3-3]
		\arrow["{\mu'_W \cdot a}", dotted, from=1-5, to=1-6]
		\arrow["{i_{C}}"{pos=0.3}, from=1-5, to=3-5]
		\arrow["{\tf(x)}"'{pos=0.2}, from=2-2, to=2-4]
		\arrow["{\iota_{A'}}"'{pos=0.2}, from=2-2, to=4-2]
		\arrow["{\tf(y')}"'{pos=0.3}, from=2-4, to=2-6]
		\arrow["{\iota_{B'}}"'{pos=0.2}, from=2-4, to=4-4]
		\arrow["{\mu'_{W'}\cdot a'}", dotted, from=2-6, to=2-7]
		\arrow["{i_{C}}", from=2-6, to=4-6]
		\arrow["x"{pos=0.2}, from=3-1, to=3-3]
		\arrow["\alpha"', from=3-1, to=4-2]
		\arrow["{\pi_A}"'{pos=0.2}, from=3-1, to=5-1]
		\arrow["y"{pos=0.2}, from=3-3, to=3-5]
		\arrow["\beta"', from=3-3, to=4-4]
		\arrow["{\pi_B}"'{pos=0.2}, from=3-3, to=5-3]
		\arrow["\delta"', dotted, from=3-5, to=3-6]
		\arrow[equals, from=3-5, to=4-6]
		\arrow["{\pi_{C}}"{pos=0.2}, from=3-5, to=5-5]
		\arrow["{x'}"'{pos=0.3}, from=4-2, to=4-4]
		\arrow["{\pi_{A'}}"'{pos=0.2}, from=4-2, to=6-2]
		\arrow["{y'}"'{pos=0.2}, from=4-4, to=4-6]
		\arrow["{\pi_{B'}}"'{pos=0.2}, from=4-4, to=6-4]
		\arrow["{\alpha\cdot\delta}", dotted, from=4-6, to=4-7]
		\arrow["{\pi_{C}}"{pos=0.3}, from=4-6, to=6-6]
		\arrow["{\lf(x)}"'{pos=0.2}, from=5-1, to=5-3]
		\arrow["{\delta_A}", dotted, from=5-1, to=6-1]
		\arrow["{\lf(\alpha)}"'{pos=0.7}, from=5-1, to=6-2]
		\arrow["{\lf(y)}"'{pos=0.3}, from=5-3, to=5-5]
		\arrow["{\delta_B}", dotted, from=5-3, to=6-3]
		\arrow["{\lf(\beta)}"', from=5-3, to=6-4]
		\arrow["{\varphi_{C,A}(\delta)}", dotted, from=5-5, to=5-6]
		\arrow["{\delta_W}", dotted, from=5-5, to=6-5]
		\arrow[equals, from=5-5, to=6-6]
		\arrow["{\lf(x')}"'{pos=0.3}, from=6-2, to=6-4]
		\arrow["{\delta_{A'}}"', dotted, from=6-2, to=7-2]
		\arrow["{\lf(y')}"'{pos=0.3}, from=6-4, to=6-6]
		\arrow["{\delta_{B'}}"', dotted, from=6-4, to=7-4]
		\arrow["{\varphi_{C,A'}(\alpha\cdot\delta)}"', dotted, from=6-6, to=6-7]
		\arrow["{\delta_{C}}"', dotted, from=6-6, to=7-6]
		\end{tikzcd}
		\end{equation}
		
	Here, $\tf(\alpha)$, $\tf(\beta)$, $\lf(\alpha)$, and $\lf(\beta)$ are determined by Remark \ref{Rk-inc-ext-tp} (2). Furthermore, since $\varphi: \Ee \to \Ee|_{\mathcal{F}}(\lf-,\lf-)$ is a natural transformation, we have 
	\[ \varphi_{C,A\prime}(\alpha\cdot\delta)=\lf(\alpha)\cdot\varphi_{C,A}(\delta). \] Consequently, the bottom-right square commutes, and the diagram can be completed to obtain a morphism of conflations by taking $1_{\lf(C)}: \lf(C)\to\lf(C)$ as the corresponding vertical morphism.
	
	We now verify that diagram \eqref{dia: cuadrado psi natural 1} commutes. To this end, let $\delta\in\mathbb{E}(C,A)$. On the one hand
		$$ \begin{aligned}
			\Ee(\tf(C), \tf(\alpha))\, \psi_{C,A} (\delta )&= \Ee(\tf(C), \tf(\alpha))(\mu_{W}' \cdot a)\\
			& = \tf(\alpha) \cdot (\mu_{W}' \cdot a)
		\end{aligned}$$

	On the other hand, according to diagram \eqref{dia: nartularidad psi 2}
	
	 $$ \begin{aligned}
	 	 \psi_{C,A'}\mathbb{E}(C,\alpha)(\delta)= & \psi_{C,A'}\bigl(\mathbb{E}(C,\alpha)(\delta)\bigr)\\
	 	& =	\psi_{C,A'}(\alpha\cdot\delta)\\
	 	& =\mu'_{W'}\cdot a'.
	 \end{aligned}
		$$
		 Now, by diagram \eqref{dia: nartularidad psi 2}, we have
		 \[ \begin{aligned}
		 	 i_A\cdot(\mu'_{W}\cdot a)&=\delta\cdot i_C,\\ 
		 	 i_{A\prime}\cdot(\mu_{W'}' \cdot a')&=(\alpha\cdot\delta)\cdot i_C. 
		 \end{aligned} \] 
		 It follows from the first equality that
	 $$\begin{aligned}
	\alpha\cdot ( \delta \cdot i_{C})& =  \alpha \cdot ( i_{A} \cdot (\mu_{W}' \cdot a))\\
	& = (\alpha \, i_{A}) \cdot (\mu_{W}' \cdot a) \\
	& = (i_{A'}\, \tf(\alpha)) \cdot (\mu_{W}' \cdot a)\\
	& = i_{A'} \cdot ( \tf(\alpha) \cdot (\mu_{W}' \cdot a)).
	 \end{aligned}$$
Comparing these equalities, we obtain that
	\begin{equation}\label{eq: cancelar iA}
		\Ee(\tf(C), i_{A'})(\tf(\alpha) \cdot (\mu'_{W} \cdot a))=\Ee(\tf(C), i_{A'})( \mu'_{W'} \cdot a').
	\end{equation}
Now, by applying \cite[Cor. 3.12]{NP19} to the conflation $\tf(A')\xrightarrow{\iota_{A'}} A'\xrightarrow{\pi_{A'}}\lf(A'),$
	  we obtain the exact sequence
	\begin{center}  
	  
 \scalebox{0.8}{$
	\begin{tikzcd}[ampersand replacement=\&]
		{ \C(\tf(C), \tf(A'))}
		\& { \C(\tf(C), A')}
		\& { \C(\tf(C), \lf(A))}
		\& { \Ee(\tf(C), \tf(A'))}
		\&\&
		{ \Ee(\tf(C), A').}
		\arrow[from=1-1, to=1-2]
		\arrow[from=1-2, to=1-3]
		\arrow[from=1-3, to=1-4]
		\arrow["{\Ee(\tf(C), i_{A'})}", from=1-4, to=1-6]
	\end{tikzcd}
	$} 
	\end{center}
	Since $(\T,\F)$ is a torsion pair, we have $\C(\tf(C),\lf(A))=0$. Consequently, the morphism $\Ee(\tf(C),i_{A'})$ is a monomorphism. Hence from \eqref{eq: cancelar iA} we get
	 \[ \tf(\alpha)\cdot(\mu'_{W}\cdot a)=\mu_{W'}'\cdot a'. \]
	 This proves that diagram \eqref{dia: cuadrado psi natural 1} commutes. In particular, the diagram \eqref{dia: nartularidad psi 2} can be completed by taking $1_{\tf(C)}:\tf(C)\longrightarrow\tf(C)$ as the corresponding morphism. Therefore 
	 $(\tf,\psi):
		(\mathcal{D},\mathbb{E},\mathfrak{s})
		\to
		(\mathcal{T},\mathbb{E}|_{\T},
		\mathfrak{s}|_{\T})$
		is an extriangulated functor.
		
		$(b) \Rightarrow (a)$ Let $(\tf,\nu):\mathcal{D}\to\mathcal{T}$ be extriangulated and $A\xrightarrow{x}B\xrightarrow{y}C\overset{\delta}{\dashrightarrow}$ be a conflation in $\mathcal{D}$. Thus we obtain the following conflation in $\mathcal{T}$
		\begin{center}
			$\tf(A)\xrightarrow{\tf(x)}\lf(B)\xrightarrow{\tf(y)}
			\tf(C)\overset{\delta'_1}{\dashrightarrow}$ where $\delta'_1:=\nu_{C,A}(\delta).$
		\end{center}		
		
		On the other hand, for each  $X\in\mathcal{D}$, we will consider the canonical conflation $\delta_X:\;\tf(X)\xrightarrow{\iota_X}X\xrightarrow{\pi_X}\lf(X)$, associated with the $s$-torsion pair $(\mathcal{T},\mathcal{F}).$  Thus, we have the following commutative diagram
		
		\begin{equation}\label{dia: completar con tor y free}
			\begin{tikzcd}[ampersand replacement=\&]
				{\tf(A)} \& {\tf(B)} \& {\tf(C)} \& {} \\
				A \& B \& C \& {} \\
				{\lf(A)} \& {\lf(B)} \& {\lf(C)} \& {} \\
				{} \& {} \& {}
				\arrow["{\tf(x)}", from=1-1, to=1-2]
				\arrow["{\iota_A}"', from=1-1, to=2-1]
				\arrow["{\tf(y)}", from=1-2, to=1-3]
				\arrow["{\iota_B}"', from=1-2, to=2-2]
				\arrow["{\iota_C}"', from=1-3, to=2-3]
				\arrow["x", from=2-1, to=2-2]
				\arrow["{\pi_A}"', from=2-1, to=3-1]
				\arrow["y", from=2-2, to=2-3]
				\arrow["{\pi_B}"', from=2-2, to=3-2]
				\arrow["\delta"', dotted, from=2-3, to=2-4]
				\arrow["{\pi_C}"', from=2-3, to=3-3]
				\arrow["{\lf(x)}"', from=3-1, to=3-2]
				\arrow["{\delta_A}", dotted, from=3-1, to=4-1]
				\arrow["{\lf(y)}"', from=3-2, to=3-3]
				\arrow["{\delta_B}", dotted, from=3-2, to=4-2]
				\arrow["{\nu_{C,A}(\delta)}", dotted, from=1-3, to=1-4]
				\arrow["{\delta_C}", dotted, from=3-3, to=4-3]
			\end{tikzcd}
		\end{equation}
		
		 Similarly to the previous implication, we claim that the sequence
		 $\lf(A)\xrightarrow{\lf(x)}\lf(B)\xrightarrow{\lf(y)}\lf(C)$
		 is a conflation. The proof follows essentially the same arguments, using the dual of \cite[Lem.~5.9]{NP19} together with \cite[Prop.~3.7]{NP19}. This yields a diagram of $\s$-conflations as follows
		 
		 \[\begin{tikzcd}[ampersand replacement=\&]
		 	\& {\tf(A)} \& {\tf(B)} \& {\tf(C)} \& {} \& \\
		 	\& A \& B \& C \& {} \\
		 	\& {\lf(A)} \& {\lf(B)} \& {\lf(C)} \& {} \\
		 	{} \& {} \& {} \& {} \&\& {}
		 	\arrow["{\tf(a)}", from=1-2, to=1-3]
		 	\arrow["{i_{A}}"', from=1-2, to=2-2]
		 	\arrow["{\tf(b)}", from=1-3, to=1-4]
		 	\arrow["{i_{B}}"', from=1-3, to=2-3]
		 	\arrow["{\nu_{C,A}}"{pos=0.3}, dashed, shorten <=2pt, shorten >=2pt, from=1-4, to=1-5]
		 	\arrow["{i_{C}}"', from=1-4, to=2-4]
		 	\arrow["a", from=2-2, to=2-3]
		 	\arrow["{v_{A}}"', from=2-2, to=3-2]
		 	\arrow["b", from=2-3, to=2-4]
		 	\arrow["{v_{B}}"', from=2-3, to=3-3]
		 		\arrow["\delta"{pos=0.3}, dashed, shorten <=2pt, shorten >=2pt, from=2-4, to=2-5]
		 	\arrow["{v_{C}}"', from=2-4, to=3-4]
		 	\arrow["{\lf(a)}"', from=3-2, to=3-3]
		 	\arrow["{\delta'_{A}}", dashed, from=3-2, to=4-2]
		 	\arrow["{\lf(b)}"', from=3-3, to=3-4]
		 	\arrow["{\delta'_{B}}", dashed, from=3-3, to=4-3]
		 	\arrow["{\nu^2_{\lf C, \lf A}.}"{pos=0.7}, dashed, shorten <=2pt, shorten >=2pt, from=3-4, to=3-5]
		 	\arrow["{\delta_{C}}", dashed, from=3-4, to=4-4]
		 \end{tikzcd}\]
		
	\end{proof}

	\begin{lemma}\label{lem:t(F)=0}
		Let $\mathcal{D}$ be an extriangualted category with negative first extension and $(\mathcal{T}, \mathcal{F})$ be an s-torsion pair in $\mathcal{D}.$ Then, for its  associatied functors $\tf:\mathcal{D}\to\mathcal{T}$ and $\lf:\mathcal{D}\to\mathcal{F},$ we have that $\tf(\F)=0,$ $\lf(\T)=0,$ $\lf|_\F\simeq 1_\F$ and $\tf|_\T\simeq 1_\T.$ 
	\end{lemma}
	\begin{proof}
		 Let $F \in \F.$ Consider the canonical conflation $\tf(F)\xrightarrow{\iota_F} F \xrightarrow{\pi_F} \lf(F)\xdashrightarrow{\delta_F}$. Since $0 \in \T \cap \F$ and 
		 $F\in \F$, we also have the conflation $0\xrightarrow{0} F \xrightarrow{1_F} F \xdashrightarrow{0}.$ Thus, by Proposition \ref{Rk-inc-ext-tp} (2), we get the following commutative diagram
			\begin{equation*}
				\begin{tikzcd}[ampersand replacement=\&]
					{\tf(F)} \& F \& {\lf(F)} \& {} \\
					0 \& F \& F \& {}
					\arrow["{{\iota_{F}}}", from=1-1, to=1-2]
					\arrow["a"', from=1-1, to=2-1]
					\arrow["{{\pi_{F}}}", from=1-2, to=1-3]
					\arrow[equals, from=1-2, to=2-2]
					\arrow["\delta_F", dotted, from=1-3, to=1-4]
					\arrow["b", from=1-3, to=2-3]
					\arrow["0"', from=2-1, to=2-2]
					\arrow["{{1_{F}}}"', from=2-2, to=2-3]
					\arrow["0"', dotted, from=2-3, to=2-4]
				\end{tikzcd}
			\end{equation*}
 and therefore $a$ and $b$ are isomorphisms. Hence 	$\tf(\F)=0$	and $\lf|_\F\simeq 1_\F.$	
	\end{proof}
	
	\begin{lemma}\label{con-0-iso}
		For an extriangulated category $\Ccal$ and a conflation $X\xrightarrow{f}Y\xrightarrow{g}Z\xdashrightarrow{\delta},$ the following conditions hold true.
		\begin{enumerate}
			\item[a)]  $X=0$ $\Leftrightarrow$ $g$ is an isomorphism.
			
			\item[b)] $Z=0$ $\Leftrightarrow$ $f$ is an isomorphism.
		\end{enumerate}
	\end{lemma}
	\begin{proof} Let  $X=0.$  Then $\delta \in \Ee(Z,X)=\Ee(Z,0)= \{0 \}.$ On the other hand $\s(0)=[0\xrightarrow{0}Z\xrightarrow{1_{Z}}Z]$, and thus there exists a morphism $a: Y \to Z$ such that the following diagram commutes
			\[\begin{tikzcd}[ampersand replacement=\&]
				{X} \& Y \& Z \& {} \\
				0 \& Z \& Z \& {.}
				\arrow["f", from=1-1, to=1-2]
				\arrow[equals, from=1-1, to=2-1]
				\arrow["g", from=1-2, to=1-3]
				\arrow["a"', from=1-2, to=2-2]
				\arrow["\delta", dotted, from=1-3, to=1-4]
				\arrow[equals, from=1-3, to=2-3]
				\arrow[from=2-1, to=2-2]
				\arrow["{1_{Z}}"', from=2-2, to=2-3]
				\arrow["0"', dotted, from=2-3, to=2-4]
			\end{tikzcd}\]
			Then, by \cite[Cor. 3.6 (1)]{NP19}, $a$ is an isomorphism. In particular 
			$g= 1_{Z} \circ a$ is an isomorphism.
			
	Let now $g$ be an isomorphism. Then, similarly as above, from $(ET3)^{op}$  and 	\cite[Cor. 3.6 (2)]{NP19}, it follows that $X\simeq 0.$
	\end{proof}

	Motivated by the classical notion of a TTF-triple and by the notion of
	an s-torsion pair introduced by
	Adachi--Enomoto--Tsukamoto, we consider the following definition.

	\begin{definition}\label{Ds-TTF} Let $\Ccal$ be an extriangulated category with negative first extension. An \textbf{s-TTF triple} in $\Ccal$ is a triple 
		$(\mathcal T,\mathcal F,\mathcal G)$ such that
		$(\mathcal T,\mathcal F)$  and  $(\mathcal F,\mathcal G)$ are s-torsion pairs in 
		$\Ccal.$
	\end{definition}

	\begin{Remark}\label{Rk-sTTF} Let $\mathcal{C}$ be an extriangulated category with negative first extension and $(\mathcal{T},\mathcal{F},\mathcal{G})$ be an s-TTF triple in $\mathcal{C}.$  By Proposition \ref{P-ortTF}, we know that the central class 
	$\mathcal{F}$ is thick and determines both extremal
	classes. More precisely, $\mathcal{T}={}^{\perp}\mathcal{F}$ and $
	\mathcal{G}=\mathcal{F}^{\perp}$.Thus, once the central class $\mathcal{F}$ is fixed, the classes $\mathcal{T}$ and $\mathcal{G}$ are uniquely determined.
	In this sense, the central class $\F$ 
	encodes the two extremal ones and it is also known as \textbf{s-TTF class}. In particular, any s-TTF class in $\Ccal$ becomes an extriangulated subcategory.
	\end{Remark}

	\section{Localization of Extriangulated Categories}
	
	This section is devoted to recalling some of the main results of
	\cite{NakaokaOgawaSakai2022}, which will provide the notation and
	preliminary results needed to state our main theorem. Throughout the
	section, let $(\mathcal{C},\mathbb{E},\mathfrak{s})$ be an
	extriangulated category. We shall also fix a collection of morphisms
	$\mathscr{S}\subseteq\operatorname{Mor}(\mathcal{C})$ satisfying the
	following condition.
	
	\begin{itemize}
		\item[(M0)] The collection $\mathscr{S}$ contains all isomorphisms
		in $\mathcal{C}$ and is closed under compositions and finite coproducts  $(f_1\coprod f_2\in
		\mathscr{S}(X_1\coprod X_2,Y_1\coprod Y_2)$ if $f_i\in\mathscr{S}(X_i,Y_i)$ for $i=1,2$).
	\end{itemize}
	
	We denote by $\widetilde{\mathcal{C}}$ the localization of
	$\mathcal{C}$ with respect to $\mathscr{S}$. Our goal is to endow
	$\widetilde{\mathcal{C}}$ with a natural extriangulated structure under
	suitable hypotheses.
	
	We begin by associating a full subcategory
	$\mathcal{N}_{\mathscr{S}}\subseteq\mathcal{C}$ to the collection
	$\mathscr{S}$.
	
	\begin{definition}\label{def: Objetos N_S de C para loca 0 to X y X to 0}
		Let $\mathscr{S}\subseteq\operatorname{Mor}(\mathcal{C})$ be as above,
		that is, assume that $(M0)$ holds. Define
		$\mathcal{N}_{\mathscr{S}}\subseteq\mathcal{C}$ to be the class of those $N\in\mathcal{C}$ such that the
		morphisms $N\to 0$ and $0\to N$ belong to $\mathscr{S}$.
	\end{definition}
	
	Notice that $\mathcal{N}_{\mathscr{S}}$ is an additive
	subcategory of $\mathcal{C}$.
	Consequently, for each pair of objects $X,Y\in\mathcal{C}$, the
	collection $\{\mathcal{N}_{\mathscr{S}}(X,Y)\}_{(X,Y)\in\Ccal^2},$ where
	\[
	\mathcal{N}_{\mathscr{S}}(X,Y)
	:=
	\left\{
	f\in\mathcal{C}(X,Y)
	\;\middle|\;
	\begin{array}{c}
		\text{there exist }N\in\mathcal{N}_{\mathscr{S}}
		\text{ and morphisms }f_1:X\to N,\\
		f_2:N\to Y\text{ such that }f=f_2\circ f_1
	\end{array}
	\right\}
	\]
	forms a two-sided ideal of $\mathcal{C}.$ Hence, we have the quotient category
	$\overline{\mathcal{C}}
	:=
	\mathcal{C}/[\mathcal{N}_{\mathscr{S}}]$ and the canonical quotient functor $p:\mathcal{C}\to
	\overline{\mathcal{C}}
	=
	\mathcal{C}/[\mathcal{N}_{\mathscr{S}}].$ For each
	$f\in\mathcal{C}(X,Y)$, we write $\overline{f}
	:=
	f+\mathcal{N}_{\mathscr{S}}(X,Y)$ for the morphism in $\overline{\mathcal{C}}$ represented by $f$.
	Moreover, we denote by $\overline{\mathscr{S}}
	\subseteq
	\operatorname{Mor}(\overline{\mathcal{C}})$
	the closure of $p(\mathscr{S})$ under composition with isomorphisms
	in $\overline{\mathcal{C}}$.

The following consequence 
will be useful in what follows. 

\begin{Remark}\label{obs:  WIC hace que S y S barra iguales}
	If $\mathcal{C}$ is idempotent complete,
	or more generally if every split monomorphism in $\mathcal{C}$ has a
	cokernel in $\mathcal{C}$, that is, if $\mathcal{C}$ is weakly idempotent
	complete, then every collection $\mathscr{S}$ satisfying {\rm(M0)} also
	satisfies condition {\rm(iv)} of \cite[Rk.~3.3]{NakaokaOgawaSakai2022}, and consequently all the equivalent
	conditions in \cite[Lem.~3.2]{NakaokaOgawaSakai2022}. In particular, if
	$(\mathcal{C},\mathbb{E},\mathfrak{s})$ satisfies {\rm(WIC)}, then every
	collection $\mathscr{S}$ satisfying {\rm(M0)} satisfies
	$p(\mathscr{S})=\overline{\mathscr{S}}$.
\end{Remark}

	The main result proved in \cite[Thm. 3.5]{NakaokaOgawaSakai2022} is the following.
	
	\begin{theorem}\label{Teo: Localización NOS}
		For an extriangulated category $(\Ccal,\mathbb{E},\mathfrak{s})$ and a class
		$\mathscr{S}\subseteq\operatorname{Mor}(\mathcal{C})$
		satisfying {\rm(M0)}, the following statements hold true.
		
		\begin{enumerate}
			\item Suppose that $\overline{\mathscr{S}}$ satisfies conditions
			{\rm(MR1)}, {\rm(MR2)}, and {\rm(MR3)} below. Then there exists a
			weakly extriangulated category
			$(\widetilde{\mathcal{C}},\widetilde{\mathbb{E}},
			\widetilde{\mathfrak{s}})$ together with an exact functor $(Q,\mu):
			(\mathcal{C},\mathbb{E},\mathfrak{s})
			\to
			(\widetilde{\mathcal{C}},\widetilde{\mathbb{E}},
			\widetilde{\mathfrak{s}})$.
			\begin{itemize}
				\item[{\rm(MR1)}] The collection
				$\overline{\mathscr{S}}\subseteq
				\operatorname{Mor}(\overline{\mathcal{C}})$ satisfies the
				$2$-out-of-$3$ property with respect to composition in
				$\overline{\mathcal{C}}$.
				
				\item[{\rm(MR2)}] The collection $\overline{\mathscr{S}}$ is a
				multiplicative system in $\overline{\mathcal{C}}$.
				
				\item[{\rm(MR3)}] Let $\left[
				A\xrightarrow{x}B\xrightarrow{y}C,\delta
				\right]$ and $\left[
				A'\xrightarrow{x'}B'\xrightarrow{y'}C',\delta'
				\right]$ be two $\mathfrak{s}$-triangles, and let
				$a\in\mathcal{C}(A,A')$ and
				$c\in\mathcal{C}(C,C')$ be morphisms satisfying $a\cdot\delta=\delta'\cdot c$.\\
				 If $\overline{a}$ and $\overline{c}$ belong to
				$\overline{\mathscr{S}}$, then there exists
				$b\in\overline{\mathscr{S}}(B,B')$ such that
				$b\circ\overline{x}
				=
				\overline{x}'\circ\overline{a}$, $\overline{c}\circ\overline{y}
				=\overline{y}'\circ b$.
			\end{itemize}
			
			\item The exact functor $(Q,\mu):
			(\mathcal{C},\mathbb{E},\mathfrak{s})
			\to
			(\widetilde{\mathcal{C}},\widetilde{\mathbb{E}},
			\widetilde{\mathfrak{s}})$ obtained in {\rm(1)} is characterized by the following universal
			property.
			
			\begin{enumerate}
				\item[{\rm(i)}] Let $(F,\varphi):
				(\mathcal{C},\mathbb{E},\mathfrak{s})
				\to
				(\mathcal{D},\mathbb{F},\mathfrak{t})$	be an exact functor such that $F(s)$ is an isomorphism for every
				$s\in\mathscr{S}$. Then there exists a unique exact functor
				
				$(\widetilde{F},\widetilde{\varphi}):
				(\widetilde{\mathcal{C}},\widetilde{\mathbb{E}},
				\widetilde{\mathfrak{s}})
				\to
				(\mathcal{D},\mathbb{F},\mathfrak{t})$
				such that
				$(F,\varphi)
				=
				(\widetilde{F},\widetilde{\varphi})\circ(Q,\mu)$.
				
				\item[{\rm(ii)}] Let $(F,\varphi),(G,\psi):
				(\mathcal{C},\mathbb{E},\mathfrak{s})
				\to
				(\mathcal{D},\mathbb{F},\mathfrak{t})$
				be exact functors such that $F(s)$ and $G(s)$ are isomorphisms
				for every $s\in\mathscr{S}$. Let $(\widetilde{F},\widetilde{\varphi}),
				(\widetilde{G},\widetilde{\psi}):
				(\widetilde{\mathcal{C}},\widetilde{\mathbb{E}},
				\widetilde{\mathfrak{s}})
				\to
				(\mathcal{D},\mathbb{F},\mathfrak{t})$ be the exact functors obtained in {\rm(i)}. Then, for every
				natural transformation $\eta:(F,\varphi)\rightarrow(G,\psi)$
				between exact functors, there exists a unique natural
				transformation
				$\widetilde{\eta}:
				(\widetilde{F},\widetilde{\varphi})
				\rightarrow
				(\widetilde{G},\widetilde{\psi})$
				between exact functors such that
				$\eta=\widetilde{\eta}\circ(Q,\mu)$.
			\end{enumerate}
			
			\item If, in addition, $\mathscr{S}$ satisfies the condition {\rm(MR4)} below, then
			$(\widetilde{\mathcal{C}},\widetilde{\mathbb{E}},
			\widetilde{\mathfrak{s}})$ is an extriangulated category.
			
			\begin{itemize}
				\item[{\rm(MR4)}] The sets $\mathrm{M}_{\mathrm{inf}}
				:=
				\left\{
				v\circ x\circ u
				\,\middle|\,
				x\text{ is an }\mathfrak{s}\text{-inflation},
				\ u,v\in\mathscr{S}
				\right\}$
				and $\mathrm{M}_{\mathrm{def}}
				:=
				\left\{
				v\circ y\circ u
				\,\middle|\,
				y\text{ is an }\mathfrak{s}\text{-deflation},
				\ u,v\in\mathscr{S}
				\right\}$ are closed under composition in $\mathcal{C}$.
			\end{itemize}
		\end{enumerate}
	\end{theorem}
	
	We shall use this result throughout the present paper. For this reason,
	we next give a brief overview of the main ideas involved in its proof.
	In particular, we focus on the construction of the localization
	$\widetilde{\mathcal{C}}$, the extension
	$\widetilde{\mathbb{E}}$, and the exact realization
	$\widetilde{\mathfrak{s}}$, as well as on the universal property in
	Theorem~\ref{Teo: Localización NOS}{\rm(2)}.

	\subsection{Construction of $\widetilde{\mathcal{C}}$}\label{subsec:constru wide C}
	
	Suppose that $\mathscr{S}$ satisfies {\rm(M0)}, {\rm(MR1)},
	{\rm(MR2)} and {\rm(MR3)}. Assume moreover that $\mathscr{S}=p^{-1}(\overline{\mathscr{S}})$ (for example, if $\mathcal{C}$ is WIC). Then any localization $Q:\mathcal{C}\to\widetilde{\mathcal{C}}$ with
	respect to $\mathscr{S}$ factors uniquely through the canonical
	quotient functor $p:\mathcal{C}\longrightarrow\overline{\mathcal{C}}$.\\
	Consequently, the localization $\widetilde{\mathcal{C}}$ can be
	constructed by first passing to the quotient category
	$\overline{\mathcal{C}}$ and then localizing with respect to the
	multiplicative system $\overline{\mathscr{S}}$.
	
	\begin{definition}\label{def: localiza Q}
		Let
		$
		\overline{Q}:\overline{\mathcal{C}}\to \widetilde{\mathcal{C}}$
		be a localization of $\overline{\mathcal{C}}$ with respect to the
		multiplicative system $\overline{\mathscr{S}}=\langle p(\mathscr{S})\rangle$,
		and consider the composition
		$Q=\overline{Q}\circ p:
		\mathcal{C}\to \widetilde{\mathcal{C}}$.
		Then $Q$ is a localization of $\mathcal{C}$ with respect to
		$\mathscr{S}$.
	\end{definition}
	
	We now consider the quotient category $\overline{\mathcal{C}}$ and the multiplicative system $\overline{\mathscr{S}}$. By {\rm(MR2)},
	$\overline{\mathscr{S}}$ admits both a left and a right calculus of
	fractions. Thus, we may realize the localization as
	\[
	\widetilde{\mathcal{C}}
	:=
	\overline{\mathscr{S}}^{-1}\overline{\mathcal{C}}.
	\]
	We shall use the following notation to describe its objects and
	morphisms.
	
	\begin{itemize}
		\item[(a)] Since $\operatorname{Ob}
		\bigl(\overline{\mathscr{S}}^{-1}\overline{\mathcal{C}}\bigr) =\operatorname{Ob}(\overline{\mathcal{C}})$, we identify the objects of the localization with the objects of
		$\overline{\mathcal{C}}$ and use the same notation for them.
		
		\item[(b)] A morphism
		$X\xrightarrow{\alpha}Y$ in
		$\overline{\mathscr{S}}^{-1}\overline{\mathcal{C}}$ is represented
		by an equivalence class of pairs of morphisms in
		$\overline{\mathcal{C}}$. We denote such a class by $[\overline{s}/\overline{f}]$,
		where $\overline{f}:X\to Y'$ and $\overline{s}:Y\to Y'$, $\overline{s}\in\overline{\mathscr{S}}$.
		Equivalently, it can be represented by the diagram
		\[
		\begin{tikzcd}[ampersand replacement=\&]
			\& {Y'} \& \\
			X \&\& Y
			\arrow["{\overline{f}}", from=2-1, to=1-2]
			\arrow["{\overline{s}}"', from=2-3, to=1-2]
		\end{tikzcd}
		\]
		or, more briefly, by the pair $(\overline{f}:X\to Y',\,\overline{s}:Y\to Y')$.

		\footnote{We use the notation $(s/f)$ so that the morphism
			$s\in\mathscr{S}$ appears on the left-hand side of the
			``fraction.'' Notice that $(s/f)$ denotes an equivalence class of
			pairs, rather than merely a pair of morphisms.}
		When no confusion can arise, we will omit the bars over morphisms in order to simplify the notation.
		
		\item[(c)] Two pairs $(f_1:X\to Y_1,s_1:Y\to Y_1)$ and 
		$(f_2:X\to Y_2,s_2:Y\to Y_2)$
		represent the same morphism if there exist a third pair
		$(f_3:X\to Y_3,s_3:Y\to Y_3)$ and morphisms $u:Y_1\to Y_3$,
		$v:Y_2\to Y_3$ such that		
		$u\circ f_1=f_3=v\circ f_2$ and $u\circ s_1=s_3=v\circ s_2$.

		Equivalently, this can be expressed by the commutative diagram
		\[
		\begin{tikzcd}[ampersand replacement=\&]
			\& {Y_1} \& \\
			X \& {Y_3} \& Y \\
			\& {Y_2} \&
			\arrow["u", from=1-2, to=2-2]
			\arrow["{f_1}", from=2-1, to=1-2]
			\arrow["{f_3}", from=2-1, to=2-2]
			\arrow["{s_1}"', from=2-3, to=1-2]
			\arrow["{s_3}", from=2-3, to=2-2]
			\arrow["{f_2}"', from=2-1, to=3-2]
			\arrow["{s_2}", from=2-3, to=3-2]
			\arrow["v"', from=3-2, to=2-2]
		\end{tikzcd}
		\]
		In particular, the equalities $u\circ s_1=s_3=v\circ s_2$ ensure that the resulting morphisms remain in the multiplicative 	system $\overline{\mathscr{S}}$.
		
		\item[(d)] The composition of the equivalence classes represented by $(f:X\to Y',s:Y\to Y')$
		and $(g:Y\to Z',t:Z\to Z') $ is defined as follows. By {\rm(MR2)}, there exist morphisms $g':Y'\to W$, $s':Z'\to W$
	with $s'\in\overline{\mathscr{S}}$ such that the diagram
		\[\begin{tikzcd}
		&& W && \\
		& {Y'} && {Z'} \\
		X && Y && Z
		\arrow["{g'}", from=2-2, to=1-3]
		\arrow["{s'}"', from=2-4, to=1-3]
		\arrow["f", from=3-1, to=2-2]
		\arrow["s"', from=3-3, to=2-2]
		\arrow["g", from=3-3, to=2-4]
		\arrow["t"', from=3-5, to=2-4]
	\end{tikzcd}\]
		commutes. We then define $[t/g]\circ[s/f]
		:=[s'\circ t\,/\,g'\circ f]$.

		\item[(e)] The identity morphism of an object $X$ in
		$\overline{\mathscr{S}}^{-1}\overline{\mathcal{C}}$ is represented by the class $[\operatorname{id}_X/\operatorname{id}_X]$,
		corresponding to the pair $(\operatorname{id}_X:X\to X,\,
		\operatorname{id}_X:X\to X)$.
	\end{itemize}

	\begin{Remark}\label{obs: sobre fracciones izq derechas}
		The condition {\rm(MR2)} in Theorem~\ref{Teo: Localización NOS}
		implies that $\overline{\mathscr{S}}$ admits both a left and a right
		calculus of fractions. Since $\widetilde{\mathcal{C}}$ is the usual
		localization of an additive category with respect to a multiplicative
		system admitting a calculus of fractions, the localizations obtained
		using left and right fractions are equivalent. We denote them by
		$\widetilde{\mathcal{C}}_{\mathrm{left}}$ and
		$\widetilde{\mathcal{C}}_{\mathrm{right}}$, respectively.
		
		Consequently, every fraction $[s/f]\in
		\widetilde{\mathcal{C}}_{\mathrm{left}}$ can also be represented by a
		fraction of the form $[f'/s']\in
		\widetilde{\mathcal{C}}_{\mathrm{right}}$, and both represent the same
		morphism in the localized category. Whenever no confusion is likely
		to arise, we shall omit the primes and simply write
		$[s/f]=[f/s]$, understanding this equality as an identification of representatives of
		the same morphism, rather than as a literal equality of pairs of
		morphisms.
	\end{Remark}
	
	\subsection{The localization functor} \label{sunsec: funtor localizacion}
	
	According to Definition~\ref{def: localiza Q}, the localization functor 
	$Q:\mathcal{C}\to\widetilde{\mathcal{C}}$ is obtained as the composition of the quotient functor $p:\mathcal{C}\to \overline{\mathcal{C}}$ with the localization functor $\overline{Q}:\overline{\mathcal{C}}\to\overline{\mathscr{S}}^{-1}\overline{\mathcal{C}}$, that is, $Q=\overline{Q}\circ p$.
	
	The functor $Q$ acts as the identity on objects; namely, $Q(X)=X$ for every object $X\in\mathcal{C}$. On morphisms, for each $f:X\to Y$, we have $Q(f)=(\overline{\operatorname{id}_{X}},\overline{f})=[\overline{\operatorname{id}_{X}}/\overline{f}]$, so that $Q(f)$ is represented by the corresponding right fraction.
	
	Moreover, if $s:Y\to Y'$ belongs to $\overline{\mathscr{S}}$, then $Q(s)=[\overline{\operatorname{id}_{Y'}}/\overline{s}]$ is an isomorphism in $\widetilde{\mathcal{C}}$. Consequently, every morphism of $\widetilde{\mathcal{C}}$ can be written in the form
	\[
	[\overline{s}/\overline{f}]
	=
	[\overline{s}/\overline{\operatorname{id}_{Y'}}]
	\circ
	[\overline{\operatorname{id}_{Y'}}/\overline{f}]
	=
	[\overline{\operatorname{id}_{Y'}}/\overline{s}]^{-1}
	\circ
	[\overline{\operatorname{id}_{Y'}}/\overline{f}]
	=
	Q(s)^{-1}\circ Q(f).
	\]
	Dually, every morphism of $\widetilde{\mathcal{C}}$ also admits a representation by right fractions. In what follows, we shall use either description interchangeably, whenever convenient.
	
	\subsection{The bifunctor $\widetilde{\mathbb{E}}$} \label{subsec: bifuntor wilde E}
	
	In this subsection, we recall the construction of the biadditive bifunctor $\widetilde{\mathbb{E}}:\widetilde{\mathcal{C}}^{\op}\times\widetilde{\mathcal{C}}\to\operatorname{Ab}$ and the natural transformation $\mu:\mathbb{E}\Rightarrow\widetilde{\mathbb{E}}\circ(Q^{\op}\times Q)$. To this end, we introduce a biadditive bifunctor $\overline{\mathbb{E}}$ defined on the quotient category, together with natural transformations $\overline{\mu}$ and $\widetilde{\pi}$, so that we obtain the following commutative diagram:
	\begin{equation}\label{dia: constr localiz}
		\xy
		(-12,14)*+{\mathcal{C}^\op\times\mathcal{C}}="0";
		(-12,0)*+{\overline{\mathcal{C}}^\op\times\overline{\mathcal{C}}}="2";
		(-12,-14)*+{\widetilde{\mathcal{C}}^\op\times\widetilde{\mathcal{C}}}="4";
		(14,0)*+{\operatorname{Ab}}="6";
		{\ar_{p^\op\times p} "0";"2"};
		{\ar^{\overline{Q}^\op\times \overline{Q}} "2";"4"};
		{\ar@/_1.5pc/_{Q^\op\times Q} (-19,11.5);(-19,-11.5)};
		{\ar^{\mathbb{E}} "0";"6"};
		{\ar^{\overline{\mathbb{E}}} "2";"6"};
		{\ar_{\widetilde{\mathbb{E}}} "4";"6"};
		{\ar@{=>}_{\widetilde{\pi}} (-4,7);(-4,3)};
		{\ar@{=>}^{\overline{\mu}}  (-4,-3);(-4,-7)};
		{\ar@{}|\circlearrowright (-17,0);(-26,0)};
		\endxy
	\end{equation}
	We then define $\mu:=(\overline{\mu}\ast(p^\op\times p))\circ \widetilde{\pi}$, where $\ast$ denotes the composition of a functor with a natural transformation.
	
	For each $A,C\in\mathcal{C}$, consider the subset $\mathcal{K}(C,A)\subseteq\mathbb{E}(C,A)$ given by $\mathcal{K}(C,A)=\{\delta\in\mathbb{E}(C,A)\mid s\cdot\delta=0\text{ for some }s\in\mathscr{S}(A,A')\}$. By \cite[Lem.~3.13]{NakaokaOgawaSakai2022}, this coincides with the set $\{\delta\in\mathbb{E}(C,A)\mid\delta\cdot t=0\text{ for some }t\in\mathscr{S}(C',C)\}$.
	
	Thus, by \cite[Prop.~3.15]{NakaokaOgawaSakai2022}, the quotient $\mathbb{E}/\mathcal{K}:\mathcal{C}^{\op}\times\mathcal{C}\to\operatorname{Ab},$ where 
	$(\mathbb{E}/\mathcal{K})(C,A):=\mathbb{E}(C,A)/\mathcal{K}(C,A)$, is a biadditive bifunctor such that $(\mathbb{E}/\mathcal{K})(N,-)=0$ and $(\mathbb{E}/\mathcal{K})(-,N)=0$ for every $N\in\mathcal{N}_{\mathscr{S}}$.
	
	\begin{definition}\label{DefEbar}
		The quotient $\mathbb{E}/\mathcal{K}$ induces a biadditive bifunctor $\overline{\mathbb{E}}:\overline{\mathcal{C}}^\op\times\overline{\mathcal{C}}\to\operatorname{Ab}$ as follows.
		
		\begin{enumerate}
			\item For any $A,C\in\mathcal{C}$, set $\overline{\mathbb{E}}(C,A):=(\mathbb{E}/\mathcal{K})(C,A)$.
			
			\item Given $\overline{\delta}\in\overline{\mathbb{E}}(C,A)$, $\overline{a}\in\overline{\mathcal{C}}(A,A')$, and $\overline{c}\in\overline{\mathcal{C}}(C',C)$, define $\overline{a}\cdot\overline{\delta}:=\overline{a\cdot\delta}$ and $\overline{\delta}\cdot\overline{c}:=\overline{\delta\cdot c}$. These operations are well defined, that is, they do not depend on the choice of representatives since $\mathcal{K}\subseteq\mathbb{E}$ is a subfunctor.
		\end{enumerate}
	\end{definition}

	Moreover, there is a natural transformation $\widetilde{\pi}:\mathbb{E}\rightarrow\overline{\mathbb{E}}\circ(p^{\op}\times p)$ whose components are given, for any $A,C\in\mathcal{C}$, by $\widetilde{\pi}_{C,A}:\mathbb{E}(C,A)\to\overline{\mathbb{E}}(C,A)$, $\delta\mapsto\overline{\delta}$. Finally, with a slight abuse of notation, we shall denote an $\overline{\mathbb{E}}$-extension $\overline{\delta}\in\overline{\mathbb{E}}(C,A)$ by $C\overset{\overline{\delta}}{\dashrightarrow}A$.
	
	With the notation introduced above, we now recall the construction of the bifunctor $\widetilde{\mathbb{E}}:\widetilde{\mathcal{C}}^{\op}\times\widetilde{\mathcal{C}}\to\operatorname{Ab}$ and the natural transformation $\overline{\mu}:\overline{\mathbb{E}}\rightarrow\widetilde{\mathbb{E}}\circ(\overline{Q}^{\op}\times\overline{Q})$. The bifunctor $\widetilde{\mathbb{E}}$ is formally defined by means of a coend and is given by $\widetilde{\mathbb{E}}=\widetilde{\mathcal{C}}\otimes_{\overline{\mathcal{C}}}\overline{\mathbb{E}}\otimes_{\overline{\mathcal{C}}}\widetilde{\mathcal{C}}$. To this end, consider the localization functor $\overline{Q}:\overline{\mathcal{C}}\to\widetilde{\mathcal{C}}$ between additive categories. By restriction of scalars along $\overline{Q}$, it induces a functor $\overline{Q}^{*}:\operatorname{Mod}(\widetilde{\mathcal{C}})\to\operatorname{Mod}(\overline{\mathcal{C}})$. In particular, the bifunctor $\widetilde{\mathcal{C}}(-,-)$ acquires an $\overline{\mathcal{C}}$-$\overline{\mathcal{C}}$-bimodule structure via $\overline{Q}^{*}$; namely, $\overline{Q}^{*}(\widetilde{\mathcal{C}}(-,-))=\widetilde{\mathcal{C}}(-,-)\circ(\overline{Q}^{\op}\times\overline{Q})$. Consequently, evaluating the resulting coend at $C,A\in\widetilde{\mathcal{C}}$, we obtain
	\begin{align*}
		\widetilde{\mathbb{E}}(C,A)
		&=
		\widetilde{\mathcal{C}}(-,C)\otimes_{\overline{\mathcal{C}}}
		\left(\overline{\mathbb{E}}\otimes_{\overline{\mathcal{C}}}\widetilde{\mathcal{C}}(A,-)\right) \\
		&\cong
		\int^{Z\in\overline{\mathcal{C}}}
		\widetilde{\mathcal{C}}(Z,C)\otimes_{\overline{\mathcal{C}}}
		\left(\overline{\mathbb{E}}\otimes_{\overline{\mathcal{C}}}\widetilde{\mathcal{C}}(A,-)\right)(Z) \\
		&\cong
		\int^{Z\in\overline{\mathcal{C}}}
		\widetilde{\mathcal{C}}(Z,C)\otimes_{\overline{\mathcal{C}}}
		\left(\overline{\mathbb{E}}(-,Z)\otimes_{\overline{\mathcal{C}}}\widetilde{\mathcal{C}}(A,-)\right) \\
		&\cong
		\int^{Z\in\overline{\mathcal{C}}}
		\left(
		\widetilde{\mathcal{C}}(Z,C)\otimes_{\overline{\mathcal{C}}}
		\int^{X\in\overline{\mathcal{C}}}
		\overline{\mathbb{E}}(X,Z)\otimes_{\overline{\mathcal{C}}}\widetilde{\mathcal{C}}(A,X)
		\right) \\
		&\cong
		\int^{Z\in\overline{\mathcal{C}}}
		\int^{X\in\overline{\mathcal{C}}}
		\widetilde{\mathcal{C}}(Z,C)\otimes_{\overline{\mathcal{C}}}
		\left(
		\overline{\mathbb{E}}(X,Z)\otimes_{\overline{\mathcal{C}}}\widetilde{\mathcal{C}}(A,X)
		\right) \\
		&=
		\int^{(Z,X)\in\overline{\mathcal{C}}^{\op}\times\overline{\mathcal{C}}}
		\widetilde{\mathcal{C}}(Z,C)\otimes_{\overline{\mathcal{C}}}
		\overline{\mathbb{E}}(X,Z)\otimes_{\overline{\mathcal{C}}}
		\widetilde{\mathcal{C}}(A,X).
	\end{align*}
	This holds for all $A,C\in\widetilde{\mathcal{C}}$.
	
	The above coend can be described as a suitable quotient of the abelian group $\coprod_{X,Z\in\overline{\mathcal{C}}}\widetilde{\mathcal{C}}(Z,C)\otimes\overline{\mathbb{E}}(X,Z)\otimes\widetilde{\mathcal{C}}(A,X)$. Since in this work we only recall the construction, we refer the reader to \cite{NakaokaOgawaSakai2022} where the well-definedness of the coend and the corresponding equivalence relation are established in detail. It is worth noting that in that work the equivalence classes are denoted by $[\FR{\overline{t}}{\overline{\delta}}{\overline{s}}]$, with the tensor product symbols omitted explicitly.
	
	We shall adopt this notation and write $[\FR{\overline{t}}{\overline{\delta}}{\overline{s}}]=[C\overset{\overline{t}}{\leftarrow}Z\overset{\overline{\delta}}{\dashrightarrow}X\overset{\overline{s}}{\leftarrow}A]\in\widetilde{\mathbb{E}}(C,A)$ for the equivalence class associated with the triple $(\FR{\overline{t}}{\overline{\delta}}{\overline{s}})=(C\overset{\overline{t}}{\leftarrow}Z\overset{\overline{\delta}}{\dashrightarrow}X\overset{\overline{s}}{\leftarrow}A)$. By definition,
	\[
	\widetilde{\mathbb{E}}(C,A)
	=
	\left\{
	[\FR{\overline{t}}{\overline{\delta}}{\overline{s}}]
	\;\middle|\;
	X,Z\in\overline{\mathcal{C}},\
	\overline{s},\overline{t}\in\overline{\mathscr{S}},\
	\overline{\delta}\in\overline{\mathbb{E}}(Z,X)
	\right\}.
	\]
	
	Now let $\alpha\in\widetilde{\mathcal{C}}(A,A')$ be a morphism represented as $\alpha=\overline{Q}(\overline{u})^{-1}\circ\overline{Q}(\overline{a})$, where $\overline{a}\in\overline{\mathcal{C}}(A,D)$ and $\overline{u}\in\overline{\mathscr{S}}(A',D)$. Consider the following commutative square in $\overline{\mathcal{C}}$:
	\begin{equation}\label{DD}
		\xymatrix{
			X \ar[r]^{\overline{a}'} & X' \\
			A \ar[u]^{\overline{s}} \ar[r]_{\overline{a}} & D \ar[u]_{\overline{s}'}
		}
	\end{equation}
	where $\overline{s}'\in\overline{\mathscr{S}}$. The right action on the bifunctor $\widetilde{\mathbb{E}}$ is then defined by $\alpha\cdot[\FR{\overline{t}}{\overline{\delta}}{\overline{s}}]:=[\FR{\overline{t}}{\overline{a}'\cdot\overline{\delta}}{\overline{s}'\circ\overline{u}}]$. The left action on 
	$\widetilde{\mathbb{E}}$ is defined dually.
	
	On the other hand, for any pair of elements $[\overline{t}_i\backslash \overline{\delta}_i/ \overline{s}_i]\in\widetilde{\mathbb{E}}(C,A)$, $i=1,2$, their sum is defined by
	\[
	[\overline{t}_1\backslash \overline{\delta}_1/\overline{s}_1]+[\overline{t}_2\backslash \overline{\delta}_2/\overline{s}_2]
	:=
	[\overline{t}\backslash(\overline{\rho_1+\rho_2})/\overline{s}],
	\]
	where common denominators (see \cite[Prop. 3.19 (2)]{NakaokaOgawaSakai2022}) are chosen so that $[t_i\backslash\delta_i/s_i]=[t\backslash\rho_i/s]$ for $i=1,2$. This operation is well defined and endows $\widetilde{\mathbb{E}}(C,A)$ with an abelian group structure whose zero element is $[\operatorname{id}\backslash 0/\operatorname{id}]$.
	
	\begin{definition}\label{Def_sN}
		Let $[\FR{\overline{t}}{\overline{\delta}}{\overline{s}}]=[C\overset{\overline{t}}{\leftarrow}Z\overset{\overline{\delta}}{\dashrightarrow}X\overset{\overline{s}}{\leftarrow}A]\in\widetilde{\mathbb{E}}(C,A)$ be a $\widetilde{\mathbb{E}}$-extension and $\mathfrak{s}(\delta)=[X\overset{x}{\longrightarrow}Y\overset{y}{\longrightarrow}Z]$ be the realization $\delta\in\mathbb{E}(C,A)$. The realization $\widetilde{\mathfrak{s}}$ of the element $[\FR{\overline{t}}{\overline{\delta}}{\overline{s}}]$ is defined by
		\[
		\widetilde{\mathfrak{s}}\left([\FR{\overline{t}}{\overline{\delta}}{\overline{s}}]\right)
		:=
		[A\xrightarrow{Q(x\circ s)}Y\xrightarrow{Q(t\circ y)}C].
		\]
	\end{definition}
	
	\subsection{The natural transformation $\varphi: \Ee \to \wE$} \label{subsec: transformacion E to wide E}
	
	Recall that $\varphi: \Ee \to \wE$ is defined as $(\ovl{\mu} \circ (p^{op} \times p)) \ast \widetilde{\pi}$, according to diagram (\refeq{dia: constr localiz}). Thus, it suffices to define a natural transformation $\ovl{\mu}: \ovl{\Ee} \to \wE$. For each pair of objects $C,A$, define the assignment 
	$\ovl{\mu}_{C,A}\co \ovl{\Ee}(C,A)\to \wE(C,A),\;\ovl{\delta}\mapsto [\ovl\id \backslash \ovl\delta/ \ovl\id]$ for every $\ovl{\delta}\in \ovl{\Ee}(C,A)$. Moreover, each $\ovl{\mu}_{C,A}$ is a monomorphism.
Consequently,  $\varphi: \Ee \to \wE$ is given componentwise by the assignment 
$\varphi_{C,A}\co \Ee(C,A)\to \wE(C,A),\;\delta\mapsto\ovl\id \backslash \ovl\delta/ \ovl\id]$.

	\subsection{Universal Property} \label{subsec: propiedad universal locali}
	We recall some aspects of the universal property of the localization functor $Q\co(\C,\Ee,\sfr)\to(\wC,\wE,\ws)$ established in Theorem~\ref{Teo: Localización NOS}(2). In particular, we are interested in describing the construction of the extriangulated functor $(\widetilde F,\widetilde\varphi)\co(\widetilde{\mathcal C},\widetilde{\mathbb E},\widetilde{\mathfrak s})\to(\mathcal D,\mathbb F,\mathfrak t)$.
	
	Let $(F,\phi)\co(\C,\Ebb,\sfr)\to(\D,\Fbb,\tfr)$ be an extriangulated functor that sends every morphism $s\in\Ssc$ to an isomorphism in $\D$. Since the morphisms in $\Ssc$ become isomorphisms in the quotient category, the functor $F$ annihilates the morphisms that are identified with zero in $\ovl{\C}$ and, consequently, factors through the quotient. We denote by $\ovl{F}\co\ovl{\C}\to\D$ the induced functor, which moreover satisfies that $\ovl{F}(s)$ is an isomorphism for every $s\in\ovl{\Ssc}$.
	
	Since $\ovl{Q}\co\ovl{\C}\to\wC$ is the localization of $\ovl{\C}$ with respect to the multiplicative system $\ovl{\Ssc}$, its universal property guarantees the existence and uniqueness of an additive functor $\widetilde{F}\co\wC\to\D$ such that the following diagram commutes:
	
\[
\begin{tikzcd}
	\C & {\ovl{\C}} & \wC \\
	\D
	\arrow["p", from=1-1, to=1-2]
	\arrow["F"', from=1-1, to=2-1]
	\arrow["{\overline{Q}}", from=1-2, to=1-3]
	\arrow["{\overline{F}}", from=1-2, to=2-1]
	\arrow["{\widetilde{F}}", bend left=20, from=1-3, to=2-1]
\end{tikzcd}
\]
	
	Equivalently, we have $F=\widetilde{F}\circ Q$, where $Q=\ovl{Q}\circ p$. In particular, for every fraction $[\overline{s}/\overline{f}]\in\wC$, we obtain the expression $\widetilde{F}([\overline{s}/\overline{f}])=F(s)^{-1}F(f)$.
	
	We now define a natural transformation $\wt{\phi}\co\wE\rightarrow\Fbb\circ(\wt{F}^{\op}\times\wt{F})$ such that the pair $(\wt{F},\wt{\phi})$ is an extriangulated functor.
	
	Let $A,C\in\C.$ Consider an arbitrary element 
	$$[\FR{\ovl{t}}{\ovl{\del}}{\ovl{s}}]=[C\ov{\ovl{t}}{\lla}Z\ov{\ovl{\del}}{\dra}X\ov{\ovl{s}}{\lla}A]\in\wE(C,A)$$ 
	represented $\sfr(\delta)=[X\ov{x}{\lra}Y\ov{y}{\lra}Z]$. Since $F(s)$ and $F(t)$ are isomorphisms in $\D$, we first apply the natural transformation $\phi_{Z,X}\co\Ebb(Z,X)\to\Fbb(FZ,FX)$ to the element $\delta$. This yields a $\tfr$-triangle $\tfr(\phi(\delta))=[FX\ov{Fx}{\lra}FY\ov{Fy}{\lra}FZ]$. We then define $\wt{\phi}_{C,A}\co\wE(C,A)\to\Fbb(F(C),F(A))$ by $\wt{\phi}_{C,A}\big([\FR{\ovl{t}}{\ovl{\del}}{\ovl{s}}]\big):=(F(s)^{-1})\cdot\phi(\delta)\cdot(F(t)^{-1})$.
	
	The well-definedness and additivity of the map $\wt{\phi}_{C,A}\co\wE(C,A)\to\Fbb(F(C),F(A))$ follow from the equivalence relation defining $\wE$ and \cite[Cor. 3.20]{NakaokaOgawaSakai2022}.
	
	\subsection{An special case of localization} 
In this section, we focus on a particular class of localizations, namely, biresolvent categories, which will play a central role in the main result of this work. We start by recalling \cite[Def. 4.1]{NakaokaOgawaSakai2022}.
	
	\begin{definition}
		Let $\mathcal{N} \subseteq \mathcal{C}$ be a thick subcategory. We associate to $\mathcal{N}$ the following classes of morphisms:
		
		\[
		\mathcal{L} = \{ f \in \mathrm{Mor}(\mathcal{C}) \mid \text{$f$ is an inflation and } Cone(f) \in \mathcal{N} \},
		\]
		
		\[
		\mathcal{R} = \{ f \in \mathrm{Mor}(\mathcal{C}) \mid \text{$f$ is a deflation and } CoCone(f) \in \mathcal{N} \}.
		\]
		
		We define $\mathscr{S}_{\mathcal{N}} \subseteq \mathrm{Mor}(\mathcal{C})$ to be the smallest class of morphisms closed under composition and containing both $\mathcal{L}$ and $\mathcal{R}$.
	\end{definition}
	
	 We will occasionally use the notation $\langle \mathcal{L} \cup \mathcal{R} \rangle$ for $\mathscr{S}_{\mathcal{N}}$. We next verify that this class satisfies the hypotheses of Theorem~\ref{Teo: Localización NOS}.
	 
	\begin{lemma}\label{S_NcumpleM0}
		Let $\mathcal{N} \subseteq \mathcal{C}$ be a thick subcategory. Then, the classes $\mathcal{L}$ and $\mathcal{R}$ defined above contain all isomorphisms and are closed under finite coproducts.
	\end{lemma}
	
	\begin{proof}
		We prove the statement for $\mathcal{L}.$ The the proof for $\mathcal{R}$ is dual and it is left to the reader.
		
		Let $f\co X\to Y$ be an isomorphism. Then $X\xrightarrow{f} Y\to 0$ is a conflation and thus $Cone(f)=0\in\mathcal{N}.$ 
		
		Now let $f_i\co X_i\to Y_i$ be in $\mathcal{L}$ for $i=1,2.$ Then, we have conflations $X_i\xrightarrow{f_i}Y_i\xrightarrow{g_i}Z_i$ with $Z_i\in\mathcal{N}$ for 
		$i=1,2.$ Since $Z_1\coprod Z_2\in\mathcal{N}$ and $Cone(f_1\coprod f_2)=Z_1\coprod Z_2,$ it follows that $f_1\coprod f_2\in\mathcal{L}.$
	\end{proof}

	\begin{lemma}\label{S_NcumpleM02}
		Let $\mathcal{N} \subseteq \mathcal{C}$ be a thick subcategory. Then, the class $\mathscr{S}_{\mathcal{N}}$ satisfies {\rm(M0)}.
	\end{lemma}
	
	\begin{proof}
		Since every isomorphism belongs to $\mathcal{L}\cap\mathcal{R}$, and $\mathcal{L}\cup\mathcal{R}\subseteq\mathscr{S}_{\mathcal{N}}$, it follows that $\mathscr{S}_{\mathcal{N}}$ contains all isomorphisms. By definition, $\mathscr{S}_{\mathcal{N}}$ is closed under compositions.
		
		It remains to prove that $\mathscr{S}_{\mathcal{N}}$ is closed under finite coproducts. Let $f,g\in\mathscr{S}_{\mathcal{N}}$. By definition, there exist factorizations $f=f_n\circ\cdots\circ f_1$ and $g=g_m\circ\cdots\circ g_1$, where each $f_i$ and $g_j$ belongs to $\mathcal{L}\cup\mathcal{R}$. Thus, it suffices to show that the coproduct of any two morphisms in $\mathcal{L}\cup\mathcal{R}$ belongs to $\mathscr{S}_{\mathcal{N}}$.
		
		We distinguish the following cases.
		
		\begin{enumerate}
			\item[a)] Suppose that $f_i\co X_i\to Y_i$, for $i=1,2$, are both in $\mathcal{L}$ (respectively, both in $\mathcal{R}$). By Lemma~\ref{S_NcumpleM0}, we have
			$f_1\coprod f_2\in\mathcal{L}$
			(respectively, $\mathcal{R}$). Hence $f_1\coprod f_2\in\mathscr{S}_{\mathcal{N}}$.
			
			\item[b)] Suppose that $f_1: X_1\to Y_1$ belongs to $\mathcal{L}$ and $f_2\co X_2\to Y_2$ belongs to $\mathcal{R}$. consider the composition
			
			\[\begin{tikzcd}[ampersand replacement=\&]
				{X_{1} \coprod X_{2}} \&\& {X_{1} \coprod Y_{2}} \&\& {Y_{1} \coprod Y_{2}}
				\arrow["{\left(\begin{smallmatrix} 1 & 0 \\ 0 & f_{2}\end{smallmatrix}\right)}", from=1-1, to=1-3]
				\arrow["{\left(\begin{smallmatrix} f_{1} & 0 \\ 0 & 1\end{smallmatrix}\right)}", from=1-3, to=1-5]
			\end{tikzcd}
			\]
			Thus
			$f_1\coprod f_2=
			\left(\begin{smallmatrix} f_1 & 0 \\ 0 & 1\end{smallmatrix}\right)
			\left(\begin{smallmatrix} 1 & 0 \\ 0 & f_2\end{smallmatrix}\right)$.
			Moreover
			$\left(\begin{smallmatrix} f_1 & 0 \\ 0 & 1\end{smallmatrix}\right)\in\mathcal{L}$
			and
			$\left(\begin{smallmatrix} 1 & 0 \\ 0 & f_2\end{smallmatrix}\right)\in\mathcal{R}$,
			since $\mathcal{L}$ and $\mathcal{R}$ contain all identity morphisms. Therefore, $f_1\coprod f_2$ is a composition of a morphism in $\mathcal{R}$ followed by a morphism in $\mathcal{L}$, and hence belongs to $\mathscr{S}_{\mathcal{N}}$.
			
			\item[c)] Suppose that $f_1\co X_1\to Y_1$ belongs to $\mathcal{R}$ and $f_2\co X_2\to Y_2$ belongs to $\mathcal{L}$. The argument is analogous to that in (b), using the factorization
			$f_1\coprod f_2=
			\left(\begin{smallmatrix} 1 & 0 \\ 0 & f_2\end{smallmatrix}\right)
			\left(\begin{smallmatrix} f_1 & 0 \\ 0 & 1\end{smallmatrix}\right)$.
		\end{enumerate}
		
		Therefore, $\mathscr{S}_{\mathcal{N}}$ contains all isomorphisms and is closed under compositions and finite direct sums. Hence, $\mathscr{S}_{\mathcal{N}}$ satisfies {\rm(M0)}.
	\end{proof}
	
In general, if $\mathcal{N}$ is a thick subcategory of $\mathcal{C}$, then the class $\langle \mathcal{L} \cup \mathcal{R} \rangle$ satisfies condition {\rm(M3)} by \cite[Lem. 4.6]{NakaokaOgawaSakai2022}, and $\mathscr{S}_{\mathcal{N}}$ forms an almost multiplicative system by \cite[Lem. 4.7]{NakaokaOgawaSakai2022}. However, when $\mathcal{N}$ is biresolving in $\mathcal{C}$, the class $\mathscr{S}_{\mathcal{N}}$ is a multiplicative system.
	
\begin{definition}
	A thick subcategory $\mathcal{N} \subseteq \mathcal{C}$ is called \textbf{biresolving} if, for every $C \in \mathcal{C}$, there exist an inflation $C \to N$ and a deflation $N' \to C$ for some $N,N' \in \mathcal{N}$.
\end{definition}
	
	The authors in \cite{NakaokaOgawaSakai2022} further prove that the class $\mathscr{S}_{\mathcal{N}}=\langle \mathcal{L} \cup \mathcal{R} \rangle$ satisfies all the hypotheses of Theorem~\ref{Teo: Localización NOS}; see \cite[Prop. 4.26]{NakaokaOgawaSakai2022}. More precisely, they establish the following result \cite[Cor. 4.27]{NakaokaOgawaSakai2022}

\begin{proposition}\label{pro: localizacion resolvente 1}
	Let $\mathcal{N} \subseteq \mathcal{C}$ be a biresolving thick subcategory. Then the localization $\mathcal{C}[\mathscr{S}_{\mathcal{N}}^{-1}]$ is a triangulated category.
\end{proposition}

	\section{Recollements of additive categories}
	
	Recollements provide a useful categorical framework for studying an
	ambient category in terms of two smaller categories. They consist of
	adjoint functors that relate the three categories and allow us to
	identify a distinguished subcategory of the ambient category as the
	kernel of a suitable functor. In this section, we recall the definition
	of a recollement of additive categories and some of its basic
	properties. These properties will be used throughout the subsequent
	sections.
	
	The following definition is taken from \cite[Def.2.3]{AMP25}.
	
	\begin{definition}\label{def: recollement aditivas}
		Let $\mathcal{A}$, $\mathcal{B}$, and $\mathcal{C}$ be additive
		categories. A \textbf{recollement} of $\mathcal{B}$ by $\mathcal{A}$
		and $\mathcal{C}$ is a diagram of additive functors
		\[
		\xymatrix{
			\mathcal{A} \ar[r]^{i_*} &
			\mathcal{B} \ar[r]^{j^*}
			\ar@/^1pc/[l]^{i^!}
			\ar@/_1pc/[l]_{i^*} &
			\mathcal{C}
			\ar@/^1pc/[l]^{j_*}
			\ar@/_1pc/[l]_{j_!}
		}
		\]
		satisfying the following conditions 
		\begin{itemize}
			\item[(AR1)] $(i^*,i_*,i^!)$ and $(j_!,j^*,j_*)$ are adjoint triples;
			\item[(AR2)] the functors $i_*$, $j_!$, and $j_*$ are fully faithful;
			\item[(AR3)] $\operatorname{Im}(i_*)=\operatorname{Ker}(j^*)$.
		\end{itemize}
	\end{definition}
	
\begin{Remark}\label{NRk-rec}
		For the recollement of additive categories given in
		Definition~\ref{def: recollement aditivas}, we recall the following basic properties.
		\begin{itemize}
			\item[(a)] It follows from (AR1) and (AR3) that $i^*\circ j_!=0$ and $i^!\circ j_!=0$.
			\item[(b)] For an adjoint pair of functors
			$F:\mathcal{X}\to\mathcal{Y}$, $G:\mathcal{Y}\to\mathcal{X},$
			we have natural transformations
			$\eta:1_{\mathcal{X}}\to G\circ F$ and 
$\varepsilon:F\circ G\to 1_{\mathcal{Y}}$,
			called the \textbf{unit} and the \textbf{counit}, respectively. They
			satisfy the triangular identities
			$\varepsilon_{F(X)}\circ F(\eta_X)=1_{F(X)}$ and 
 $G(\varepsilon_Y)\circ\eta_{G(Y)}=1_{G(Y)}$ for every $X\in\mathcal{X}$ and $Y\in\mathcal{Y},$ see
			\cite[Thm.~3.1.5]{MacLane}.
			
			We use the following notation for the units and counits induced by the adjoint pairs in the given recollement
			\[
			\begin{aligned}
				&\eta_1:1_{\mathcal{B}}\to i_*\circ i^*,
				&&\varepsilon_1:i^*\circ i_*\to 1_{\mathcal{A}},\\
				&\eta_2:1_{\mathcal{A}}\to i^!\circ i_*,
				&&\varepsilon_2:i_*\circ i^!\to 1_{\mathcal{B}},\\
				&\eta_3:1_{\mathcal{C}}\to j^*\circ j_!,
				&&\varepsilon_3:j_!\circ j^*\to 1_{\mathcal{B}},\\
				&\eta_4:1_{\mathcal{B}}\to j_*\circ j^*,
				&&\varepsilon_4:j^*\circ j_*\to 1_{\mathcal{C}}.
			\end{aligned}
			\]
			
		Since $i_*$, $j_!$, and $j_*$ are fully faithful, the morphisms $\varepsilon_1$, $\eta_2$, $\eta_3$, and $\varepsilon_4$ are natural isomorphisms; see \cite[Prop.~3.4.1]{MacLane}. In particular, for every $A\in\mathcal{A}$ and $C\in\mathcal{C}$, the morphisms
		\[
		\varepsilon_1(A),\qquad
		\eta_2(A),\qquad
		\eta_3(C),\qquad
		\varepsilon_4(C)
		\]
		are isomorphisms. Moreover, by the triangular identities for the corresponding adjunctions, the morphisms
		\[
		\eta_1(i_*A),\qquad
		\varepsilon_2(i_* A),\qquad
		\varepsilon_3(j_!C),\qquad
		\eta_4(j_*C)
		\]
		are also isomorphisms whenever they are defined.
		\end{itemize}
	\end{Remark}
	
	The following is a standard result in category theory, see \cite[Prop. 3.4.1, p.~114]{Borceux1} and its dual.
	
	\begin{proposition}\label{pro: F fiel y pleno sii unidad iso}
		For a functor $F\colon\mathcal{X}\to\mathcal{Y},$ the following statements hold true.
		\begin{enumerate}
			\item[(a)] Let $G\colon\mathcal{Y}\to\mathcal{X}$ be a left adjoint of $F$. Then $F$ is full and faithful if, and only if, the counit $\varepsilon\colon G\circ F\xrightarrow{\sim} 1_{\mathcal{X}}.$ 
			\item[(b)] Let $H\colon\mathcal{Y}\to\mathcal{X}$ be a right adjoint of $F$. Then $F$ is full and faithful if, and only if, the unit $\eta\colon 1_{\mathcal{X}}\xrightarrow{\sim} H\circ F.$ 
		\end{enumerate}
	\end{proposition}

	\begin{lemma}\label{lem: adjunto fiel y pleno}
		Let $F\colon\mathcal{D}\to\mathcal{E}$ be a functor which admits a left adjoint $L\colon\mathcal{E}\to\mathcal{D}$ and a right adjoint $R\colon\mathcal{E}\to\mathcal{D}$. Then $L$ is full and faithful if, and only if, $R$ is full and faithful.
	\end{lemma}
	
	\begin{proof}
		Let $\eta\colon 1_{\mathcal{E}}\to F\circ L$ be the unit of the adjunction $(L,F)$ and let $\varepsilon'\colon F\circ R\to 1_{\mathcal{E}}$ be the counit of the adjunction $(F,R)$. For any objects $X,Y\in\mathcal{E}$, consider the following  diagram
		\[
		\begin{tikzcd}
			\operatorname{Hom}_{\mathcal{D}}(L(X),R(Y))
			\arrow[r,"\cong"]
			\arrow[d,"\cong"']
			&
			\operatorname{Hom}_{\mathcal{E}}(X,F(R(Y)))
			\arrow[d,"{\varepsilon'_Y\circ-}"]
			\\
			\operatorname{Hom}_{\mathcal{E}}(F(L(X)),Y)
			\arrow[r,"{-\circ\eta_X}"']
			&
			\operatorname{Hom}_{\mathcal{E}}(X,Y),
		\end{tikzcd}
		\]
		where the horizontal isomorphism (induced by the adjunction $(L,F)$) is given by $f\mapsto F(f)\circ\eta_X$, while the vertical isomorphism (induced by the adjunction $(F,R)$) is given by $f\mapsto\varepsilon'_Y\circ F(f)$. Thus, for every $f\colon L(X)\to R(Y)$, both paths in the above diagram send $f$ to $\varepsilon'_Y\circ F(f)\circ\eta_X$, and hence the diagram commutes.
		
		Since the horizontal and vertical maps in the above diagram are isomorphisms, the natural transformations $\eta$ and $\varepsilon'$ are isomorphisms simultaneously. Indeed, by Yoneda's lemma, the right-hand vertical map is an isomorphism for all $X,Y$ if and only if $\varepsilon'_Y$ is an isomorphism for every $Y.$ Similarly, the bottom horizontal map is an isomorphism for all $X,Y$ if and only if $\eta_X$ is an isomorphism for every $X$. Hence, $\eta$ is a natural isomorphism if and only if $\varepsilon'$ is a natural isomorphism. By Proposition~\ref{pro: F fiel y pleno sii unidad iso}, this is equivalent to $L$ being full and faithful if and only if $R$ is full and faithful.
	\end{proof}
	
	\begin{lemma}\label{LPadjPerp} Let $(F:\X\to \Y,\,G:\Y\to \X)$ be an adjoint par of additive functors between additive categories. Then $\Ker(F)={}^\perp \im(G)$ and  $\Ker(G)=\im(F)^\perp.$
	\end{lemma}
	\begin{proof} It follows Yoneda's lemma and $\Hom_\Y(FX,Y)\simeq\Hom_\X(X,GY)$ $\forall\,X\in\X,$ $\forall\,Y\in\Y.$	
	\end{proof}
	
Based on the definition of recollement in triangulated categories,  we introduce the following notion of s-recollement in extriangulated categories.

\begin{definition}\label{def: recollement alejandro 1}
	Let $\mathcal{A},$ $\Bcal$ and $\mathcal{C}$ be extriangulated categories with $\Bcal$ having negative first extension. An 
	\textbf{s-recollement} of $\mathcal{B}$ by $\mathcal{A}$ and $\mathcal{C}$ is a diagram
	\[
	\begin{tikzcd}[ampersand replacement=\&, column sep=large]
		\mathcal{A} \&\& \mathcal{B} \&\& \mathcal{C}
		\arrow["{i_*}", from=1-1, to=1-3]
		\arrow["{i^!}", bend left=45, from=1-3, to=1-1]
		\arrow["{i^*}"', bend right=45, from=1-3, to=1-1]
		\arrow["{j^*}", from=1-3, to=1-5]
		\arrow["{j_*}", bend left=45, from=1-5, to=1-3]
		\arrow["{j_!,}"', bend right=45, from=1-5, to=1-3]
	\end{tikzcd}
	\]
	satisfying the following conditions.
	\begin{itemize}
	\item[(sER1)] $(i^*,i_*,i^!)$ and $(j_!,j^*,j_*)$ are extriangulated adjoint triples, see Definition \ref{extadjtD}.
	\item[(sER2)] The functors $i_*$, $j_!$, and $j_*$ are fully faithful.
	\item[(sER3)] $\operatorname{Im}(i_*)=\operatorname{Ker}(j^*).$
	\item[(sER4)] $\Ee^{-1}(\Ker(i^*),\im(i_*))=0$ and $\Ee^{-1}(\im(i_*),\Ker(i^!))=0.$ Moreover, for any $X\in\D,$ the unite $\eta_X:X\to i_* i^*X$ is a deflation and the co-unite $\varepsilon_X:i_* i^!X\to X$ is an inflation.
	\end{itemize}
\end{definition}

The term ``s-recollement'' given in the above definition can be explained by the following result.

\begin{proposition}\label{MPsrec} Let $\D$ and $D'$ be extriangulated categories such that $\D$ admits negative first extension. Consider the diagram of extriangulated functors
\begin{equation*}
		\begin{tikzcd}[ampersand replacement=\&, column sep=large]
			{\mathcal{D}'}
			\&\&
			{\mathcal{D}}
			\arrow["{i_{*}}", from=1-1, to=1-3]
			\arrow["{i^{!}}", bend left=35, from=1-3, to=1-1]
			\arrow["{i^{*}}"', bend right=35, from=1-3, to=1-1]
		\end{tikzcd}
	\end{equation*}		
such that $(i^*,i_*,i^!)$ is an adjoint triple and $i_*$ is full and faitfull. Then, the following conditions are equivalent.
\begin{itemize}
\item[(a)] $\im(i_*)$ is closed under extensions and condition (sER4) given in Definition \ref{def: recollement alejandro 1} holds true.
\item[(b)] $(\Ker(i^*),\im(i_*),\Ker(i^!))$ is an s-TTF triple in $\mathcal{D},$ see Definition \ref{defi: s torsion}.
\end{itemize} 
 Moreover if one of the equivalent conditions holds true, then 
for  each $X\in\D$ there exist conflations
\begin{equation*}
	\begin{tikzcd}
		{X_{0}} & X & {i_{\ast}i^{*}(X)} & {}
		\arrow[from=1-1, to=1-2]
		\arrow["{\eta_1(X)}", from=1-2, to=1-3]
		\arrow["{\delta_{X}}", dotted, from=1-3, to=1-4]
	\end{tikzcd}\\
	\begin{tikzcd}[ampersand replacement=\&]
		{i_{*}i^!(X)} \& X \& {X'_{0}} \& {}
		\arrow["{\varepsilon_2(X)}", from=1-1, to=1-2]
		\arrow[from=1-2, to=1-3]
		\arrow["{\delta'_{X}}", dotted, from=1-3, to=1-4]
	\end{tikzcd}
\end{equation*}
where $X_0\in\Ker(i^*)$ and $X'_0\in \Ker(i^!).$ 
\end{proposition}
\begin{proof} By Lemma \ref{LPadjPerp} we have that $\Ker(i^*)={}^\perp \im(i_*)$ and $\Ker(i^!)=\im(i_*)^\perp.$ Then the result follows from \cite[lem. 5.5 (b,c)]{AMP25}.
\end{proof}

Since every extriangulated functor is exact on both sides, in the sense of \cite[Def.~5.1]{AMP25}, and every s-torsion pair is, in particular, a torsion pair, 
we get from Proposition \ref{MPsrec} that any s-recollement is a recollement as defined in \cite[Def. 5.3]{AMP25}. Thus, all the results established in \cite{AMP25} for recollments remain valid for s-recollements. We also recall that, for a class $\Mcal\subseteq \Mor(\Ccal)$ we denote by $\langle \Mcal\rangle$ de closure of $\Mcal$ under compositions with isomorphisms in $\Ccal.$

		\begin{lemma}\label{lem: para supra de adjunto izq} 
		Let $(i^*:\Bcal\to \Acal, i_*:\Acal\to \Bcal)$ be an extriangulated adjoint pair (see Definition \ref{extadjpD}) with $\Bcal$ having negative first extension and such that \((\Ker(i^{*}),\im(i_{*}))\) is an s-torsion pair in $\Bcal.$
		Let \(\E:=\im(i_{*})\) and consider the classes
		\
		
		$(1)\;\Ll= 
		\{f\in \Mor(\B)\mid 
		f \text{ is an }$ $\text{inflation and } 
		Cone(f)\in \E\},$ and 
		
		$(2)\;\R= 
		\{f\in \Mor(\B)\mid 
		f \text{ is a }\text{deflation and }$ $ 
		CoCone(f)\in\E\}$.
		
		 If \(X\in{}^{\perp}\E\), then 
		\[ 
		\Ho_{\B}(X,s): 
		\Ho_{\B}(X,B_1) 
		\xrightarrow{\sim}
		\Ho_{\B}(X,B_2), 
		\] 
		for every \(s:B_1\to B_2\) belonging to \(\langle\Ll\cup\R\rangle\). 
	\end{lemma} 
	\begin{proof} We have the extriangulated adjoint pair
\[ 
		\begin{tikzcd}[ampersand replacement=\&] 
			{(\mathcal{A}, \mathbb{F}, \mathfrak{r})} 
			\&\& 
			{(\mathcal{B}, \mathbb{E}, \mathfrak{s})} 
			\arrow["{(i_{*} , \varphi^{2})}"', bend right=15, from=1-1, to=1-3] 
			\arrow["{(i^{*}, \varphi^{1})}"', bend right=15, from=1-3, to=1-1] 
		\end{tikzcd} 
		\] 	
	
		It suffices to prove the assertion for \(s \in \R \cup \Ll\).

		 Let \(s \in \Ll\).  Then, there is a conflation  $B_1\xrightarrow{s} B_2\to E$ with $E\in \E.$
			Applying $\B(X,-)$ to the preceding conflation, we obtain the following exact sequence
			\[\begin{tikzcd}[ampersand replacement=\&] 
				{\Ee^{-1}(X,E)} \& {\B(X,B_1)} \&\& {\B(X,B_2)} \& {\B(X,E)} 
				\arrow[from=1-1, to=1-2] 
				\arrow["{\B(X,s)}", from=1-2, to=1-4] 
				\arrow[from=1-4, to=1-5] 
			\end{tikzcd}\] 
			Since $(\Ker(i^{*}), \im(i_{*}))$ is an s-torsion pair, it follows that $\Ee^{-1}(X,E)=0$ and $\B(X,E)=0$. Hence, $\B(X,s)$ is an isomorphism.

			 Let \(s\in\R\). Then, there is a conflation $E\to B_1\xrightarrow{s} B_2$ with $E\in \E.$ Applying the functor \(\B(X,-)\) to the preceding conflation, we obtain the exact sequence 
			\[ 
			\begin{tikzcd}[ampersand replacement=\&] 
				{\B(X,E)} 
				\& {\B(X, B_1)} 
				\&\& {\B(X, B_2)} 
				\& {\Ee(X,E)} 
				\arrow[from=1-1, to=1-2] 
				\arrow["{\B(X,s)}", from=1-2, to=1-4] 
				\arrow[from=1-4, to=1-5] 
			\end{tikzcd} 
			\] 
			
			Since \((\Ker(i^{*}),\im(i_{*}))\) is an s-torsion pair,  we have that $\B(X,E)=0$. 
			
			On the other hand, since \(E\in\E=\im(i_{*})\), there exists  \(A\in\A\) such that \(i_{*}(A)\simeq E\). Moreover,  by \cite[Lemma 3.8]{BGLS} we obtain 
			\[ 
			\Ee(X,E) 
			\simeq 
			\Ee(X,i_{*}A) 
			\simeq 
			\mathbb{F}(i^{*}X,A) 
			=0, 
			\] 
			where the last equality follows from \(X\in\Ker(i^{*})\). Consequently $\B(X,s)$ is an isomorphism. 
	\end{proof}
	
The proof of the following result is similar to the proof of Lemma~\ref{lem: para supra de adjunto izq}.

\begin{lemma}\label{lem: para supra de adjunto derecho}
Let $(i_*:\Acal\to \Bcal, i^!:\Bcal\to \Acal)$ be an extriangulated adjoint pair (see Definition \ref{extadjpD}) with $\Bcal$ having negative first extension and such that \((\im(i_{*}), \Ker(i^{!}))\) is an s-torsion pair in $\Bcal.$
		Let \(\E:=\im(i_{*})\) and consider the classes
		\
		
		$(1)\;\Ll= 
		\{f\in \Mor(\B)\mid 
		f \text{ is an }$ $\text{inflation and } 
		Cone(f)\in \E\},$ and 
		
		$(2)\;\R= 
		\{f\in \Mor(\B)\mid 
		f \text{ is a }\text{deflation and }$ $ 
		CoCone(f)\in\E\}$.
		
		 If \(Y\in\E^{\perp}\),  then 
		 $$\Ho_{\B}(s,Y)\colon\Ho_{\B}(B_2,Y)\xrightarrow{\sim}\Ho_{\B}(B_1,Y)$$ 
		  for every $s\colon B_1\to B_2$ belonging to  $\langle\Ll\cup\R\rangle$.
\end{lemma}

\section{Main result}
The following result provides a construction of an s-recollement from a suitable   extriangulated adjoint triple and its associated s-TTF triple. It can be viewed as a generalization, to the setting of extriangulated categories, of the construction of recollements appearing in the work of Cline--Parshall--Scott \cite[Thm. 2.1]{Cline-Parshall-Scott}. More precisely, starting with a fully faithful extriangulated functor $i_*$ admitting both a left and a right extriangulated adjoint, and assuming that the corresponding kernels and image form an s-TTF triple, we construct the third term of the s-recollement as an extriangulated localization of $\mathcal{D}$. The resulting functors $j_!$, $j^*$, and $j_*$ then provide the remaining part of the s-recollement.

\begin{theorem}\label{teo: completar recoll}
	Let $(\D,\Ee,\s)$ and $(D',\Ff,\rr)$ be extriangulated categories such that $\D$ satisfies condition (WIC) and admits negative first extension. Consider the diagram of functors
\begin{equation*}
		\begin{tikzcd}[ampersand replacement=\&, column sep=large]
			{\mathcal{D}'}
			\&\&
			{\mathcal{D}}
			\arrow["{i_{*}}", from=1-1, to=1-3]
			\arrow["{i^{!}}", bend left=35, from=1-3, to=1-1]
			\arrow["{i^{*}}"', bend right=35, from=1-3, to=1-1]
		\end{tikzcd}
	\end{equation*}		
which satisfies the following conditions	
	\begin{enumerate}
		\item[(a)] $(i^{*},i_{*},i^{!})$ is an extriangulated adjoint triple;
		\item[(b)] $i_{*}$ is full and faithful;
		\item[(c)] $(\Ker(i^{*}),\im(i_{*}),\Ker(i^{!}))$ is an s-TTF triple in $\D.$
	\end{enumerate}
	Then $\E:= \im(i_{*})$ is a thick extriangulated subcategory of $\D$ and there exists an extriangulated category $(\D/\E,\Gg,\uu)$ which is triangulated. Moreover,   there are extriangulated functors $j^{*}\colon\D\to \D/\E$, $j_{!}\colon\D/\E\to \D$, and $j_{*}\colon\D/\E\to \D$ such that:
	\begin{enumerate}
		\item $(j_{!},j^{*},j_{*})$ is an extriangulated adjoint triple;
		\item $j_{!}$ and $j_{*}$ are full and faithful;
		\item $\Ker(j^*)=\im(i_*),$ $\im(j_{!})=\Ker(i^*)$ and $\Ker(i^{!})=\im(j_*).$
	\end{enumerate}
	In particular, we obtain the following s-recollement of extriangulated categories
	\begin{equation}\label{dia: recollment 1}
		\begin{tikzcd}[ampersand replacement=\&, column sep=large]
			{\mathcal{D}'}
			\&\&
			{\mathcal{D}}
			\&\&
			{\mathcal{D}/\mathcal{E}}
			\arrow["{i_{*}}", from=1-1, to=1-3]
			\arrow["{i^{!}}", bend left=35, from=1-3, to=1-1]
			\arrow["{i^{*}}"', bend right=35, from=1-3, to=1-1]
			\arrow["{j^{*}}", from=1-3, to=1-5]
			\arrow["{j_{*}}", bend left=35, from=1-5, to=1-3]
			\arrow["{j_{!}}"', bend right=35, from=1-5, to=1-3]
		\end{tikzcd}
	\end{equation}
\end{theorem}

\begin{proof}
Since $\E:= \im(i_{*}),$  by (c) and Proposition \ref{P-ortTF},  we get that $\E$ is a thick extriangulated subcategory of $\D.$ Now consider the following classes of morphisms
\[
\mathcal{L}=\{f\in \operatorname{Mor}(\D) \mid f \text{ is an }\text{inflation with } \operatorname{Cone}(f)\in \mathcal{\E}\},
\]
\[
\mathcal{R}=\{f\in \operatorname{Mor}(\D) \mid f \text{ is a }\text{deflation with } \operatorname{CoCone}(f)\in \mathcal{E}\}.
\]

We set $S_{\mathcal{E}}:=\langle \mathcal{L}\cup \mathcal{R}\rangle$ which is the smallest class of morphisms closed under compositions of elements in $\mathcal{L}\cup \mathcal{R}$ and isomorphisms. In what follows, we use de notation introduced in Remark \ref{NRk-rec} for the unite and the co-unite of adjoint pairs. By Proposition \ref{MPsrec}, for  each $X\in\D$ there exist conflations
\begin{equation}\label{diag: sucesion de torsion -1}
	\begin{tikzcd}
		{X_{0}} & X & {i_{\ast}i^{*}(X)} & {}
		\arrow[from=1-1, to=1-2]
		\arrow["{\eta_1(X)}", from=1-2, to=1-3]
		\arrow["{\delta_{X}}", dotted, from=1-3, to=1-4]
	\end{tikzcd}\\
	\begin{tikzcd}[ampersand replacement=\&]
		{i_{*}i^!(X)} \& X \& {X'_{0}} \& {}
		\arrow["{\varepsilon_2(X)}", from=1-1, to=1-2]
		\arrow[from=1-2, to=1-3]
		\arrow["{\delta'_{X}}", dotted, from=1-3, to=1-4]
	\end{tikzcd}
\end{equation}
where $i^*(X_0)=0$ and $i^!(X'_0)=0.$ In particular, we have that $\mathcal{E}$ is a biresolving category. Hence by Lemma \ref{S_NcumpleM02} and \cite[Prop. 4.26]{NakaokaOgawaSakai2022},  $S_{\mathcal{E}}$ satisfies conditions: (M0), (M1), (M2), (M3) and (M4). Therefore, by \cite[Cors. 3.4 and 4.27]{NakaokaOgawaSakai2022} there exists an extriangulated category $(\mathcal{D}/\mathcal{E}, \mathbb{G}, \mathfrak{u})$ which is also triangulated, and the extriangulated functor
\[
(Q, \psi^{1}): (\D, \Ee, \s)\to (\D/\E, \Gg, \uu)
\]
defined as follows
\begin{enumerate}
	\item[(i)] For $f\in\D(X,Y),$ we have  $Q(f):=[\overline{1_X}/ \overline{f}]$ where $\overline{f} \in \D/ \N_{S_{\mathcal{E}}}(X,Y)$ (see Subsection \ref{sunsec: funtor localizacion});
	\item[(ii)] the natural transformation $\psi^{1}:\Ee\to \Gg\circ ({Q}^{op}\times Q),$ where for each pair of objects $A,C\in\D$, we set
	\[
	\psi^{1}_{C,A}(\delta):=[\overline{1_C}\backslash \overline{\delta}/\overline{1_A}]
	\]
	and $\overline{\delta}\in\overline{\Ee}(C,A):= \Ee(C,A)/ \mathcal{K}(C,A)$ (see Subsection \ref{subsec: transformacion E to wide E}).
\end{enumerate}

We set $j^*:=Q:\D\to \D/\E.$ Let us show, firstly, that  $\operatorname{Ker}(j^*)=\mathcal E:=\operatorname{Im}(i_*).$ By Remark \ref{Rthick-sub} (2) , 
$\operatorname{Ker}(j^*)$ is a thick extriangulated subcategory of $\D$. Moreover, since $j^*$ is the localization functor associated with the multiplicative system $\mathcal S_{\mathcal E}=\langle\mathcal L\cup\mathcal R\rangle$, its kernel is the smallest thick extriangulated subcategory of $\D$ containing $\mathcal E$.
On the other hand,   by the construction of the localization $j^*:=Q:\D\to \D/\E.,$ every object of $\mathcal E$ is isomorphic to the zero object in $\D/\mathcal E.$ Hence $\mathcal E\subseteq\operatorname{Ker}(j^*)$ and therefore $\operatorname{Ker}(j^*)=\mathcal E.$

In what follows, we find the left and right adjoints of $j^*:=Q:\D\to \D/\E.$ We divide it into several steps.

\textbf{(Step 1) There exists an extriangulated functor $(\Phi, \overline{\psi^{2}}): \D \longrightarrow \D$.}\\
Let $X\in\D$ and consider the unit $\eta_1:1_{\D}\to i_*i^*$ of the extriangulated adjoint pair $(i^{*},i_{*}).$ By (\ref{diag: sucesion de torsion -1}) there is a conflation
\begin{equation}\label{diag: sucesion de torsion 1}
	\begin{tikzcd}
		{X_{0}} & X & {i_{\ast}i^{*}(X)} & {}
		\arrow[from=1-1, to=1-2]
		\arrow["{\eta_1(X)}", from=1-2, to=1-3]
		\arrow["{\delta_{X}}", dotted, from=1-3, to=1-4]
	\end{tikzcd}
\end{equation}
such that $X_{0}\in\operatorname{Ker}(i^{*}).$ Hence, \eqref{diag: sucesion de torsion 1} corresponds to the canonical sequence associated with the object $X$ via the s-torsion pair $(\operatorname{Ker}(i^{*}),\operatorname{Im}(i_{*})).$ Thus $\tf(X)=X_0$ where $\tf:\D\to \operatorname{Ker}(i^{*})$ is the torsion functor associated to the s-torsion pair $(\operatorname{Ker}(i^{*}),\operatorname{Im}(i_{*})),$ see Remark \ref{Rk-inc-ext-tp} (2).

Consider the inclusion functor $i:\operatorname{Ker}(i^{*})\to \D$ and the extriangulated functors  $(i^{*}, \varphi^{2}):(\mathcal{D}, \mathbb{E}, \mathfrak{s})\to (\mathcal{D}', \mathbb{F}, \mathfrak{r})$ and $(i_{*},\varphi^{1}):(\mathcal{D}', \mathbb{F}, \mathfrak{r})\to (\mathcal{D}, \mathbb{E}, \mathfrak{s}).$  It follows from Remark \ref{Rk-inc-ext-tp} that 
\[
(i, \mu^1):(\operatorname{Ker}(i^*), \Ee\mid_{\operatorname{Ker}(i^*)}, \s\mid_{\operatorname{Ker}(i^*)})\to (\D, \Ee, \s)
\]
is an extriangulated functor where $\mu^1$ is the identity. We set  $\Phi:=i\circ\tf:\D\to \D$. Explicitly, for $f:X\to Y$ in $\D,$ we have $\Phi(f)=i(\tf(f))=\tf(f)$, where $\tf(f)$ is the unique morphism making the following diagram commute
\begin{equation}\label{diag: funtor torsion 3 con eta}
	\begin{tikzcd}
		{\tf(X):=X_{0}} & X & {i_{\ast}i^{*}(X)} & {} \\
		{\tf(Y):=Y_{0}} & Y & {i_{\ast}i^{*}(Y)} & {}
		\arrow[from=1-1, to=1-2]
		\arrow["{\tf(f)}"', from=1-1, to=2-1]
		\arrow["{\eta_1(X)}", from=1-2, to=1-3]
		\arrow["f"', from=1-2, to=2-2]
		\arrow["{\delta_{X}}", dotted, from=1-3, to=1-4]
		\arrow["{i_{\ast}i^{*}(f)}", from=1-3, to=2-3]
		\arrow[from=2-1, to=2-2]
		\arrow["{\eta_1(Y)}"', from=2-2, to=2-3]
		\arrow["{\delta_{Y}}"', dotted, from=2-3, to=2-4]
	\end{tikzcd}
\end{equation}

Since $(i^{*},\varphi^{2})$ and $(i_{*},\varphi^{1})$ are extriangulated functors, their composition $(i_{*}i^{*},\mu^{3})$ is also an extriangulated functor. The natural transformation $\mu^{3}$ associated with this composition is given by the following composition of natural transformations:
\[
\begin{tikzcd}[ampersand replacement=\&]
	{\mathbb{E}} \& {\mathbb{F}\circ({i^{*}}^{op} \times i^{*})} \&\& {(\mathbb{E}\circ({i_{*}}^{op} \times i_{*})) \circ ({i^{*}}^{op} \times i^{*})}
	\arrow["{\varphi^2}", Rightarrow, from=1-1, to=1-2]
	\arrow["{\varphi^1 \blacktriangle ({i^{*}}^{op} \times i^{*})}", Rightarrow, from=1-2, to=1-4]
\end{tikzcd}
\]
where $\blacktriangle$ denotes \textbf{whiskering}; that is,
\[
\mu^3_{C,A}
=\left[\left(\varphi^1\blacktriangle({i^{*}}^{op}\times i^{*})\right)\circ\varphi^2\right]_{C,A}
=\left(\varphi^1\blacktriangle({i^{*}}^{op}\times i^{*})\right)_{C,A}\circ\varphi^2_{C,A}
=\varphi^1_{i^*(C),i^*(A)}\circ\varphi^2_{C,A}
\]
for all $A,C\in\D$. Now, since $(i_{*}i^*,\mu^{3}), (1_{\D},1_{\Ee}):(\D,\Ee,\s)\longrightarrow(\D,\Ee,\s)$ are extriangulated functors,\footnote{Notice that $i_{*}i^*$ is the composition of the torsion-free functor $\lf:\D\to \operatorname{Im}(i_{*})$ with the inclusion functor from $\operatorname{Im}(i_{*})$ into $\D$.} Proposition \ref{pro:torsion_funtor_extriangulado} implies that the torsion functor 
$$(\tf,\mu^{2}):(\D,\Ee,\s)\to (\operatorname{Ker}(i^{*}),\Ee|_{\Ker(i^{*})},\s|_{\operatorname{Ker}(i^{*})})$$
 is an extriangulated functor.\footnote{For a precise definition of $\mu^{2}$, see the proof of Proposition \ref{pro:torsion_funtor_extriangulado}. Essentially, $\mu^{2}$ is obtained by restriction, since in that proof it is defined as $\mu'_{W}\cdot a'$, where $a'$ is an isomorphism.}

Consequently, $\Phi$ is the composition of two extriangulated functors and is therefore itself an extriangulated functor. Notice that the associated natural transformation of $\Phi$ is $\overline{\psi^{2}}:\Ee\Longrightarrow\Ee(\Phi^{op}-,\Phi-)$ which is  the composition of the natural transformations
\[
\begin{tikzcd}[ampersand replacement=\&]
	{\mathbb{E}} \& {\mathbb{E}\circ({\tf}^{op}\times\tf)} \&\& {(\mathbb{E}\circ({\tf}^{op}\times\tf))\circ({i}^{op}\times i).}
	\arrow["{\mu^2}", Rightarrow, from=1-1, to=1-2]
	\arrow["{\mu^1\blacktriangle({i}^{op}\times i)}", Rightarrow, from=1-2, to=1-4]
\end{tikzcd}
\]	

\textbf{(Step 2) The functor $\Phi$ from  (Step 1) induces an extriangulated functor $j_!:=\overline{\Phi}:\D/\E\to \D$.}

Indeed, since $\D/\E$ is a localization, it suffices to define $\overline{\Phi}$ by means of the universal property discussed in Subsection \ref{subsec: propiedad universal locali}. We first show that $\Phi(f)=\tf(f)$ is an isomorphism for every $f\in S_{\E}$. Since $S_{\E}$ is the smallest subclass of morphisms closed under composition and containing both $\mathcal{L}$ and $\mathcal{R}$, it suffices to show that $\tf(f)$ is an isomorphism for every $f\in\mathcal{L}\cup\mathcal{R}$.\\

Let $f:X \to Y$ be a morphism in $\mathcal{L}$. By definition, there exists a conflation
\[
\begin{tikzcd}[ampersand replacement=\&]
	X \& Y \& Z \& {}
	\arrow["f", from=1-1, to=1-2]
	\arrow["g", from=1-2, to=1-3]
	\arrow["\delta", dotted, from=1-3, to=1-4]
\end{tikzcd}
\]
with $Z\in\E$. By Proposition \ref{pro:torsion_funtor_extriangulado}, we have the following commutative diagram of conflations:

\begin{equation}\label{dia: f es iso para Phi}
	\begin{tikzcd}[ampersand replacement=\&]
		{\tf(X)} \& {\tf(Y)} \& {\tf(Z)} \& {} \\
		X \& Y \& Z \& {} \\
		{i_{*}i^{*}(X)} \& {i_{*}i^{*}(Y)} \& {i_{*}i^{*}(Z)} \& {} \\
		{} \& {} \& {}
		\arrow["{\tf(f)}", from=1-1, to=1-2]
		\arrow["{{{\iota_X}}}"', from=1-1, to=2-1]
		\arrow["{\tf(g)}", from=1-2, to=1-3]
		\arrow["{{{\iota_Y}}}"', from=1-2, to=2-2]
		\arrow["{{{ \mu^2_{C,A}(\delta)}}}", dotted, from=1-3, to=1-4]
		\arrow["{\iota_{Z}}", from=1-3, to=2-3]
		\arrow["x", from=2-1, to=2-2]
		\arrow["{\eta_1(X)}"', from=2-1, to=3-1]
		\arrow["y", from=2-2, to=2-3]
		\arrow["{\eta_1(Y)}"', from=2-2, to=3-2]
		\arrow["\delta"', dotted, from=2-3, to=2-4]
		\arrow["{\eta_1(Z)}", from=2-3, to=3-3]
		\arrow["{{{\lf(x)}}}"', from=3-1, to=3-2]
		\arrow["{{{\delta_{X}}}}", dotted, from=3-1, to=4-1]
		\arrow["{{{\lf(y)}}}"', from=3-2, to=3-3]
		\arrow["{{{\delta_Y}}}", dotted, from=3-2, to=4-2]
		\arrow["{{{ \mu^3_{C,A}(\delta)}}}"', dotted, from=3-3, to=3-4]
		\arrow["{{\delta_{Z}}}", dotted, from=3-3, to=4-3]
	\end{tikzcd} 
\end{equation}

Since $Z\in\E:=\operatorname{Im}(i_{*})$ and $(\operatorname{Ker}(i^{*}),\operatorname{Im}(i_{*}))$ is an s-torsion pair, Lemma \ref{lem:t(F)=0}(a) implies that $t(Z)=0$. Consequently, Lemma \ref{con-0-iso}  implies that $\tf(f)$ is an isomorphism. Similarly, we get that $\tf(f)$ is an isomorphism for $f\in\mathcal{R}.$
Therefore, by the universal property of localization given in Subsection \ref{subsec: propiedad universal locali}, there exists a unique extriangulated functor
\[
(\overline{\Phi},\psi^2):(\D/\E,\Gg,\uu)\longrightarrow(\D,\Ee,\s)
\]
such that the following diagram commutes
\begin{equation}\label{diag:  Phi pasa a j* de la localizacion}
	\begin{tikzcd}[ampersand replacement=\&]
		{(\D, \Ee, \s)} \&\& {(\D,\Ee,\s)} \\
		\\
		{(\D / \E, \Gg, \uu).}
		\arrow["{(\Phi,\overline{\psi^2})}", from=1-1, to=1-3]
		\arrow["{(Q, \psi^1)}"', from=1-1, to=3-1]
		\arrow["{(\overline{\Phi}, \psi^2)}"', from=3-1, to=1-3]
	\end{tikzcd}
\end{equation}
More explicitly, the functor $j_!:=\overline{\Phi}$ and the natural transformation $\psi^{2}$ are defined as follows:

\begin{enumerate}
	\item[a)] Let $\alpha\in\operatorname{Hom}_{\D/\E}(X,Y)$ be represented by a diagram of the form
	\[
	\begin{tikzcd}[ampersand replacement=\&]
		\& Z \& \\
		X \&\& Y
		\arrow["s"', from=1-2, to=2-1]
		\arrow["g", from=1-2, to=2-3]
	\end{tikzcd}
	\]
	with $s\in S_{\E}$. Then $\alpha=Q(g)\circ Q(s)^{-1}$ and consequently $\overline{\Phi}(\alpha):=\Phi(g)\Phi(s)^{-1}$.
	
	\item[b)] Let
	$[\overline{t}\di\overline{\delta}/\overline{s}]
	=
	\left[
	\begin{tikzcd}[ampersand replacement=\&]
		C \& Z \& X \& A
		\arrow["{\overline{t}}"', from=1-2, to=1-1]
		\arrow["{\overline{\delta}}", dotted, from=1-2, to=1-3]
		\arrow["{\overline{s}}"', from=1-4, to=1-3]
	\end{tikzcd}
	\right]$ be an element of $\Gg(C,A)$. We define
	\[
	\psi^{2}([\overline{t}\di\overline{\delta}/\overline{s}])
	:=
	\Phi(s)^{-1}\cdot\overline{\psi^{2}}(\delta)\cdot\Phi(t)^{-1}.
	\]
\end{enumerate}
In particular $\im(j_!)=\im(\Phi)=\im(\tf)=\Ker(i^*).$
	
	\textbf{(Step 3) The functor  $j_!:=\overline{\Phi}$ defined in (Step 2) is a left adjoint of $j^*:=Q:\D\to \D/\E$.}
	\
	
	To this end, we show that there is a functorial isomorphism
	\[
	\Hom_{\D}(\overline{\Phi}(X),Y)\simeq\Hom_{\D/\E}(X,Q(Y))
	\]
	for $X\in\D/\E$ and $Y\in\D$. By construction, $\overline{\Phi}(X)=\tf(X)=X_{0}$ and $Q(Y)=Y$. We show there is a functorial isomorphism $\Hom_{\D}(X_{0},Y)\simeq\Hom_{\D/\E}(X,Y)$.\\
	We first define an isomorphism $\Gamma_{X,Y}:\Hom_{\mathcal D}(X_{0},Y)
	\to
	\Hom_{\mathcal D/\mathcal E}(X,Y)$.
	By construction (see \eqref{diag: sucesion de torsion 1}), we have a conflation in $\mathcal D$
	\begin{equation}
		\begin{tikzcd}
			{X_{0}} & X & {i_{\ast}i^{*}(X)} & {}
			\arrow["{u_X}", from=1-1, to=1-2]
			\arrow["{\eta_1(X)}", from=1-2, to=1-3]
			\arrow["{\delta_{X}}", dotted, from=1-3, to=1-4]
		\end{tikzcd}
	\end{equation}
	such that $i^{*}(X_{0})=0$. Hence, by the construction of $S_{\E}$, we have $u_X\in S_{\E}$. Given $\theta\in\Hom_{\mathcal D}(X_{0},Y)$, we can form the following morphism $\alpha$ from $X$ to $Y$ in $\mathcal D/\E$
	\[
	\begin{tikzcd}[ampersand replacement=\&]
		\& {X_{0}} \& \\
		X \&\& Y.
		\arrow["{u_{X}}"', from=1-2, to=2-1]
		\arrow["\theta", from=1-2, to=2-3]
		\arrow[
		"\alpha"',
		decorate,
		decoration={snake, amplitude=.4mm, segment length=2mm},
		from=2-1,
		to=2-3
		]
	\end{tikzcd}
	\]
	Thus $\alpha:=Q(\theta)Q(u_X)^{-1}:X
	\mathrel{\tikz[baseline=-0.5ex]{
			\draw[->,decorate,decoration={snake, amplitude=.4mm, segment length=2mm}]
			(0,0) -- (1.5,0);
	}}Y$. More precisely, according to the construction in Subsection \ref{subsec:constru wide C}, we have
	\[
	\alpha=Q(\overline{\theta})Q(\overline{u_X})^{-1}.
	\]
	However, since $\D$ is an extriangulated category satisfying condition (WIC), we get   from Remark \ref{obs:  WIC hace que S y S barra iguales} that $p(\mathcal{S}_{\E})=\overline{S_{\E}}$ (see the construction in Subsection \ref{subsec:constru wide C}). Therefore, throughout the proof, we will simply write $f$ instead of $\overline{f}$.
	
	We thus define
	\[
	\Gamma_{X,Y}(\theta):=Q(\theta)Q(u_X)^{-1}.
	\]
	Notice that $\Gamma_{X,Y}$ is a morphism of abelian groups since it is defined in terms of $Q$, which is an additive functor.
	
	\textbf{(3a) $\Gamma_{X,Y}$ is surjective.}
	
	Let $\rho\in\operatorname{Hom}_{\D/\E}(X,Y).$ Then $\rho$ is represented by a diagram of the form
	\[
	\begin{tikzcd}
		& Z \arrow[dl, "s"'] \arrow[dr, "g"] & \\
		X \arrow[
		"\rho"',
		decorate,
		decoration={snake, amplitude=.4mm, segment length=2mm},
		from=2-1,
		to=2-3
		]& & Y
	\end{tikzcd}
	\]
	with $s\in S_{\E}$; that is, $\rho=Q(g)Q(s)^{-1}=[s/g]$.
	
	Since $S_{\E}=\langle\R\cup\Ll\rangle$, it suffices to consider the cases $s\in\Ll$ and $s\in\R$. Suppose first that $s\in\Ll$. By Lemma \ref{lem: para supra de adjunto izq}, we get 
	$$\Hom_{\D}(X_{0},s):
	\Hom_{\D}(X_{0},Z)
	\xrightarrow{\sim}
	\Hom_{\D}(X_{0},X).$$ 
	 Hence, there exists a unique $f:X_{0}\to Z$ such that $sf=u_X.$
	Then  we obtain the following commutative diagram
	\[
	\begin{tikzcd}[ampersand replacement=\&]
		\& {X_{0}} \& \\
		\& Z \\
		X \&\& Y.
		\arrow["f"{pos=0.6}, from=1-2, to=2-2]
		\arrow["{u_{X}=sf}"', from=1-2, to=3-1]
		\arrow["gf", from=1-2, to=3-3]
		\arrow["s", from=2-2, to=3-1]
		\arrow["g"', from=2-2, to=3-3]
	\end{tikzcd}
	\]
	Therefore the following diagram commutes
	\[
	\begin{tikzcd}[ampersand replacement=\&]
		\&\& {X_{0}} \&\& \\
		\& {X_{0}} \&\& Z \\
		X \&\&\&\& Y
		\arrow["{1_{X_{0}}}"', from=1-3, to=2-2]
		\arrow["f"{pos=0.6}, from=1-3, to=2-4]
		\arrow["{u_{X}}"', from=2-2, to=3-1]
		\arrow["gf"', from=2-2, to=3-5]
		\arrow["s"'{pos=0.7}, from=2-4, to=3-1]
		\arrow["g", from=2-4, to=3-5]
	\end{tikzcd}
	\]
	which shows that $\rho=[s/g]=[u_X/gf].$ Hence
	\[
	\Gamma_{X,Y}(gf)
	=
	Q(gf)Q(u_X)^{-1}
	=
	Q(g)Q(s)^{-1}
	=
	\rho.
	\]
	Thus, $\Gamma_{X,Y}$ is surjective when $s\in\Ll$. The case $s\in\R$ is proved analogously, again using Lemma \ref{lem: para supra de adjunto izq}. Hence, $\Gamma_{X,Y}$ is surjective.\\
	
	\textbf{(3b) $\Gamma_{X,Y}$ is injective.}
	
	Let $\theta:X_{0}\to Y$ be a morphism such that
	\[
	\Gamma_{X,Y}(\theta)
	=
	Q(\theta)\circ Q(u_X)^{-1}
	=
	0.
	\]
	Since $Q(u_X)$ is an isomorphism, it follows that $Q(\theta)=0$. The morphism $Q(\theta)$ is represented by the following diagram:
	\[
	\begin{tikzcd}
		& X_{0} \arrow[dl, "1_{X_{0}}"'] \arrow[dr, "\theta"] & \\
		X_{0} && Y.
	\end{tikzcd}
	\]
	
	On the other hand, the zero morphism $0\in\Hom_{\D/\E}(X_{0},Y)$ is represented by the diagram
	\[
	\begin{tikzcd}
		& X_{0} \arrow[dl, "1_{X_{0}}"'] \arrow[dr, "0"] & \\
		X_{0} && Y.
	\end{tikzcd}
	\]
	
	Since $Q(\theta)=0$, there exist morphisms $s':Z\to X_{0}$ and $\beta:Z\to X_{0}$ such that the following diagram commutes
	\[
	\begin{tikzcd}[ampersand replacement=\&]
		\& Z \& \\
		{X_{0}} \&\& {X_{0}} \\
		{X_{0}} \&\& Y
		\arrow["{s'}"', from=1-2, to=2-1]
		\arrow["\beta", from=1-2, to=2-3]
		\arrow["{1_{X_{0}}}"', from=2-1, to=3-1]
		\arrow["\theta"'{pos=0.2}, from=2-1, to=3-3]
		\arrow["{1_{X_{0}}}"'{pos=0.3}, from=2-3, to=3-1]
		\arrow["0", from=2-3, to=3-3]
	\end{tikzcd}
	\]
	and the morphism $s'=1_{X_{0}}\circ s'=\beta$ belongs to $S_{\E}$. In particular, we have $\theta s'=0$ with $s'\in S_{\E}$.
	
	By \cite[Prop. 4.41]{NakaokaOgawaSakai2022}, we know that $S_{\E}$ is a multiplicative system compatible with the extriangulated structure of $\D$. In particular, $S_{\E}$ is a multiplicative system. Thus, from $\theta s'=0$, there exists  $s'':Y\to Y'$ such that $s''\in S_{\E}$ and
	\[
	s''\theta=0.
	\]
	
	By Lemma \ref{lem: para supra de adjunto izq}, we know that 
	$$\Hom_{\D}(X_{0},s''):
	\Hom_{\D}(X_{0},Y)
	\xrightarrow{\sim}
	\Hom_{\D}(X_{0},Y').$$ 
	Since $\Hom_{\D}(X_{0},s'')(\theta)=s''\theta=0$, we conclude that $\theta=0$. Therefore, $\Gamma_{X,Y}$ is injective.
	
	Combining this with $(3a)$, we conclude that $\Gamma_{X,Y}$ is an isomorphism of abelian groups.

	\textbf{(3c) $\Gamma_{X,Y}$ is functorial in the first variable.}
	
	Let $\alpha:X
	\mathrel{\tikz[baseline=-0.5ex]{
			\draw[->,decorate,decoration={snake, amplitude=.4mm, segment length=2mm}]
			(0,0) -- (0.7,0);
	}}X'$ be a morphism in $\D/\E$. Let us show that the following diagram commutes
	\begin{equation}\label{dia: para funtorialidad Phi 1}
		\begin{tikzcd}[ampersand replacement=\&]
			{\D(\overline{\Phi}(X'),Y)} \&\& {\D/\E(X',Y)} \\
			{\D(\overline{\Phi}(X),Y)} \&\& {\D/\E(X,Y)}
			\arrow["{\Gamma_{X',Y}}", from=1-1, to=1-3]
			\arrow["{\D(\overline{\Phi}(\alpha),Y)}"', from=1-1, to=2-1]
			\arrow["{\D/\E(\alpha,Y)}", from=1-3, to=2-3]
			\arrow["{\Gamma_{X,Y}}"', from=2-1, to=2-3]
		\end{tikzcd}
	\end{equation}
	
	Indeed, suppose that $\alpha$ is represented by the following diagram
	\[
	\begin{tikzcd}[ampersand replacement=\&]
		\& Z \& \\
		X \&\& {X',}
		\arrow["{s_{\alpha}}"', from=1-2, to=2-1]
		\arrow["{f_{\alpha}}", from=1-2, to=2-3]
	\end{tikzcd}
	\]
	that is, $\alpha=Q(f_{\alpha})Q(s_{\alpha})^{-1}$. Consequently,
	\[
	\overline{\Phi}(\alpha)
	=
	\overline{\Phi}([s_{\alpha}\di f_{\alpha}])
	=
	\overline{\Phi}\bigl(Q(f_{\alpha})Q(s_{\alpha})^{-1}\bigr)
	=
	\Phi(f_{\alpha})\circ\Phi(s_{\alpha})^{-1}.
	\]
	Now, we recall (see \eqref{diag: funtor torsion 3 con eta}) that 
	$\Phi(h)=\tf(h)$ for any $h\in\mathrm{Mor}(\D).$ Hence, we obtain the following commutative diagram
	\begin{equation}\label{diag: funtor torsion 3 con eta2}
		\begin{tikzcd}[ampersand replacement=\&]
			{X_{0}} \& X \& {i_{*}i^{*}(X)} \& {} \\
			{Z_{0}} \& Z \& {i_{\ast}i^{*}(Z)} \& {} \\
			{X'_{0}} \& {X'} \& {i_{\ast}i^{*}(X')} \& {}
			\arrow["u_{X}", from=1-1, to=1-2]
			\arrow["{\eta_1(X)}", from=1-2, to=1-3]
			\arrow["{\delta_{X}}", dotted, from=1-3, to=1-4]
			\arrow["{\Phi(s_{\alpha})}", from=2-1, to=1-1]
			\arrow["u_{Z}", from=2-1, to=2-2]
			\arrow["{\Phi(f_{\alpha})}"', from=2-1, to=3-1]
			\arrow["{s_{\alpha}}", from=2-2, to=1-2]
			\arrow["{\eta_1(Z)}", from=2-2, to=2-3]
			\arrow["{f_{\alpha}}"', from=2-2, to=3-2]
			\arrow["{i_{\ast}i^{*}(s_{\alpha})}"', from=2-3, to=1-3]
			\arrow["{\delta_{Z}}", dotted, from=2-3, to=2-4]
			\arrow["{i_{\ast}i^{*}(f_{\alpha})}", from=2-3, to=3-3]
			\arrow["u_{X'}", from=3-1, to=3-2]
			\arrow["{\eta_1(X')}"', from=3-2, to=3-3]
			\arrow["{\delta_{X'}}"', dotted, from=3-3, to=3-4]
		\end{tikzcd}
	\end{equation}
	
	Thus, diagram \eqref{dia: para funtorialidad Phi 1} can be rewritten as
	\begin{equation}\label{dia: para funtorialidad Phi 1-1}
		\begin{tikzcd}[ampersand replacement=\&]
			{\D(X_{0}',Y)} \&\& {\D/\E(X',Y)} \\
			{\D(X_{0},Y)} \&\& {\D/\E(X,Y)}
			\arrow["{\Gamma_{X',Y}}", from=1-1, to=1-3]
			\arrow["{\D(\Phi(f_{\alpha})\Phi(s_{\alpha})^{-1},Y)}"', from=1-1, to=2-1]
			\arrow["{\D/\E(Q(f_{\alpha})Q(s_{\alpha})^{-1},Y)}", from=1-3, to=2-3]
			\arrow["{\Gamma_{X,Y}}"', from=2-1, to=2-3]
		\end{tikzcd}
	\end{equation}
	
	We now show that diagram \eqref{dia: para funtorialidad Phi 1-1} commutes. Let $h\in\D(X_{0}',Y)$. On the one hand, we have
	\begin{align*}
		\D/\E(\alpha,Y)\circ\Gamma_{X_{0}',Y}(h)
		&=
		\D/\E(\alpha,Y)\bigl([u_{X'}/h]\bigr)\\
		&=
		\D/\E(\alpha,Y)\bigl(Q(h)Q(u_{X'})^{-1}\bigr)
		\\
		&=
		\bigl(Q(h)Q(u_{X'})^{-1}\bigr)\circ\alpha\\
		&=
		Q(h)Q(u_{X'})^{-1}Q(f_{\alpha})Q(s_{\alpha})^{-1}.
	\end{align*}
	
	On the other hand, we know that
	\begin{align*}
		\Gamma_{X,Y}\circ\D(\overline{\Phi}(\alpha),Y)(h)
		&=
		\Gamma_{X,Y}\bigl(h\circ\overline{\Phi}(\alpha)\bigr)\\
		&=
		\Gamma_{X,Y}\bigl(h\circ\Phi(f_{\alpha})\circ\Phi(s_{\alpha})^{-1}\bigr)\\
		&=
		\bigl[u_X\,/\,(h\circ\Phi(f_{\alpha})\circ\Phi(s_{\alpha})^{-1})\bigr]\\
		&=
		Q(h)Q(\Phi(f_{\alpha}))Q(\Phi(s_{\alpha})^{-1})Q(u_X)^{-1}.
	\end{align*}
	
	From \eqref{diag: funtor torsion 3 con eta2}, we have
	$Q(u_X)Q(\Phi(s_{\alpha}))
	=
	Q(s_{\alpha})Q(u_Z)$ and $Q(u_{X'})Q(\Phi(f_{\alpha}))
	=
	Q(f_{\alpha})Q(u_Z)$.
	Therefore, $Q(\Phi(s_{\alpha})^{-1})
	=
	Q(u_Z)^{-1}Q(s_{\alpha})^{-1}Q(u_X)$
	and $Q(\Phi(f_{\alpha}))
	=
	Q(u_{X'})^{-1}Q(f_{\alpha})Q(u_Z)$.
	Consequently,
	\begin{align*}
		\Gamma_{X,Y}\circ \D(\overline{\Phi}&(\alpha),Y)(h)
		=
		Q(h)Q(\Phi(f_{\alpha}))Q(\Phi(s_{\alpha})^{-1})Q(u_X)^{-1}\\
		&=
		Q(h)\bigl(Q(u_{X'})^{-1}Q(f_{\alpha})Q(u_Z)\bigr)
		\bigl(Q(u_Z)^{-1}Q(s_{\alpha})^{-1}Q(u_X)\bigr)
		Q(u_X)^{-1}\\
		&=
		Q(h)Q(u_{X'})^{-1}Q(f_{\alpha})Q(s_{\alpha})^{-1}.
	\end{align*}
	Thus $\D/\E(\alpha,Y)\circ\Gamma_{X',Y}(h)
	=
	\Gamma_{X,Y}\circ\D(\overline{\Phi}(\alpha),Y)(h)$ and hence diagram \eqref{dia: para funtorialidad Phi 1-1} commutes. Therefore, diagram \eqref{dia: para funtorialidad Phi 1} also commutes.

		\textbf{(3d) $\Gamma_{X,Y}$ is functorial in the second variable.}\\
		Let $\beta: Y \to Y'$ be a morphism in $\D$. Let us show that the following diagram commutes
		\begin{equation}\label{dia: funtorial Gamma 2}
			\begin{tikzcd} 
				{\D(X, Y)} & & {\D / \E (X,Q(Y))} \\
				{\D(X, Y')} & & {\D / \E (X,Q(Y'))}.
				\arrow["{{\Gamma_{X,Y}}}", from=1-1, to=1-3]
				\arrow["{{\D(X, \beta)}}"', from=1-1, to=2-1] 
				\arrow["{{\D / \E (X,Q(\beta))}}", from=1-3, to=2-3] 
				\arrow["{{\Gamma_{X,Y'}}}"', from=2-1, to=2-3] 
			\end{tikzcd}
		\end{equation}
		
		Indeed, let $h \in \D(X,Y)$. On the one hand,  we have
		$$\begin{aligned}
			\Gamma_{X,Y'} \circ \D(X, \beta)(h)
			&= \Gamma_{X,Y'}(\beta h)\\
			&= Q(\beta h)Q(u_{X})^{-1}.
		\end{aligned}$$
		
		On the other hand,
	$$	\begin{aligned}
			\D/\E(X,Q(\beta)) \circ \Gamma_{X,Y}(h)
			&= \D/\E(X,Q(\beta))(Q(h)Q(u_{X})^{-1})\\
			&= Q(\beta)Q(h)Q(u_{X})^{-1}\\
			&= Q(\beta h)Q(u_{X})^{-1}.
		\end{aligned}$$
		Thus, the diagram \eqref{dia: funtorial Gamma 2} commutes.
		
		\textbf{(Step 4) $j_!:=\overline{\Phi}$ defined in (Step 2) is full and faithful.}\\
		By Proposition \ref{pro: F fiel y pleno  sii unidad iso} (b), it is enough to show that the unit $\eta: 1_{\D / \E} \to Q\circ\overline{\Phi}$ of the adjunction $(\overline{\Phi},Q)$ is an isomorphism.
		
		Consider the functorial isomorphism constructed above
		$$\Gamma_{X,\overline{\Phi}(X)}:
		\Ho_{\D}(\overline{\Phi}(X),\overline{\Phi}(X))
		\to
		\Ho_{\D/\E}\left(X,(Q\circ\overline{\Phi})(X)\right).$$
		The unit of the adjunction $(\overline{\Phi},Q)$ is defined by $\eta_{X}:=\Gamma_{X,\overline{\Phi}(X)}(1_{\overline{\Phi}(X)})$ \cite[Thm. 3.1.5]{Borceux1}, p.~99). Notice that $Q(\overline{\Phi}(X))=Q(X_{0})=X_{0}$. Therefore, by construction, $\Gamma_{X,\overline{\Phi}(X)}(1_{\overline{\Phi}(X)})$ is given by the fraction $[1_{X_{0}}\di u_{X}]$, represented by the following diagram
		$$\begin{tikzcd}[ampersand replacement=\&]
			\& X_{0} \& \\
			X \&\& X_{0}.
			\arrow["u_{X}"', from=1-2, to=2-1]
			\arrow["1_{X_{0}}", from=1-2, to=2-3]
		\end{tikzcd}$$
		
		Since $u_{X} \in S(\E)$, it follows that $Q(u_{X})=[1_{X_{0}}\di u_{X}]=Q(u_{X})$ is invertible. Thus
		$$
		\eta_X=\Gamma_{X,\overline{\Phi}(X)}(1_{\overline{\Phi}(X)})
		=Q(1_{X_{0}})Q(u_{X})^{-1}
		=Q(u_{X})^{-1}
		$$
		is an isomorphism in $\D/\E$. Consequently, $\eta$ is an isomorphism and therefore $\overline{\Phi}$ is full and faithful.	
		
			\textbf{(Step 5) The transformations $\eta: 1_{\D / \E}\to j^*\,j_!$ and
				$\varepsilon: j_!\, j^* \to 1_{\D}$ are natural extriangulated
				transformations.} 
				 Indeed, by \cite[Theorem 3.1.5]{Borceux1},
				the morphisms $\eta_{X}:=\Gamma_{X,\overline{\Phi}(X)}
				(1_{\overline{\Phi}(X)})$ and
				$\varepsilon_{Y}:= \Gamma_{QY,Y}^{-1}(1_{Q(Y)})$
				define natural transformations. In particular, by construction,
				$\eta_{X}=Q(u_{X})^{-1}$ and
				$\varepsilon_{Y}=1_{Y}\, u_{Y}=u_{Y}$.
				
				Since $(j^*:=Q, \psi^1)$ and
				$(j_!:=\overline{\Phi}, \psi^2)$ are extriangulated functors,
				their composition $(Q\, \overline{\Phi}, \lambda^1)$ is also an
				extriangulated functor, where, for each $X,Y \in \D/\E$,
				the morphism
				$\lambda^1_{Y,X}: \mathbb{G}(Y,X) \to
				\mathbb{G}(\overline{\Phi} \, Q Y, \, \overline{\Phi} \, Q X)$
				is defined as follows.

		\begin{equation}\label{eq: lamb 1}
			\begin{aligned}
				\lambda^1_{Y,X}([\overline{t}\di\overline{\delta}/\overline{s}])
				&= \psi^1_{\overline{\Phi Y}, \overline{\Phi}X}
				\left(\overline{\psi}^2([\overline{t}\di\overline{\delta}/\overline{s}])\right)\\
				&= \psi^1_{Y_{0}, X_{0}}
				\left(\overline{\psi}^2([\overline{t}\di\overline{\delta}/\overline{s}])\right)\\
				&= \psi^1_{Y_{0}, X_{0}}
				\left(\Phi(s)^{-1}\cdot
				\overline{\psi^{2}}(\delta)\cdot\Phi(t)^{-1}\right)\\
				& = [\overline{1_{X_0}}\di \Phi(s)^{-1}\cdot
				\overline{\psi^{2}}(\delta)\cdot\Phi(t)^{-1} / \overline{1_{Y_0}}]\\
				& = [\overline{1_{X_0}}\di (i\, \tf(s))^{-1}\cdot
				\overline{\psi^{2}}(\delta)\cdot (i\, \tf(t))^{-1} / \overline{1_{Y_0}}]\\
				& = [\overline{1_{X_0}}\di \tf(s)^{-1}\cdot
				\overline{\psi^{2}}(\delta)\cdot  \tf(t)^{-1} / \overline{1_{Y_0}}] & (see \  Step  \ 1)\\
				& = [\overline{1_{X_0}}\di \tf(s)^{-1}\cdot
				\mu^1_{\tf C, \tf A} \, \mu^2_{C,A}(\delta)\cdot  \tf(t)^{-1} / \overline{1_{Y_0}}] & ( \mu^1_{\tf C, \tf A} \  inclusion)\\
				& = [\overline{1_{X_0}}\di \tf(s)^{-1}\cdot \mu^2_{C,A}(\delta)\cdot  \tf(t)^{-1} / \overline{1_{Y_0}}]
			\end{aligned}
		\end{equation}
		
		for every element $[\overline{t}\di\overline{\delta}/\overline{s}]$
		of $\mathbb{G}(Y,X)$. Similarly, $(1_{\D/ \E}, id):
		\D/\E \to \D/\E$ is an extriangulated functor. We now show that
		$\eta: (1_{\D / \E}, id) \to (j^*\,j_!, \lambda^1)$
		is a natural extriangulated transformation in the sense of
		\cite[Definition 2.11]{NakaokaOgawaSakai2022}.\\
	\textit{Claim:} The following diagram commutes.
		\begin{equation}\label{dia: eta-nat}
			\begin{tikzcd}[ampersand replacement=\&]
				{\mathbb{G}(Y,X)} \& {\mathbb{G}(Y,X)} \\
				{\mathbb{G}(Q \overline{\Phi}Y, Q \overline{\Phi}X)} \& {\mathbb{G}(Y, Q \overline{\Phi}X)}
				\arrow["{id_{Y,X}}", from=1-1, to=1-2]
				\arrow["{\lambda^1_{Y,X}}"', from=1-1, to=2-1]
				\arrow["{\mathbb{G}(Y, \eta_{X})}", from=1-2, to=2-2]
				\arrow["{\mathbb{G}(\eta_{Y}, Q \overline{\Phi}X)}"', from=2-1, to=2-2]
			\end{tikzcd}
		\end{equation}

		Indeed, let $[\overline{t}\di\overline{\delta}/\overline{s}]$ be an element of $\mathbb{G}(Y,X)$. On the one hand,
		
	 \begin{align*}
		 	\mathbb{G}(Y, \eta_{X}) \, id_{Y,X} ([\overline{t}\di\overline{\delta}/\overline{s}])& = \mathbb{G}(Y, \eta_X)([\overline{t}\di\overline{\delta}/\overline{s}])\\
		 	& =[\overline{t}\di\overline{\delta}/\overline{s}] \cdot \eta_{X}\\
		 	& = [\overline{t}\di\overline{\delta}/\overline{s}] \cdot Q(u_{X})^{-1}\\
		 	& = [\overline{t}\di\overline{\delta}/\overline{s}] \cdot [\overline{1_{X_{0}}} / \overline{u_{X}}].
		 \end{align*}

		 The last action is computed according to \cite[Definition 3.23]{NakaokaOgawaSakai2022}, using the right fraction $[\overline{u_{X}}/ \overline{1_{X}}]$, which is equivalent to the left fraction $[\overline{1_{X_{0}}} / \overline{u_{X}}]$ via the following diagram.

		 \[\begin{tikzcd}[ampersand replacement=\&]
		 	{X_{0}} \& {X_{0}} \\
		 	X \& X
		 	\arrow["{1_{X_{0}}}", from=1-1, to=1-2]
		 	\arrow["{u_{X}}"', from=1-1, to=2-1]
		 	\arrow["{u_{X}}", from=1-2, to=2-2]
		 	\arrow["{1_{X}}"', from=2-1, to=2-2]
		 \end{tikzcd}\]
		 
		 obtained from condition $(M2)$ of $\mathcal{S}_{\E}$. To compute $[\overline{t}\di\overline{\delta}/\overline{s}] \cdot [\overline{u_{X}}/ \overline{1_{X}}]$, it suffices to consider the following diagram 
		 \[\begin{tikzcd}[ampersand replacement=\&]
		 	X \& X \&\&\& \\
		 	A \& A \& B \& C \& {} \\
		 	\&\&\& Y
		 	\arrow["s"', from=1-1, to=2-1]
		 	\arrow["{1_{X}}"', from=1-2, to=1-1]
		 	\arrow["s", from=1-2, to=2-2]
		 	\arrow["{1_{A}}", from=2-2, to=2-1]
		 	\arrow["a", from=2-2, to=2-3]
		 	\arrow["b", from=2-3, to=2-4]
		 	\arrow["\delta", dashed, from=2-4, to=2-5]
		 	\arrow["t", from=2-4, to=3-4]
		 \end{tikzcd}\]
		 
		 and so,
		 $$\begin{aligned}
		 	[\overline{t}\di\overline{\delta}/\overline{s}] \cdot [\overline{1_{X_{0}}} / \overline{u_{X}}] & = [\overline{t}\di\overline{\delta}/\overline{s}] \cdot [\overline{u_{X}}/ \overline{1_{X}}]\\
		 	& = [\overline{t}\di\overline{1_{A} \cdot \delta}/\overline{s \, u_{X}}] \\
		 	& = [\overline{t}\di\overline{\delta}/\overline{s\, u_{X}}].
		 \end{aligned}$$

		 On the other hand, 
		 \begin{equation}\label{eq: eta-extr-lad2}
		 	\begin{aligned}
		 	\mathbb{G}(\eta_{Y},Q\, \overline{\Phi} X)\, \lambda^1_{Y,X}([\overline{t}\di\overline{\delta}/\overline{s}]) & = \mathbb{G}(\eta_{Y},Q\, \overline{\Phi} X) \left( [\overline{1_{X_0}}\di \Phi(s)^{-1}\cdot
		 	\overline{\psi^{2}}(\delta)\cdot\Phi(t)^{-1} / \overline{1_{Y_0}}] \right)\\
		 	& = \mathbb{G}(\eta_{Y},Q\, \overline{\Phi} X) \left( [\overline{1_{X_0}}\di \Phi(s)^{-1}\cdot
		 	\overline{\psi^{2}}(\delta)\cdot\Phi(t)^{-1} / \overline{1_{Y_0}}] \right)\\
		 	& = \mathbb{G}(\eta_{Y},Q\, \overline{\Phi} X) [\overline{1_{X_0}}\di \tf(s)^{-1}\cdot \mu^2_{C,A}(\delta)\cdot  \tf(t)^{-1} / \overline{1_{Y_0}}]\\
		 	& = [\overline{1_{X_0}}\di \tf(s)^{-1}\cdot \mu^2_{C,A}(\delta)\cdot  \tf(t)^{-1} / \overline{1_{Y_0}}] \cdot \eta_{Y}\\
		 	& = [\overline{1_{X_0}}\di \tf(s)^{-1}\cdot \mu^2_{C,A}(\delta)\cdot  \tf(t)^{-1} / \overline{1_{Y_0}}] \cdot [\overline{1_{Y_0}}/ \overline{u_Y}].
		 \end{aligned}
		\end{equation} 
		
		We now consider the following canonical diagram of conflations (see Proposition \ref{pro:torsion_funtor_extriangulado})
		 \begin{equation}\label{dia: nat-ext-eta}
		 \begin{tikzcd}[ampersand replacement=\&]
		 	\&\&\&\&\& {} \\
		 	{X_{0}} \& {A_{0}} \& {B_{0}} \& {C_{0}} \& {Y_{0}} \& {} \\
		 	X \& A \& B \& C \& Y \& {} \\
		 	{i^*i_*X} \& {i^*i_*A} \& {i^*i_*B} \& {i^*i_*C} \& {i^*i_*Y.} \\
		 	{} \& {} \& {} \& {} \& {}
		 	\arrow["{\tf(s)}", from=2-1, to=2-2]
		 	\arrow["{u_{X}}"', from=2-1, to=3-1]
		 	\arrow["{\tf(a)}", from=2-2, to=2-3]
		 	\arrow["{u_{A}}"', from=2-2, to=3-2]
		 	\arrow["{\tf(b)}", from=2-3, to=2-4]
		 	\arrow["{u_{B}}"', from=2-3, to=3-3]
		 	\arrow["{\mu_{C,A}^2(\delta)}"{pos=0.3}, dashed, shorten <=2pt, shorten >=60pt, from=2-4, to=1-6]
		 	\arrow["{\tf(t)}"', from=2-4, to=2-5]
		 	\arrow["{u_{C}}"', from=2-4, to=3-4]
		 	\arrow["{u_{Y}}", from=2-5, to=3-5]
		 	\arrow["s"', from=3-1, to=3-2]
		 	\arrow["{\eta_{X}}"', from=3-1, to=4-1]
		 	\arrow["a"', from=3-2, to=3-3]
		 	\arrow["{\eta_{A}}"', from=3-2, to=4-2]
		 	\arrow["b"', from=3-3, to=3-4]
		 	\arrow["{\eta_{B}}"', from=3-3, to=4-3]
		  \arrow["\delta"{pos=0.2}, dashed, shorten <=2pt, shorten >=60pt, from=3-4, to=2-6]
		 	\arrow["t"{description, pos=0.6}, from=3-4, to=3-5]
		 	\arrow["{\eta_{C}}"', from=3-4, to=4-4]
		 	\arrow["{\eta_{Y}}", from=3-5, to=4-5]
		 	\arrow["{i^*i_*s}"', from=4-1, to=4-2]
		 	\arrow["{\delta_{X}}"', dashed, from=4-1, to=5-1]
		 	\arrow["{i^*i_*a}"', from=4-2, to=4-3]
		 	\arrow["{\delta_{A}}"', dashed, from=4-2, to=5-2]
		 	\arrow["{i^*i_*b}"', from=4-3, to=4-4]
		 	\arrow["{\delta_{B}}"', dashed, from=4-3, to=5-3]
		 	\arrow["{\mu^3_{C,A}(\delta)}"{pos=0.2}, dashed, shorten <=2pt, shorten >=60pt, from=4-4, to=3-6]
		 	\arrow["{i^*i_*t}"', from=4-4, to=4-5]
		 	\arrow["{\delta_{C}}"', dashed, from=4-4, to=5-4]
		 	\arrow["{\delta_{Y}}", dashed, from=4-5, to=5-5]
		 \end{tikzcd}
		\end{equation}

Since $t,s \in \mathcal{S}_{\E}$, it follows that
$\overline{\Phi}(t)$ and $\overline{\Phi}(s)$ are isomorphisms in
$\D/\E$. By the definition of $\overline{\Phi}$, we conclude that
$\tf(s)^{-1}: A_{0}\to X_{0}$ and
$\tf(t)^{-1}: Y_{0}\to C_{0}$ are isomorphisms.
Hence, by \cite[Proposition 3.7]{NP19}, the element
$\tf(s)^{-1} \cdot \mu^2_{C,A}(\delta) \cdot \tf(t)^{-1}$
is realized by the following conflation
		\[\begin{tikzcd}[ampersand replacement=\&]
			{X_{0}} \&\& {B_{0}} \&\& {Y_{0}} \& {}
			\arrow["{\tf(a)\tf(s)}", from=1-1, to=1-3]
			\arrow["{\tf(t)\tf(b)}", from=1-3, to=1-5]
			\arrow["{\tf(s)^{-1} \cdot \mu^2_{C,A}(\delta) \cdot \tf(t)^{-1}}"{pos=1}, dashed, from=1-5, to=1-6]
		\end{tikzcd}\]
		With the above in place, we perform the action
		$[\overline{1_{X_0}}\di \tf(s)^{-1}\cdot \mu^2_{C,A}(\delta)\cdot  \tf(t)^{-1} / \overline{1_{Y_0}}] \cdot [\overline{1_{Y_0}}/ \overline{u_Y}]$,
		according to \cite[Definition 3.23]{NakaokaOgawaSakai2022}, using the left fraction $[\overline{1_{Y_{0}}} / \overline{u_{Y}}]$. This yields the following diagram
		
		\[\begin{tikzcd}[ampersand replacement=\&]
			{X_0} \&\&\&\&\&\& {} \\
			{X_{0}} \&\& {B_{0}} \&\& {Y_{0}} \& {Y_{0}} \\
			\&\&\&\& {Y_{0}} \& {Y_{0}}
			\arrow["{1_{X_0}}"'{pos=0.4}, from=1-1, to=2-1]
			\arrow["{\tf(a)\tf(s)}", from=2-1, to=2-3]
			\arrow["{\tf(t)\tf(b)}", from=2-3, to=2-5]
			\arrow["{\tf(s)^{-1} \cdot \mu^2_{C,A}(\delta) \cdot \tf(t)^{-1}}"{pos=0.7}, dashed, shorten <=2pt, shorten >=40pt, from=2-5, to=1-7]
			\arrow["{1_{Y_0}}"', from=2-5, to=3-5]
			\arrow["{1_{Y_0}}", from=2-6, to=2-5]
			\arrow["{1_{Y_0}}", from=2-6, to=3-6]
			\arrow["{1_{Y_0}}", from=3-6, to=3-5]
		\end{tikzcd}\]
		
		Consequently, we have
		$[\overline{1_{X_0}}\di \tf(s)^{-1}\cdot \mu^2_{C,A}(\delta)\cdot  \tf(t)^{-1} / \overline{1_{Y_0}}] \cdot [\overline{1_{Y_0}}/ \overline{u_Y}]= [\overline{1_{X_0}}\di \tf(s)^{-1}\cdot \mu^2_{C,A}(\delta)\cdot  \tf(t)^{-1} / \overline{u_Y}]$.
		Therefore, it follows from \eqref{eq: eta-extr-lad2} that
		$$
		\mathbb{G}(\eta_{Y},Q\, \overline{\Phi} X)\, \lambda^1_{Y,X}([\overline{t}\di\overline{\delta}/\overline{s}])= [\overline{1_{X_0}}\di \tf(s)^{-1}\cdot \mu^2_{C,A}(\delta)\cdot  \tf(t)^{-1} / \overline{u_Y}].
		$$
		
		We now show that diagram \eqref{dia: eta-nat} commutes. To this end, it suffices to verify that
		$[\overline{1_{X_0}}\di \tf(s)^{-1}\cdot \mu^2_{C,A}(\delta)\cdot  \tf(t)^{-1} / \overline{u_Y}]= [\overline{s \, u_{X}}\di \delta / \overline{t}]$,
		in the sense of \cite[Proposition 3.9]{NakaokaOgawaSakai2022}. Indeed, we have the following diagram
		 \[\begin{tikzcd}[ampersand replacement=\&]
		 	\&\& {Y_{0}} \&\& {X_{0}} \&\& \\
		 	Y \&\& {Y_{0}} \&\& A \&\& {X_0} \\
		 	\&\& C \&\& A
		 	\arrow["{\tf(s)^{-1}\cdot \mu^2_{C,A}(\delta)\cdot  \tf(t)^{-1}}", dashed, from=1-3, to=1-5]
		 	\arrow["{u_{Y}}"', from=1-3, to=2-1]
		 	\arrow["{s\,u_X}", from=1-5, to=2-5]
		 	\arrow["{1_{Y_0}}", from=2-3, to=1-3]
		 	\arrow["{u_{Y}}"', from=2-3, to=2-1]
		 	\arrow[dashed, from=2-3, to=2-5]
		 	\arrow["{u_{C}\tf(t)^{-1}}"'{pos=0.2}, from=2-3, to=3-3]
		 	\arrow["{1_{X_0}}"', from=2-7, to=1-5]
		 	\arrow["{s\,u_X}"', from=2-7, to=2-5]
		 	\arrow["{s\, u_X}", from=2-7, to=3-5]
		 	\arrow["t"{pos=0.4}, from=3-3, to=2-1]
		 	\arrow["\delta"', dashed, from=3-3, to=3-5]
		 	\arrow["{1_{A}}"', from=3-5, to=2-5]
		 \end{tikzcd}\]	 
		 from which it follows that
		 $1_{A}\cdot \delta \cdot (u_{C}\,\tf(t)^{-1})=
		 \delta \cdot (u_{C}\, \tf(t)^{-1})$.
		 Now, using diagram \eqref{dia: nat-ext-eta}, we obtain
		 
		  \begin{align*}
		 	(s\, u_{X}) \cdot
		 	\left(\tf(s)^{-1}\cdot \mu^2_{C,A}(\delta)\cdot \tf(t)^{-1}\right)
		 	\cdot 1_{Y_{0}}
		 	&= (s\, u_{X}) \cdot
		 	\left(\tf(s)^{-1}\cdot \mu^2_{C,A}(\delta)\cdot \tf(t)^{-1}\right)\\
		 	&= (u_{A}\, \tf(s)) \cdot
		 	\left(\tf(s)^{-1}\cdot \mu^2_{C,A}(\delta)\cdot \tf(t)^{-1}\right)\\
		 	&= u_{A} \cdot \left(
		 	\tf(s) \cdot
		 	\left(\tf(s)^{-1}\cdot \mu^2_{C,A}(\delta)\cdot \tf(t)^{-1}\right)
		 	\right)\\
		 	&= u_{A} \cdot \left(
		 	(\tf(s)\,\tf(s)^{-1})\cdot \mu^2_{C,A}(\delta)\cdot \tf(t)^{-1}
		 	\right)\\
		 	&= u_{A} \cdot \left(
		 	1_{C_0}\cdot \mu^2_{C,A}(\delta)\cdot \tf(t)^{-1}
		 	\right)\\
		 	&= u_{A} \cdot \left(
		 	\mu^2_{C,A}(\delta)\cdot \tf(t)^{-1}
		 	\right)\\
		 	&= (u_{A} \cdot \mu^2_{C,A}(\delta))\cdot \tf(t)^{-1}\\
		 	&= (\delta \cdot u_{C}) \cdot \tf(t)^{-1}\\
		 	&= \delta \cdot (u_{C}\,\tf(t)^{-1})\\
		 	&= 1_{A} \cdot \delta \cdot (u_{C}\,\tf(t)^{-1}).
		 \end{align*}
		 
		 This proves that $[\overline{1_{X_0}}\di \tf(s)^{-1}\cdot \mu^2_{C,A}(\delta)\cdot  \tf(t)^{-1} / \overline{u_Y}]= [\overline{s \, u_{X}}\di \delta / \overline{t}]$, and hence \eqref{dia: eta-nat} commutes. Therefore, $\eta: (Id_{\D/ \E}, id) \to (Q\, \overline{\Phi}, \lambda^1)$ is a natural extriangulated transformation.\\
		 
		 We now show that $\varepsilon\colon(\overline{\Phi}\, Q, \lam^{2})
		 \to (\operatorname{Id}_{\B},\operatorname{id})$ is an extriangulated
		 natural transformation; that is, the following diagram commutes

		  \begin{equation}\label{dia: eps-extr}
		  	\begin{tikzcd}[ampersand replacement=\&]
		  		{\Ee(C,A)} \&\& {\Ee(\overline{\Phi}\, QC, \overline{\Phi}\, QA)} \\
		  		{\Ee(C,A)} \&\& {\Ee(\overline{\Phi}\, QC, A)}
		  		\arrow["{\lambda^2_{C,A}}", from=1-1, to=1-3]
		  		\arrow["id"', from=1-1, to=2-1]
		  		\arrow["{\Ee(\overline{\Phi}\, QC, \varepsilon_{A})}", from=1-3, to=2-3]
		  		\arrow["{\Ee(\varepsilon_{C}, A)}"', from=2-1, to=2-3]
		  	\end{tikzcd}
		  \end{equation}
		  where $\lambda^{2}\colon \Ee(C,A)\to
		  \Ee(\overline{\Phi}\, QC,\overline{\Phi}\, QA)$ is defined as the
		  composition
		  $\psi^{2}\circ\bigl(\psi^{1}\blacktriangle(Q^{op}\times Q)\bigr)$.
		  Indeed, let $\delta\in\Ee(C,A)$. On the one hand,
		  \begin{align*}
		  	\Ee( \overline{\Phi}\, QC, \varepsilon_{A})\, \lambda^{2}_{C,A}(\delta)& = \Ee(C_{0}, \varepsilon_{A})(\lambda_{C,A}^{2}(\delta))\\
		  	& = \Ee(C_{0}, \varepsilon_{A})\left( \psi^2_{QC,QA} \psi^1_{C,A}(\delta) \right)\\
		  	& = \Ee(C_{0}, \varepsilon_{A}) \left(\psi^2_{QC,QA} ([\overline{1_{C}}\di \delta / \overline{1_{A}}]) \right)\\
		  	& = \Ee(C_{0}, \varepsilon_{A})\left(  \Phi(1_{A})^{-1} \cdot \overline{\Phi}^{2}_{C,A}(\delta) \cdot \Phi(1_{C})^{-1}\right)\\
		  	& = \Ee(C_{0}, \varepsilon_{A})\left( \overline{\Phi}^{2}_{C,A}(\delta) \right)\\
		  	& = \Ee(C_{0}, \varepsilon_{A})\left( \mu^{1}_{\tf C,\tf A} \mu^2_{C,A}(\delta) \right)\\
		  	& = \Ee(C_{0}, \varepsilon_{A})\left( \mu^{1}_{\tf C,\tf A} \mu^2_{C,A}(\delta) \right)\\
		  	& = \Ee(C_{0}, \varepsilon_{A})\left( \mu^{1}_{\tf C,\tf A} \mu^2_{C,A}(\delta) \right)\\
		  	& = \Ee(C_{0}, \varepsilon_{A})\left( \mu^{1}_{C_{0}, A_{0}} \mu^2_{C,A}(\delta) \right)\\
		  	& = \Ee(C_{0}, \varepsilon_{A})\left( \mu^2_{C,A}(\delta) \right)\\
		  	& = \varepsilon_{A} \cdot \mu^2_{C,A}(\delta) &\hspace{-2cm}( \text{since }   \varepsilon_{A}= \Gamma^{-1}_{QA,A}(1_{QA}))\\
		  	& = u_{A} \cdot \mu^2_{C,A}(\delta).
		  \end{align*}
		  
		  On the other hand, 
		  \begin{align*}
		   \Ee(\varepsilon_{C}, A) \, id_{C,A}(\delta) & = \Ee(\varepsilon_{C}, A)(\delta)\\
		   & = \delta \cdot \varepsilon_{C}\\
		   & = \delta \cdot u_{C} & (\text{see \ref{dia: nat-ext-eta} } )\\
		   & = u_{A} \cdot \mu_{C,A}^{2}(\delta).
		  \end{align*}
		  Thus, diagram \eqref{dia: eps-extr} commutes, proving that
		  $\varepsilon\colon(\overline{\Phi}\, Q,\lam^{2})
		  \to(\operatorname{Id}_{\B},\operatorname{id})$
		  is an extriangulated natural transformation.

		\textbf{(Step 6) Construction of the right adjoint $j_*$ of $j^*:=Q$.}
		The proof is similar to the one given  for the left adjoint of $Q.$
		
		\textbf{(6a) There exists an extriangulated functor $(\Lambda,\overline{\psi^3}):\D \to \D$.}
		
		 Consider the co-unite $\varepsilon_{2}: i_{*}i^{!} \to 1_{\D}$ of the adjunction $(i_{*},i^{!}).$ Let $X\in\D.$ By (\ref{diag: sucesion de torsion -1})  there exists a conflation
		\begin{equation}\label{diag: sucesion de torsion 2-1}
			\begin{tikzcd}
				i_{*}i^! X & X & X'_{0} & {}
				\arrow["{\varepsilon_2(X)}",from=1-1, to=1-2]
				\arrow[ from=1-2, to=1-3]
				\arrow["{\delta'_{X}}", dotted, from=1-3, to=1-4]
			\end{tikzcd}
		\end{equation}
		such that $X'_{0}\in\operatorname{Ker}(i^{!}).$ Consequently, \eqref{diag: sucesion de torsion 2-1} is the canonical conflation associated with the object $X$ via the s-torsion pair $(\operatorname{Im}(i_{*}),\operatorname{Ker}(i^{!})).$ Thus $\lf(X)=X'_0$ where $\lf: \D \to  \operatorname{Ker}(i^{!})$ is the torsion-free functor given by $(\operatorname{Im}(i_{*}),\operatorname{Ker}(i^{!})),$ see 
Remark \ref{Rk-inc-ext-tp} (2).
		
		Consider the inclusion functor $i: \operatorname{Ker}(i^{!}) \to \D$ and the extriangulated functor $(i^{!},\varphi^{3}):(\mathcal{D},\mathbb{E},\mathfrak{s}) \to (\mathcal{D}',\mathbb{F},\mathfrak{r}).$ Moreover, by Remark \ref{Rk-inc-ext-tp} (1), we have the extriangulated functor
		\[
		(i,\nu_{1}):(\operatorname{Ker}(i^!),\Ee|_{\operatorname{Ker}(i^!)},\s|_{\operatorname{Ker}(i^!)}) \longrightarrow (\D,\Ee,\s)
		\]
where $\nu_1$ is the identity. We set $\Lambda:=i\circ\lf:\D\to \D.$
		Explicitly, for $f:X\to Y$ in $\D$, we have 
		$
		\Lambda(f)=i(\lf(f))=\lf(f),
		$
		where $\lf(f)$ is the unique morphism making the following diagram commute
		\begin{equation}\label{diag: funtor torsion 3 con varepsilon}
			\begin{tikzcd}[ampersand replacement=\&]
				{i_{\ast}i^{!}(X)} \& X \& {X'_{0}=\lf(X)} \& {} \\
				{i_{\ast}i^{!}(Y)} \& Y \& {Y'_{0}=\lf(Y)} \& {}
				\arrow["{\varepsilon_2(X)}", from=1-1, to=1-2]
				\arrow["{{i_{\ast}i^{!}(f)}}"', from=1-1, to=2-1]
				\arrow[from=1-2, to=1-3]
				\arrow["f"', from=1-2, to=2-2]
				\arrow["{{\delta'_{X}}}", dotted, from=1-3, to=1-4]
				\arrow["{\lf(f)}", from=1-3, to=2-3]
				\arrow["{\varepsilon_2(Y)}"', from=2-1, to=2-2]
				\arrow[from=2-2, to=2-3]
				\arrow["{{\delta'_{Y}}}"', dotted, from=2-3, to=2-4]
			\end{tikzcd}
		\end{equation}
		
		Since $(i^{!},\varphi^{3})$ and $(i_{*},\varphi^{1})$ are extriangulated functors, their composition $(i_{*}i^{!},\nu_{3})$ is also an extriangulated functor. The natural transformation $\nu_{3}$ is given by the following composition of natural transformations:
		\[
		\begin{tikzcd}[ampersand replacement=\&]
			{\mathbb{E}} \& {\mathbb{F}\circ({i^{!}}^{op} \times i^{!})} \&\& {(\mathbb{E}\circ({i_{*}}^{op} \times i_{*})) \circ ({i^{!}}^{op} \times i^{!})}
			\arrow["{\varphi^3}", Rightarrow, from=1-1, to=1-2]
			\arrow["{\varphi^1 \blacktriangle ({i^{!}}^{op} \times i^{!})}", Rightarrow, from=1-2, to=1-4]
		\end{tikzcd}
		\]
		where $\blacktriangle$ denotes whiskering, that is,
		\[
		\nu^3_{C,A}
		=
		\left[
		\left(\varphi^1 \blacktriangle ({i^{!}}^{op} \times i^{!})\right)
		\circ \varphi^3
		\right]_{C,A}
		=
		\left(\varphi^1 \blacktriangle ({i^{!}}^{op} \times i^{!})\right)_{C,A}
		\circ \varphi^3_{C,A}
		=
		\varphi^{1}_{i^!(C),i^!(A)}\circ\varphi^3_{C,A}
		\]
		for all $A,C\in\D$.
		
		Now, since $(i_{*}i^!,\nu^{3}),\;
		(1_{\D},1_{\Ee}):(\D,\Ee,\s)\to (\D,\Ee,\s)$ are extriangulated functors,\footnote{Notice that $i_{*}i^!$ corresponds to the composition of the torsion functor $\tf:\D\to \operatorname{Im}(i_{*})$ with the inclusion functor of $\operatorname{Im}(i_{*})$ into $\D$.} Proposition \ref{pro:torsion_funtor_extriangulado} implies that the torsion-free functor
		\[
		(\lf,\nu^{2}):(\D,\Ee,\s)\to
		(\operatorname{Ker}(i^{!}),\Ee_{\mid_{\operatorname{Ker}(i^{!})}},\s_{\mid_{\operatorname{Ker}(i^{!})}})
		\]
		is an extriangulated functor. 
		
		Consequently, $\Lambda$ is the composition of $(\lf,\nu^2)$ and $(i,\nu^1)$, which are both extriangulated functors. Therefore, $\Lambda$ is also an extriangulated functor, with associated natural transformation $\overline{\psi^{3}}:\Ee\to\Ee(\Lambda^{op}-,\Lambda-)$  determined by the composition of natural transformations
		\[
		\begin{tikzcd}[ampersand replacement=\&]
			{\mathbb{E}} \& {\mathbb{E}\circ({\lf}^{op} \times \lf)} \&\& {(\mathbb{E}\circ({i}^{op} \times i)) \circ ({\lf}^{op} \times \lf).}
			\arrow["{\nu^2}", Rightarrow, from=1-1, to=1-2]
			\arrow["{\nu^1 \blacktriangle ({\lf}^{op} \times \lf)}", Rightarrow, from=1-2, to=1-4]
		\end{tikzcd}
		\]
	
	\textbf{(5b) The functor $\Lambda$ constructed in (5a) induces the extriangulated functor $j_*:=\overline{\Lambda}:\D/\E\longrightarrow\D$.}
	
	We first show that $\Lambda(f)$ is an isomorphism for every $f\in S_{\E}$. In fact, it is enough to verify that $\Lambda(f)$ is an isomorphism for every $f\in\mathcal{L}\cup\mathcal{R}$.
	
	Let $f:X\to Y$ be a morphism in $\mathcal{L}.$ By definition, there exists a conflation
	\[
	\begin{tikzcd}[ampersand replacement=\&]
		X \& Y \& Z \& {}
		\arrow["f", from=1-1, to=1-2]
		\arrow["g", from=1-2, to=1-3]
		\arrow["\delta", dotted, from=1-3, to=1-4]
	\end{tikzcd}
	\]
	with $Z\in\E$. By Proposition \ref{pro:torsion_funtor_extriangulado}, we obtain the following commutative diagram of conflations
	\begin{equation}\label{dia: f es iso para Phi}
		\begin{tikzcd}[ampersand replacement=\&]
			{i_{*}i^{!}(X)} \& {i_{*}i^{!}(Y)} \& {i_{*}i^{!}(Z)} \& {} \\
			X \& Y \& Z \& {} \\
			{\lf(X)} \& {\lf(Y)} \& {\lf(Z)} \& {} \\
			{} \& {} \& {}
			\arrow["{{\tf(f)}}", from=1-1, to=1-2]
			\arrow["{\varepsilon_2(X)}"', from=1-1, to=2-1]
			\arrow["{{\tf(g)}}", from=1-2, to=1-3]
			\arrow["{\varepsilon_2(Y)}"', from=1-2, to=2-2]
			\arrow["{{{{\nu^2_{C,A}(\delta)}}}}", dotted, from=1-3, to=1-4]
			\arrow["{\varepsilon_2(Z)}", from=1-3, to=2-3]
			\arrow["f", from=2-1, to=2-2]
			\arrow["{{p_{X}}}"', from=2-1, to=3-1]
			\arrow["g", from=2-2, to=2-3]
			\arrow["{p_{Y}}"', from=2-2, to=3-2]
			\arrow["\delta"', dotted, from=2-3, to=2-4]
			\arrow["{{p_{Z}}}", from=2-3, to=3-3]
			\arrow["{{{{\lf(f)}}}}"', from=3-1, to=3-2]
			\arrow["{{{{\delta'_{X}}}}}", dotted, from=3-1, to=4-1]
			\arrow["{{{{\lf(g)}}}}"', from=3-2, to=3-3]
			\arrow["{{{{\delta'_Y}}}}", dotted, from=3-2, to=4-2]
			\arrow["{{{{\nu^3_{C,A}(\delta)}}}}"', dotted, from=3-3, to=3-4]
			\arrow["{{{\delta'_{Z}}}}", dotted, from=3-3, to=4-3]
		\end{tikzcd}
	\end{equation}
	
	Since $Z\in\E=\operatorname{Im}(i_{*})$ and $(\operatorname{Im}(i_{*}),\operatorname{Ker}(i^{!}))$ is an s-torsion pair, Lemma \ref{lem:t(F)=0}implies that $\lf(Z)=0.$ Then, by  Lemma \ref{con-0-iso} we get that $\lf(f)$ is an isomorphism. Similarly, if $f\in\mathcal{R}$, it can be shown that $\tf(f)$ is also an isomorphism. Therefore $\Lambda(f)$ is an isomorphism for every $f\in S_{\E}.$ Therefore,  by the universal property of localization given in Subsection \ref{subsec: propiedad universal locali}, there exists a unique extriangulated functor $(\overline{\Lambda},\psi^3):(\D/\E,\Gg,\uu)\longrightarrow(\D,\Ee,\s)$ such that the following diagram commutes
	\begin{equation}\label{diag:  Lambda pasa a j! de la localizacion}
		\begin{tikzcd}[ampersand replacement=\&]
			{(\D,\Ee,\s)} \&\& {(\D,\Ee,\s)} \\
			\\
			{(\D/\E,\Gg,\uu)}
			\arrow["{(\Lambda,\overline{\psi^3})}", from=1-1, to=1-3]
			\arrow["{(Q,\psi^1)}"', from=1-1, to=3-1]
			\arrow["{(\overline{\Lambda},\psi^3)}"', from=3-1, to=1-3]
		\end{tikzcd}
	\end{equation}
	
	More explicity the functor  $j_*:=\overline{\Lambda}$ and the natural transformation $\psi^3$ are defined as follows:
	
	\begin{enumerate}
		\item[a')] Let $\alpha\in\Hom_{\D/\E}(X,Y)$ be represented by a diagram of the form
		\[
		\begin{tikzcd}[ampersand replacement=\&]
			\& Z \& \\
			X \&\& Y
			\arrow["s"', from=1-2, to=2-1]
			\arrow["g", from=1-2, to=2-3]
		\end{tikzcd}
		\]
		where $s\in S_{\E}$. Then $\alpha=Q(g)\circ Q(s)^{-1}$ and consequently
		\[
		\overline{\Lambda}(\alpha):=\Lambda(g)\Lambda(s)^{-1}.
		\]
		
		\item[b')] Let $[\overline{t}\di\overline{\delta}/\overline{s}]
		=
		\left[
		\begin{tikzcd}[ampersand replacement=\&]
			C \& Z \& X \& A
			\arrow["{\overline{t}}"', from=1-2, to=1-1]
			\arrow["{\overline{\delta}}", dotted, from=1-2, to=1-3]
			\arrow["{\overline{s}}"', from=1-4, to=1-3]
		\end{tikzcd}
		\right]$ be an element of $\Gg(C,A)$. We define $$\psi^{3}([\overline{t}\di\overline{\delta}/\overline{s}])
		:=
		\Lambda(s)^{-1}\cdot\overline{\psi^{3}}(\delta)\cdot\Lambda(t)^{-1}.$$
	\end{enumerate}
	In particular $\im(j_*)=\im(\Lambda)=\Ker(i^!).$

	\textbf{(5c) $j_*:=\overline{\Lambda}$ is a right adjoint of $j^*:=Q:\D\to \D/\E$.}\\
	We show that there exists a functorial isomorphism 
	$$\Theta_{X,Y}:\Ho_{\D/\E}(Q(X),Y)\xrightarrow{\sim}\Ho_{\D}(X,\overline{\Lambda}(Y))$$
	for every $Y\in\D/\E$ and $X\in\D$. By construction, $\overline{\Lambda}(Y):=\lf(Y)=Y'_{0}$ and $Q(X)=X$. Thus, it suffices to prove that there is a functorial isomorphism  $\Ho_{\D/\E}(X,Y)\xrightarrow{\sim}\Ho_{\D}(X,Y'_{0}).$ To this end, we define an isomorphism of groups
	\[
	\Theta_{X,Y}:
	\Hom_{\D/\E}(X,Y)
	\to
	\Hom_{\D}(X,Y'_{0}).
	\]
	
	Indeed, let $\alpha:X\to Y$ be a morphism in $\D/\E$ represented by the following diagram
	\[
	\begin{tikzcd}[ampersand replacement=\&]
		\& Z \& \\
		X \&\& Y.
		\arrow["f", from=2-1, to=1-2]
		\arrow["s"', from=2-3, to=1-2]
		\arrow[
		"\alpha"',
		decorate,
		decoration={snake, amplitude=.4mm, segment length=2mm},
		from=2-1,
		to=2-3
		]
	\end{tikzcd}
	\]
	that is, $\alpha=Q(s)^{-1}Q(f).$ By \eqref{diag: sucesion de torsion 2-1}, we have  the conflation in $\D$
	\begin{equation}\label{diag: sucesion de X y eta 1}
		\begin{tikzcd}
			i_{*}i^!Y & Y & Y'_{0} & {}
			\arrow["{\varepsilon_2(Y)}",from=1-1, to=1-2]
			\arrow["{v_{Y}}", from=1-2, to=1-3]
			\arrow["{\delta'_{Y}}", dotted, from=1-3, to=1-4]
		\end{tikzcd}
	\end{equation}
	such that $i^{!}(Y'_{0})=0$ and $v_Y\in S(\mathcal E)$. Since $Y'_{0}\in\operatorname{Ker}(i^!)=\E^\perp$, Lemma \ref{lem: para supra de adjunto derecho} implies that
	$
	\Hom_{\D}(s,Y'_{0}):
	\Hom_{\D}(Z,Y'_{0})
	\xrightarrow{\sim}
	\Hom_{\D}(Y,Y'_{0}).$
	Hence, for $v_Y:Y\to Y'_{0}$, there is a unique $m_z:Z\to Y'_{0}$ such that 
	$m_z\circ s=v_Y$. We therefore define
	\[
	\Theta_{X,Y}(\alpha):=m_z\circ f.
	\]
	Let us show that $\Theta_{X,Y}(\alpha)$ is well defined. Inded, let 
	$$Q(s_2)^{-1}Q(f_2)=\alpha=Q(s_1)^{-1}Q(f_1),$$ 
	that is, $\alpha$ is represented by 
	$X\xrightarrow{f_1}Z_1\xleftarrow{s_1}Y$ and  $X\xrightarrow{f_2}Z_2\xleftarrow{s_2}Y.$ Let $m_i:Z_i\to Y'_0$ be such that $v_Y=m_i\circ s_i$ for $i=1,2.$ We need to show that $m_2\circ f_2=m_1\circ f_1.$ Using that $m_1\,s_1=m_2\,s_2$ and the two equalites above defining $\alpha,$ it follows that $Q(m_2f_2)=Q(m_1f_1)$ and thus $Q(h)=0$ for 
	$h:=m_2f_2-m_1f_1.$ But, from $Q(h)=0,$ it follows that there is some $s:W\to X$ in $\mathcal{S}_\E$ such that $hs=0.$ Now, by Lemma \ref{lem: para supra de adjunto derecho} we know that $\Hom_\D(s,Y'_0):\Hom_\D(X,Y'_0)\xrightarrow{\sim}\Hom_\D(W,Y'_0).$ Hence the equality $\Hom_\D(s,Y'_0)(h)=hs=0$ implies that $h=0$ and thus $\Theta_{X,Y}(\alpha)$ is well defined.
	
	Notice that, depending on the chosen convention for the localization, morphisms in $\D/\E$ can be represented by either left or right fractions. The representation above uses the convention opposite to the one used previously. The two descriptions are equivalent, and we use the present representation for the construction of $\Theta_{X,Y}$.
	
	\textbf{(5d) $\Theta_{X,Y}:\Hom_{\D/\E}(X,Y)\to \Hom_{\D}(X,Y'_{0})$ is surjective.}
	
	Let $g\in\operatorname{Hom}_{\D}(X,Y'_{0})$. Consider the morphism $\rho:X\to Y$ represented by the following diagram
	\[
	\begin{tikzcd}
		& Y'_{0} & \\
		X
		\arrow[ru,"g"]
		\arrow[
		"\rho"',
		decorate,
		decoration={snake, amplitude=.4mm, segment length=2mm},
		from=2-1,
		to=2-3
		]
		& & Y
		\arrow[lu,"v_{Y}"']
	\end{tikzcd}
	\]
	with $v_{Y}\in S_{\E}$. Thus $\rho=Q(v_{Y})^{-1}Q(g)=[g/v_{Y}]$. By definition of $\Theta_{X,Y}$, we have $\Theta_{X,Y}(\rho)=1_{Y'_{0}}\circ g=g$. Hence $\Theta_{X,Y}$ is surjective.
	
	\textbf{(5e) $\Theta_{X,Y}:\Hom_{\D/\E}(X,Y)\to\Hom_{\D}(X,Y'_{0})$ is a morphism of abelian groups.}
	
	Indeed, let $\alpha_{1},\alpha_{2}$ in $\operatorname{Hom}_{\D/\E}(X,Y)$ be represented by the diagrams
	\[
	\begin{tikzcd}
		& Z &&&& Z & \\
		X && Y && X && Y.
		\arrow["{f_{1}}", from=2-1, to=1-2]
		\arrow[
		"{\alpha_{1}}"',
		decorate,
		decoration={snake, amplitude=.4mm, segment length=2mm},
		from=2-1,
		to=2-3
		]
		\arrow["{s_{1}}"', from=2-3, to=1-2]
		\arrow["{f_{2}}", from=2-5, to=1-6]
		\arrow[
		"{\alpha_{2}}"',
		decorate,
		decoration={snake, amplitude=.4mm, segment length=2mm},
		from=2-5,
		to=2-7]
		\arrow["{s_{2}}"', from=2-7, to=1-6]
	\end{tikzcd}
	\]
	where $s_{1},s_{2}\in S(\E)$. Since $S(\E)$ is a multiplicative system, there exist morphisms $t_{1},t_{2}:Z\longrightarrow W$ in $S(\E)$ such that $t_{1}s_{2}=t_{2}s_{1}$. Consequently, we have the following diagrams
	\[
	\begin{tikzcd}
		&& W &&&&& W && \\
		& Z && W &&& Z && W \\
		X &&&& Y & X &&&& Y
		\arrow["{t_{2}}", from=2-2, to=1-3]
		\arrow["{1_{W}}"', from=2-4, to=1-3]
		\arrow["{t_{1}}", from=2-7, to=1-8]
		\arrow["{1_{W}}"', from=2-9, to=1-8]
		\arrow["{f_{1}}", from=3-1, to=2-2]
		\arrow["{t_{2}f_{1}}"', from=3-1, to=2-4]
		\arrow["{s_{1}}"{pos=0.3}, from=3-5, to=2-2]
		\arrow["{t_{1}s_{2}}"', from=3-5, to=2-4]
		\arrow["{f_{2}}", from=3-6, to=2-7]
		\arrow["{t_{1}f_{2}}"', from=3-6, to=2-9]
		\arrow["{s_{2}}"{pos=0.3}, from=3-10, to=2-7]
		\arrow["{t_{2}s_{1}}"', from=3-10, to=2-9]
	\end{tikzcd}
	\]
	Thus, $\alpha_{1}$ and $\alpha_{2}$ can be represented using a common denominator. More precisely,
	$\alpha_{1}=[f_{1}/s_{1}]=[t_{2}f_{1}/t_{1}s_{2}]$ and
	$\alpha_{2}=[f_{2}/s_{2}]=[t_{1}f_{2}/t_{2}s_{1}]$.
	
	Therefore, without loss of generality, by  \cite[Prop. 3.19]{NakaokaOgawaSakai2022} we may assume that $\alpha_{1}=[f_{1}/s]$ and $\alpha_{2}=[f_{2}/s].$ Then
	$\alpha_{1}+\alpha_{2}=[(f_{1}+f_{2})/s]$ and hence
	\[
	\Theta_{X,Y}(\alpha_{1}+\alpha_{2})
	=m\circ(f_{1}+f_{2})
	=m\circ f_{1}+m\circ f_{2},
	\]
	where $m:Z\to Y'_{0}$ is the unique morphism such that $m\circ s=v_{Y}.$
	On the other hand,
	$\Theta_{X,Y}(\alpha_{1})=m\circ f_{1}$ and
	$\Theta_{X,Y}(\alpha_{2})=m\circ f_{2}$. Consequently,
	 $\Theta_{X,Y}$ is a morphism of abelian groups.
	
	\textbf{(5f) $\Theta_{X,Y}:\Hom_{\D/\E}(X,Y)\to\Hom_{\D}(X,Y'_{0})$ is natural in both variables.}\\ 
	The proof is similar to the one given for $\Gamma_{X,Y}:\Hom_\D(X_0,Y)\to \Hom_{\D/\E}(X,Y)$ in  (3c) and (3d).
	
	\textbf{(5g) $j_*:=\overline{\Lambda}$ is faithful and full.}\\
	Since $(j_!,j^*,j_*)$ is an adjoint triple and $j_!$ is full and faitfull, it follows from Lemma \ref{lem: adjunto fiel y pleno} that $j_*$ is full and faithful.
	
		\textbf{(5h) The transformations $\eta_{4}: 1_{\D}\to j_*\,j^*$ and
		$\varepsilon_{4}: j^*\, j_* \to 1_{\D/ \E}$ are natural extriangulated transformations.} 
	Indeed, by \cite[Theorem 3.1.5]{Borceux1},
	the morphisms $(\eta_{4})_{X}:=\Theta_{X,j^{*}(X)}
	(1_{j^*(X)})$ and
	$\varepsilon_{Y}:= \Theta_{j_{*}Y,Y}^{-1}(1_{j_{*}(Y)})=\Theta_{j_{*}Y,Y}^{-1}(1_{Y_{0}'})$ define natural transformations. In particular, by construction,
	$(\eta_{4})_{X}=v_{X}$ and
	$(\varepsilon_4)_{Y}=[ 1_{Y_{0}'}/ v_{y}]=Q(v_{Y})^{-1}$.
	
	Since $(j^*:=Q, \psi^1)$ and
	$(j_*:=\overline{\Lambda}, \psi^3)$ are extriangulated functors,
	their composition $(Q\, \overline{\Lambda}, \lambda^3)$ is also an
	extriangulated functor, where, for each $C,A \in \D/\E$,
	the morphism
	$\lambda^3_{C,A}: \mathbb{G}(C,A) \to
	\mathbb{G}(\overline{\Lambda} \, Q C, \, \overline{\Lambda} \, Q A)$
	is defined as follows.

	\begin{equation}\label{eq: lamb 3}
		\begin{aligned}
			\lambda^3_{C,A}([\overline{t}\di\overline{\delta}/\overline{s}])
			&= \psi^1_{j_{*}C,j_{*}A}
			\left(\psi^3_{C,A}([\overline{t}\di\overline{\delta}/\overline{s}])\right)\\
			&= \psi^1_{C_{0}', A_{0}'}
			\left(\psi^3([\overline{t}\di\overline{\delta}/\overline{s}])\right)\\
			&= \psi^1_{C_{0}', A_{0}'}
			\left(\Lambda(s)^{-1}\cdot
			\overline{\psi}^{3}(\delta)\cdot \Lambda(t)^{-1} \right)\\
			& = [\overline{1_{C'_0}}\di \Lambda(s)^{-1}\cdot
			\overline{\psi^{3}}(\delta)\cdot\Lambda(t)^{-1} / \overline{1_{A'_0}}]\\
			& = [\overline{1_{C'_0}}\di (i\, \lf(s))^{-1}\cdot
			\overline{\psi^{3}}(\delta)\cdot (i\, \lf(t))^{-1} / \overline{1_{A'_0}}]\\
			& = [\overline{1_{C'_0}}\di \lf(s)^{-1}\cdot
			\overline{\psi^{3}}(\delta)\cdot  \lf(t)^{-1} / \overline{1_{A'_0}}] & (see \  Step  \ 5a)\\
			& = [\overline{1_{C'_0}}\di \lf(s)^{-1}\cdot
			\nu^1_{\lf Z, \lf X} \, \nu^2_{Z,X}(\delta)\cdot  \lf(t)^{-1} / \overline{1_{A'_0}}] & ( \nu^1_{\lf Z, \lf X} \  inclusion)\\
			& = [\overline{1_{C'_0}}\di \lf(s)^{-1}\cdot \nu^2_{Z,X}(\delta)\cdot  \lf(t)^{-1} / \overline{1_{A'_0}}]
		\end{aligned}
	\end{equation}
	
	for every element $[\overline{t}\di\overline{\delta}/\overline{s}]$
	of $\mathbb{G}(C,A)$. Similarly, $(1_{\D/ \E}, id):
	\D/\E \to \D/\E$ is an extriangulated functor. We now show that
	$\varepsilon_{4}:  (j^*\,j_*, \lambda^3) \to (1_{\D / \E}, id)$ is a natural extriangulated transformation in the sense of
	\cite[Definition 2.11]{NakaokaOgawaSakai2022}.\\
	\textit{Claim:} The following diagram commutes.
	\begin{equation}\label{dia: epsilon4-nat}
	\begin{tikzcd}[ampersand replacement=\&]
		{\mathbb{G}(C,A)} \&\& {\mathbb{G}(Q \overline{\Lambda}C, Q \overline{\Lambda}A)} \\
		{\mathbb{G}(C,A)} \&\& {\mathbb{G}(Q \overline{\Lambda}C, A)}
		\arrow["{\lambda^3_{C,A}}", from=1-1, to=1-3]
		\arrow["{id_{C,A}}"', from=1-1, to=2-1]
		\arrow["{\mathbb{G}(Q \overline{\Lambda}C, \varepsilon_{A})}", from=1-3, to=2-3]
		\arrow["{\mathbb{G}(\varepsilon_{C},A)}"', from=2-1, to=2-3]
	\end{tikzcd}
	\end{equation}

	Indeed, let $[\overline{t}\di\overline{\delta}/\overline{s}]$ be an element of $\mathbb{G}(C,A)$. On the one hand,
	
	\begin{align*}
		\mathbb{G}(\varepsilon_{C},A) \, id_{C,A} ([\overline{t}\di\overline{\delta}/\overline{s}])& = \mathbb{G}(\varepsilon_{C},A)([\overline{t}\di\overline{\delta}/\overline{s}])\\
		& =[\overline{t}\di\overline{\delta}/\overline{s}] \cdot \varepsilon_{C}\\
		& = [\overline{t}\di\overline{\delta}/\overline{s}] \cdot Q(v_{C})^{-1}\\
		& = [\overline{t}\di\overline{\delta}/\overline{s}] \cdot [\overline{1_{C'_{0}}}/ \overline{v_{C}}].
	\end{align*}

	The last action is computed according to \cite[Definition 3.23]{NakaokaOgawaSakai2022}, using the left fraction $[\overline{v_{C}}/ \overline{1_{C}}]$, which is equivalent to the right  fraction $[\overline{1_{C'_{0}}} / \overline{v_{C}}]$ via the following diagram.

\[\begin{tikzcd}[ampersand replacement=\&]
	C \& C \\
	{C'_{0}} \& {C'_{0}}
	\arrow["{1_{C}}", from=1-1, to=1-2]
	\arrow["{v_{C}}"', from=1-1, to=2-1]
	\arrow["{v_{C}}", from=1-2, to=2-2]
	\arrow["{1_{C'_{0}}}"', from=2-1, to=2-2]
\end{tikzcd}\]
	
	obtained from condition $(M2)$ of $\mathcal{S}_{\E}$. To compute $[\overline{t}\di\overline{\delta}/\overline{s}] \cdot [\overline{v_{C}}/ \overline{1_{C}}]$, it suffices to consider the following diagram 
\[
\begin{tikzcd}[ampersand replacement=\&]
	\&\&\&\& {} \\
	A \& X \& Y \& Z \& Z \\
	\&\&\& C \& C
	\arrow["s", from=2-1, to=2-2]
	\arrow["a", from=2-2, to=2-3]
	\arrow["b", from=2-3, to=2-4]
	\arrow["\delta"{pos=0.7}, dashed, shorten <=2pt, shorten >=15pt, from=2-4, to=1-5]
	\arrow["t"', from=2-4, to=3-4]
	\arrow["{1_Z}"', from=2-5, to=2-4]
	\arrow["t", from=2-5, to=3-5]
	\arrow["{1_C}", from=3-5, to=3-4]
\end{tikzcd}
\]
	
	and so,
	$$\begin{aligned}
		[\overline{t}\di\overline{\delta}/\overline{s}] \cdot [\overline{1_{C'_{0}}}/ \overline{v_{C}}] & =  [\overline{t}\di\overline{\delta}/\overline{s}] \cdot [\overline{v_{C}}/ \overline{1_{C}}] \\
		& = [\overline{ v_{C}\,t} \di\overline{1_{Z} \cdot \delta}/\overline{s}] \\
		& = [\overline{ v_{C}\,t} \di\overline{\delta}/ \overline{s}].
	\end{aligned}$$

	On the other hand, 
	\begin{equation}\label{eq: eps-extr-lad3}
		\begin{aligned}
			\mathbb{G}(Q\, \overline{\Lambda} C, \varepsilon_{A})\, \lambda^3_{C,A}([\overline{t}\di\overline{\delta}/\overline{s}]) & = \mathbb{G}(Q\, \overline{\Lambda} C, \varepsilon_{A}) \big( [\overline{1_{C'_0}}\di \lf(s)^{-1}\cdot \nu^2_{Z,X}(\delta)\cdot  \lf(t)^{-1} / \overline{1_{A'_0}}]\big)\\
			& = [\overline{1_{C'_0}}\di \lf(s)^{-1}\cdot \nu^2_{Z,X}(\delta)\cdot  \lf(t)^{-1} / \overline{1_{A'_0}}] \cdot  [\overline{1_{A'_0}}/ \overline{v_A}]. 
		\end{aligned}
	\end{equation} 
	
	We now consider the following canonical diagram of conflations (see Proposition \ref{pro:torsion_funtor_extriangulado})
	\begin{equation}\label{nat-ext-eps}
		\begin{tikzcd}[ampersand replacement=\&]
			\&\&\&\& {} \& \\
			{i_{*}i^{!}A} \& {i_{*}i^{!}X} \& {i_{*}i^{!}Y} \& {i_{*}i^{!}Z} \& {} \& {i_{*}i^{!}C} \\
			A \& X \& Y \& Z \& {} \& C \\
			{A'_{0}} \& {X_{0}'} \& {Y_{0}'} \& {Z_{0}'} \&\& {C_{0}'} \\
			{} \& {} \& {} \& {} \&\& {}
			\arrow["{i_{*}i^{!}s}", from=2-1, to=2-2]
			\arrow["{(\varepsilon_{2})_{A}}"', from=2-1, to=3-1]
			\arrow["{i_{*}i^{!}a}", from=2-2, to=2-3]
			\arrow["{(\varepsilon_{2})_{X}}"', from=2-2, to=3-2]
			\arrow["{i_{*}i^{!}b}", from=2-3, to=2-4]
			\arrow["{(\varepsilon_{2})_{B}}"', from=2-3, to=3-3]
			\arrow["{\nu^3_{Z,X}}"{pos=0.3}, dashed, shorten <=2pt, shorten >=2pt, from=2-4, to=1-5]
			\arrow["{i_{*}i^{!}c}", from=2-4, to=2-6]
			\arrow["{(\varepsilon_{2})_{Z}}"', from=2-4, to=3-4]
			\arrow["{(\varepsilon_{2})_{C}}", from=2-6, to=3-6]
			\arrow["s", from=3-1, to=3-2]
			\arrow["{v_{A}}"', from=3-1, to=4-1]
			\arrow["a", from=3-2, to=3-3]
			\arrow["{v_{X}}"', from=3-2, to=4-2]
			\arrow["b", from=3-3, to=3-4]
			\arrow["{v_{Y}}"', from=3-3, to=4-3]
			\arrow["\delta"{pos=0.3}, dashed, shorten <=2pt, shorten >=10pt, from=3-4, to=2-5]
			\arrow["c", from=3-4, to=3-6]
			\arrow["{v_{Z}}"', from=3-4, to=4-4]
			\arrow["{v_{C}}", from=3-6, to=4-6]
			\arrow["{\lf(s)}"', from=4-1, to=4-2]
			\arrow["{\delta'_{A}}"', dashed, from=4-1, to=5-1]
			\arrow["{\lf(a)}"', from=4-2, to=4-3]
			\arrow["{\delta'_{X}}", dashed, from=4-2, to=5-2]
			\arrow["{\lf(b)}"', from=4-3, to=4-4]
			\arrow["{\delta'_{Y}}", dashed, from=4-3, to=5-3]
			\arrow["{\nu^2_{\lf Z, \lf X}}"'{pos=0.3}, dashed, shorten <=2pt, shorten >=10pt, from=4-4, to=3-5]
			\arrow["{\lf(t)}"', from=4-4, to=4-6]
			\arrow["{\delta'_{Z}}", dashed, from=4-4, to=5-4]
			\arrow["{\delta'_{C}}", dashed, from=4-6, to=5-6]
		\end{tikzcd}
	\end{equation}
	
	Since $t,s \in \mathcal{S}_{\E}$, it follows that
	$\overline{\Lambda}(t)$ and $\overline{\Lambda}(s)$ are isomorphisms in
	$\D/\E$. By the definition of $\overline{\Lambda}$, we conclude that
	$\lf(s): A'_{0}\to X'_{0}$ and
	$\lf(t): Z'_{0}\to C'_{0}$ are isomorphisms.
	Hence, by \cite[Proposition 3.7]{NP19}, the element
	$\lf(s)^{-1} \cdot \nu^2_{Z,X}(\delta) \cdot \lf(t)^{-1}$
	is realized by the following conflation
	\[\begin{tikzcd}[ampersand replacement=\&]
		{A'_{0}} \&\& {Y'_{0}} \&\& {C'_{0}} \& {}
		\arrow["{\lf(a)\lf(s)}", from=1-1, to=1-3]
		\arrow["{\lf(t)\lf(b)}", from=1-3, to=1-5]
		\arrow["{\lf(s)^{-1} \cdot \nu^2_{Z,X}(\delta) \cdot \lf(t)^{-1}}"{pos=1}, dashed, from=1-5, to=1-6]
	\end{tikzcd}\]
	With the above in place, we perform the action
	$  [\overline{1_{A'_0}}/ \overline{v_A}] \cdot [\overline{1_{C'_0}}\di \lf(s)^{-1}\cdot \nu^2_{C,A}(\delta)\cdot  \lf(t)^{-1} / \overline{1_{A'_0}}]$,
	according to \cite[Definition 3.23]{NakaokaOgawaSakai2022}, using the left fraction $[\overline{v_{C}} / \overline{1_{C}}]$. This yields the following diagram
	
\[\begin{tikzcd}[ampersand replacement=\&]
	{A'_0} \& {A'_0} \&\&\&\&\& \\
	{A'_0} \& {A'_{0}} \&\& {Y'_{0}} \&\& {C'_{0}} \& {} \\
	\&\&\&\&\& {C'_{0}}
	\arrow["{1_{A'_0}}"', from=1-1, to=2-1]
	\arrow["{1_{A'_0}}"', from=1-2, to=1-1]
	\arrow["{1_{A'_0}}"{pos=0.4}, from=1-2, to=2-2]
	\arrow["{1_{A'_0}}"', from=2-2, to=2-1]
	\arrow["{\lf(a)\lf(s)}", from=2-2, to=2-4]
	\arrow["{\lf(t)\lf(b)}", from=2-4, to=2-6]
	\arrow["{\lf(s)^{-1} \cdot \nu^2_{Z,X}(\delta) \cdot \lf(t)^{-1}}"{pos=1}, dashed, from=2-6, to=2-7]
	\arrow["{1_{C'_0}}"', from=2-6, to=3-6]
\end{tikzcd}\]
	Consequently, we have
	$[\overline{1_{A'_0}}/ \overline{v_A}] \cdot [\overline{1_{C'_0}}\di \lf(s)^{-1}\cdot \nu^2_{C,A}(\delta)\cdot  \lf(t)^{-1} / \overline{1_{A'_0}}]= [\overline{1_{C'_0}}\di \lf(s)^{-1}\cdot \nu^2_{Z,X}(\delta)\cdot  \lf(t)^{-1} / \overline{v_A}]$.
	Therefore, it follows from \eqref{eq: eps-extr-lad3} that
	
	We now show that diagram \eqref{dia: epsilon4-nat} commutes. To this end, it suffices to verify that
	$[\overline{1_{C'_0}}\di \lf(s)^{-1}\cdot \nu^2_{Z,X}(\delta)\cdot  \lf(t)^{-1} / \overline{v_A}]= [\overline{v_{C} t}\di \delta / \overline{s}]$,
	in the sense of \cite[Proposition 3.9]{NakaokaOgawaSakai2022}. Indeed, we have the following diagram
\[\begin{tikzcd}[ampersand replacement=\&]
	\&\& {C'_{0}} \&\& {A'_{0}} \&\& \\
	{C'_{0}} \&\& Z \&\& {A'_{0}} \&\& A \\
	\&\& Z \&\& X
	\arrow["{\lf(s)^{-1}\cdot \nu^2_{Z,X}(\delta)\cdot  \lf(t)^{-1}}", dashed, from=1-3, to=1-5]
	\arrow["{1_{C'_0}}"', from=1-3, to=2-1]
	\arrow["{1_{A'_{0}}}", from=1-5, to=2-5]
	\arrow["{\lf(t)v_Z}"', from=2-3, to=1-3]
	\arrow["{v_Ct}"', from=2-3, to=2-1]
	\arrow[dashed, from=2-3, to=2-5]
	\arrow["{1_{Z}}"{pos=0.2}, from=2-3, to=3-3]
	\arrow["{v_A}"', from=2-7, to=1-5]
	\arrow["{v_{A}}"', from=2-7, to=2-5]
		\arrow["\xi", dashed, from=2-3, to=2-5]
	\arrow["s", from=2-7, to=3-5]
	\arrow["{v_{C}t}"{pos=0.4}, from=3-3, to=2-1]
	\arrow["\delta"', dashed, from=3-3, to=3-5]
	\arrow["{\lf(s)^{-1}v_{X}}"{pos=0.7}, from=3-5, to=2-5]
\end{tikzcd}\]
	 
	from which it follows that
	$ \xi= (\lf(s)^{-1}v_X) \cdot \delta \cdot 1_{Z}= \lf(s)^{-1} \cdot (v_X \cdot \delta )$.
	Now, using diagram \eqref{nat-ext-eps}, we obtain
	
	\begin{align*}
		1_{A'_0} \cdot \big( \lf(s)^{-1}\cdot \nu^2_{\lf Z, \lf X}(\delta)\cdot  \lf(t)^{-1} \big) \cdot (\lf(t)\, v_{Z}) & = \big( (\lf(s)^{-1}\cdot \nu^2_{\lf Z, \lf X}(\delta)\cdot  \lf(t)^{-1}) \cdot \lf(t) \big)   \cdot  v_{Z}\\
		 & = \big( (\lf(s)^{-1}\cdot \nu^2_{\lf Z, \lf X}(\delta)) \cdot  (\lf(t)^{-1} \lf(t)) \big)   \cdot  v_{Z}\\
		 & = (\lf(s)^{-1}\cdot\nu^2_{\lf Z, \lf X}(\delta)  \big)   \cdot  v_{Z}\\
		 & = \lf(s)^{-1}\cdot  (\nu^2_{\lf Z, \lf X}(\delta)    \cdot  v_{Z})\\
		 & =  \lf(s)^{-1}\cdot (v_X \cdot \delta)
	\end{align*}
	
	This proves that \eqref{dia: eta-nat} commutes. Therefore, $\varepsilon_{4} : (j^* j_*, \lambda^3)  \to (Id_{\D/ \E}, id)$ is a natural extriangulated transformation.\\
	
	We now show that $\eta_{4}\colon(\operatorname{Id}_{\D},\operatorname{id}) \to (j_* j^*, \lambda^{4})$ is an extriangulated natural transformation	where $\lambda^{4}\colon \Ee(C,A)\to
	\Ee(j_* j^* C, j_* j^*  A)$ is defined as the
	composition
	$\psi^{3}\circ\bigl(\psi^{1}\blacktriangle(Q^{op}\times Q)\bigr)$, that is,
	
	\begin{align*}
		\lambda_{Z,X}^4(\delta)& = \psi^3_{QZ,QX} \psi^1_{Z,X}(\delta)\\
		& = \psi^3_{QZ,QX}([\overline{1_{Z}}\di \overline{\delta} / \overline{1_{X}}] )\\
		& = \overline{\Lambda}(1_{X})^{-1}\cdot \overline{\psi}_{Z,X}(\delta) \cdot \overline{\Lambda}(1_{Z})^{-1}\\
		& =\overline{\psi}_{Z,X}(\delta) & (\text{see Step $(6a)$ } )\\
		& =\nu^1_{\lf Z, \lf X} \nu^2_{Z,X}(\delta) & (\nu^1 \text{ is inclusion})\\
		& = \nu^2_{ Z,  X}(\delta).	
	\end{align*}
	for every $ \delta \in \Ee(Z,X)$.\\
	
	\textit{Claim:} The following diagram commutes
	
	\begin{equation}\label{dia: eps-extr2}
		\begin{tikzcd}[ampersand replacement=\&]
			{\Ee(Z,X)} \&\& {\Ee(Z,X)} \\
			{\Ee(j_{*}j^*Z, j_{*}j^*X)} \&\& {\Ee(Z, j_{*}j^{*}X)}
			\arrow["id", from=1-1, to=1-3]
			\arrow["{\lambda^4_{Z,X}}"', from=1-1, to=2-1]
			\arrow["{\Ee(Z, \eta_{X})}", from=1-3, to=2-3]
			\arrow["{\Ee(\eta_{Z}, j_{*}j^{*}X )}"', from=2-1, to=2-3]
		\end{tikzcd}
	\end{equation}

	Indeed, let $\delta\in\Ee(Z,X)$, we have 
	\begin{align*}
		\Ee(Z, \eta_{X})\, id (\delta)& = \Ee(Z, \eta_{X})(\delta)\\
		& =\eta_{X}\cdot \delta\\
		& = \Theta_{X,j^{*}(X)}(1_{j^*X}) \cdot \delta\\
		& = v_{X} \cdot \delta & ( \text{see } \ref{nat-ext-eps})\\
		& =  \nu^2_{ Z,  X}(\delta) \cdot v_{Z}\\
		& =  \nu^2_{ Z,  X}(\delta) \cdot \Theta_{Z,j^{*}(Z)}(1_{j^*Z}) \\
		& =  \nu^2_{ Z,  X}(\delta) \cdot \eta_{Z}\\
		& = \Ee(\eta_{Z}, \lf X) \nu_{ Z,  X}^2(\delta)\\
		& = \Ee(\eta_{Z}, \lf X) \nu_{\lf Z, \lf X}^1 \nu_{ Z,  X}^2(\delta)\\
		& = \Ee(\eta_{Z}, \lf X) \lambda^4_{Z,X} (\delta).
	\end{align*}	
	
	Thus, diagram \eqref{dia: eps-extr2} commutes, proving that
	$\eta_{4}\colon(\operatorname{Id}_{\D},\operatorname{id}) \to (j_* j^*, \lambda^{4})$
	is an extriangulated natural transformation.
\end{proof}

  Theorem~\ref{teo: completar recoll} can be applied when the categories $\D'$ and $\D$ are triangulated, see Example \ref{Eextf=tf}. However, in the triangulated setting, condition~(c) of Theorem~\ref{teo: completar recoll} can be omitted since it follows from the first two conditions, as we show below.

\begin{lemma}\label{Lsttf-tc}
	Let $\mathcal{D}$ and $\mathcal{D}'$ be triangulated categories, and let $i_*:\mathcal{D}'\to\mathcal{D}$ be a fully faithful triangulated functor admitting a triangulated left
	adjoint $i^*$ and a triangulated right adjoint $i^!.$ Then, for  
 $\mathcal{E}:=\operatorname{Im}(i_*),$ the following statements hold true.
\begin{itemize}
\item[(a)] $\Ker(i^*)={}^{\perp}\mathcal{E}$ and $\Ker(i^!)=\mathcal{E}^{\perp}.$
\item[(b)] For each $X\in\mathcal{D}$ there exists a distinguished triangle
	\[
  K_{X} 	\to  X \xrightarrow{\eta_X} i_*i^*X \to\Sigma K_{X},
	\]
	where $K_X\in\Ker(i^*)$ and $\eta:1_\D\to i_* i^*$ is the unite of the adjoint pair $(i^*,i_*).$
\item[(c)] For or every $X\in\mathcal{D}$ there exists a distinguished triangle
	\[
	i_*i^!X\xrightarrow{\varepsilon_X}X\longrightarrow C_X\longrightarrow\Sigma i_*i^!X,
	\]
	where $C_X\in\Ker(i^!)$ and $\varepsilon:i_* i^!\to 1_\D$ is the co-unite of the adjoint pair $(i_*, i^!).$
\end{itemize} 
\end{lemma}

\begin{proof} Notice that (a) follows from Lemma \ref{LPadjPerp}. We only prove (b) since (c) follows similarly.

	Let $X\in\D.$   Complete $\eta_X:X \to i_*i^*X$ to a distinguished triangle
		\[
	K_{X} 	\longrightarrow X \xrightarrow{\eta_x} i_*i^*X \longrightarrow\Sigma K_{X}
	\]
	Applying $i^*$ to the above triangle we get the distinguished triangle
	\[
i^*K_{X} 	\longrightarrow  i^*X \xrightarrow{i^*(\eta_X)}  i^* i_*i^*X \longrightarrow\Sigma i^* K_{X}.
\]
	Because $i_*$ is fully faithful, the co-unit $\varepsilon': i^*i_* \to 1_{\D'}$ is  an isomorphism. Hence,  the triangular identity $\varepsilon'_{i^*X}\circ i^*(\eta_X)= 1_{i^*X}$ implies  that  $i^*(\eta_X)$ is an isomorphism. Thus $i^*K_X=0$ and so $K_X\in\Ker(i^*).$
\end{proof}

\begin{proposition}\label{prop:TTF triple from adjoint triple2}
	Let $\mathcal{D}$ and $\mathcal{D}'$ be triangulated categories viewed as extriangulated categories with negative first extensions as in Example 
	\ref{triang-ext-1Ex}.
	 Let $i_*:\mathcal{D}'\to \mathcal{D}$ be a fully faithful triangulated functor admitting triangulated adjoints $i^*\dashv i_*\dashv i^!$. Then 
	 $\bigl(\Ker(i^*),\operatorname{Im}(i_*),\Ker(i^!)\bigr)$ is an s-TTF triple in $\mathcal{D}.$ 
\end{proposition}

\begin{proof} By Lemma \ref{Lsttf-tc} we get that  $\bigl(\Ker(i^*),\operatorname{Im}(i_*)\bigr)$ and  $\bigl(\operatorname{Im}(i_*),\Ker(i^!)\bigr)$ are torsion pairs in $\D.$ To finish the proof, we need  to show that $\Ee^{-1}(\Ker(i^*),\operatorname{Im}(i_*))=0$ and $\Ee^{-1}(\operatorname{Im}(i_*),\Ker(i^!))=0.$ In order to do that,  it is enough to show that $\Sigma\Ker(i^*)\subseteq\Ker(i^*)$ and $\Sigma\operatorname{Im}(i_*)\subseteq\operatorname{Im}(i_*)$.
	
	Let $X\in\Ker(i^*).$ Since $i^*$ is a triangulated functor, we have $i^*(\Sigma X)\cong\Sigma(i^*X)=\Sigma 0=0$ and thus $\Sigma\Ker(i^*)\subseteq\Ker(i^*)$.
	
	Now let $Z\in\operatorname{Im}(i_*).$ Then $Z\simeq i_*Z'$  for some $Z'\in\mathcal{D}'.$ Since $i_*$ is triangulated, we have $\Sigma Z\simeq\Sigma i_*Z'\simeq i_*(\Sigma Z')$ and thus $\Sigma Z\in \operatorname{Im}(i_*).$
\end{proof}

 The preceding proposition shows that, in the triangulated setting, there is no need to assume separately that $\bigl(\Ker(i^*),\operatorname{Im}(i_*),\Ker(i^!)\bigr)$ is an s-TTF triple. Indeed, this condition follows automatically from the existence of the adjoint triple $i^*\dashv i_*\dashv i^!$. Consequently, as an immediate application of Theorem~\ref{teo: completar recoll}, we obtain the following result.
 
 \begin{corollary}\cite[Thm. 2.1]{Cline-Parshall-Scott}
  Let $\D$ and $D'$ be triangulated categories such that $\D$ satisfies condition (WIC). Consider the diagram of triangulated functors
\begin{equation*}
		\begin{tikzcd}[ampersand replacement=\&, column sep=large]
			{\mathcal{D}'}
			\&\&
			{\mathcal{D}}
			\arrow["{i_{*}}", from=1-1, to=1-3]
			\arrow["{i^{!}}", bend left=35, from=1-3, to=1-1]
			\arrow["{i^{*}}"', bend right=35, from=1-3, to=1-1]
		\end{tikzcd}
	\end{equation*}		
which satisfies the following conditions	
	\begin{enumerate}
		\item[(a)] $(i^{*},i_{*},i^{!})$ is an  adjoint triple;
		\item[(b)] $i_{*}$ is full and faithful.
	\end{enumerate}
	Then $\E:= \im(i_{*})$ is a thick triangulated subcategory of $\D$ and there exist triangulated functors $j^{*}\colon\D\to \D/\E$, $j_{!}\colon\D/\E\to \D$ and $j_{*}\colon\D/\E\to \D$ such that:
	\begin{enumerate}
		\item $(j_{!},j^{*},j_{*})$ is an adjoint triple;
		\item $j_{!}$ and $j_{*}$ are full and faithful;
		\item $\Ker(j^*)=\im(i_*),$ $\im(j_{!})=\Ker(i^*)$ and $\Ker(i^{!})=\im(j_*).$
		\item $(\Ker(i^{*}),\im(i_{*}),\Ker(i^{!}))$ is an s-TTF triple in $\D.$
	\end{enumerate}
	In particular, we obtain the following recollement of triangulated categories
	\begin{equation}
		\begin{tikzcd}[ampersand replacement=\&, column sep=large]
			{\mathcal{D}'}
			\&\&
			{\mathcal{D}}
			\&\&
			{\mathcal{D}/\mathcal{E}}
			\arrow["{i_{*}}", from=1-1, to=1-3]
			\arrow["{i^{!}}", bend left=35, from=1-3, to=1-1]
			\arrow["{i^{*}}"', bend right=35, from=1-3, to=1-1]
			\arrow["{j^{*}}", from=1-3, to=1-5]
			\arrow["{j_{*}}", bend left=35, from=1-5, to=1-3]
			\arrow["{j_{!}}"', bend right=35, from=1-5, to=1-3]
		\end{tikzcd}
	\end{equation}
 \end{corollary}

\newpage
	
	{\footnotesize

		}

	{\footnotesize

			\vskip3mm \noindent {\bf Octavio Mendoza Hernández}:\\ 
		{\it Instituto de Matem\'aticas, Universidad Nacional Aut\'onoma de M\'exico\\
			Circuito Exterior, Ciudad Universitaria
			C.P. 04510, Ciudad de M\'exico, MEXICO.\\ {\tt omendoza@matem.unam.mx}}
		
		\vskip3mm \noindent {\bf Valente Santiago Vargas}:\\ 
		{\it Departamento de Matem\'aticas, Facultad de Ciencias, Universidad Nacional Aut\'onoma de M\'exico\\
			Circuito Exterior, Ciudad Universitaria
			C.P. 04510, Ciudad de M\'exico, MEXICO.\\ {\tt valente.santiago@ciencias.unam.mx}}

		\vskip3mm \noindent  {\bf Edgar Omar Velasco P\'aez}:\\ 
		{\it Facultad de Ciencias, Universidad Nacional Aut\'onoma de M\'exico\\
			Circuito Exterior, Ciudad Universitaria,
			C.P. 04510, Ciudad de M\'exico, MEXICO.\\ {\tt edgar\underline{\,\,}bkz13@ciencias.unam.mx }

		\end{document}